%% file: lcbc.tex
\documentclass[preprint,11pt]{elsarticle}
\usepackage[bookmarks=true,colorlinks=true,linkcolor=blue]{hyperref}

\usepackage{amsmath,amssymb,mathabx}
\usepackage{graphicx,esint,ulem}

\biboptions{sort&compress}
\usepackage[usenames,dvipsnames,svgnames,table]{xcolor}

\usepackage{color}
\definecolor{jwbGreen}{rgb}{0, .6, 0}
\definecolor{jbaPurple}{HTML}{6600FF}
\definecolor{purple}{rgb}{.7, 0., .8}
\definecolor{pinegreen}{rgb}{0.0, 0.47, 0.44}

\newcommand{\blue}{\color{blue}}

\usepackage{calc}
\usepackage[margin=1.in]{geometry}

\usepackage{algorithm,algpseudocode,caption}
\usepackage{algorithmicx}
\algrenewcommand\alglinenumber[1]{\footnotesize #1:} 
\newcommand{\algFontSize}{\footnotesize}

\usepackage{tcolorbox}
\usepackage[english]{babel}
\usepackage{amsmath}
\usepackage{amsfonts}
\usepackage{graphicx}
\usepackage{setspace}
\usepackage{bm}

\usepackage{tabu}

\usepackage{multirow}

\usepackage{mdframed}

\allowdisplaybreaks

\usepackage{tikz}

\newtheorem{theorem}{Theorem}

\newenvironment{proof}[1][Proof]{\begin{trivlist}
\item[\hskip \labelsep {\bfseries #1.}]}{\end{trivlist}}

\input tex/defs.tex
\input tex/trimFig.tex
\usepackage{xargs}
\input tex/plotFigureMacros.tex

\newcommand{\wdh}[1]{{\color{DarkBlue}wdh: #1}}

\title{adegdm}

\begin{document}

\begin{frontmatter}
 \title{Local Compatibility Boundary Conditions for High-Order Accurate Finite-Difference Approximations of PDEs}

\author[rpi]{Nour G.~Al Hassanieh\fnref{NSFgrants}}
\ead{alhasn@rpi.edu}

\author[rpi]{Jeffrey W.~Banks}
\ead{banksj3@rpi.edu}

\author[rpi]{William D.~Henshaw\corref{cor}\fnref{NSFgrants}}
\ead{henshw@rpi.edu}

\author[rpi]{Donald~W.~Schwendeman\fnref{NSFgrants}}
\ead{schwed@rpi.edu}

\address[rpi]{Department of Mathematical Sciences, Rensselaer Polytechnic Institute, Troy, NY 12180, USA}

\cortext[cor]{Corresponding author}

\fntext[NSFgrants]{Research supported by the National Science Foundation under grants DMS-1519934 and DMS-1818926.}

\fntext[RTG]{This work was partially funded by the NSF Research Training Group Grant DMS-1344962.}



\begin{abstract}
We describe a new approach to derive numerical approximations of boundary conditions for high-order accurate finite-difference approximations. The approach, called the Local Compatibility Boundary Condition (LCBC) method, uses boundary conditions and compatibility boundary conditions derived from the governing equations, as well as interior and boundary grid values, to construct a local polynomial, whose degree matches the order of accuracy of the interior scheme, centered at each boundary point.  The local polynomial is then used to derive a discrete formula for each ghost point in terms of the data.  This approach leads to centered approximations that are generally more accurate and stable than one-sided approximations.  Moreover, the stencil approximations are local since they do not couple to neighboring ghost-point values which can occur with traditional compatibility conditions.  The local polynomial is derived using continuous operators and derivatives which enables the automatic construction of stencil approximations at different orders of accuracy.  The LCBC method is developed here for problems governed by second-order partial differential equations, and it is verified for a wide range of sample problems, both time-dependent and time-independent, in two space dimensions and for schemes up to sixth-order accuracy.
\end{abstract}


\end{frontmatter}

\tableofcontents


\setlength{\parskip}{2pt}

\input tex/intro


\input tex/MotivationalProblem


\input tex/ModelEqn


\input tex/LCBCapproach


\input tex/analysis


\input tex/results


\input tex/conclusion

\appendix
\input tex/CornerAppendix.tex 
\input tex/fourthOrderWaveStability.tex
\input tex/timeStepRestrictions

\bibliographystyle{elsart-num}
\bibliography{LCBCbib,henshaw,henshawPapers}
 
\end{document}

%% file: tex/defs.tex
\renewcommand{\url}[1]{}
\newcommand{\citeCount}[1]{}

\newcommand{\bogus}[1]{{}}

\renewcommand{\Re}{{\rm Re}}

\newcommand{\mni}{\medskip\noindent}

\newcommand{\p}{\partial}

\newcommand{\f}[2]{\frac{#1}{#2}}

\def\ba#1\ea{\begin{align}#1\end{align}}

\def\bas#1\eas{\begin{align*}#1\end{align*}}

\def\bat#1\eat{\begin{alignat}{3}#1\end{alignat}}

\def\bats#1\eats{\begin{alignat*}{3}#1\end{alignat*}}

\newcommand{\bse}{\begin{subequations}}
\newcommand{\ese}{\end{subequations}}

\newcommand{\Nx}{N_x}
\newcommand{\Ny}{N_y}
\newcommand{\Nt}{N_t}
\newcommand{\dzy}{D_{0y}}
\newcommand{\dzx}{D_{0x}}
\newcommand{\dpmx}{D_{+x}D_{-x}}
\newcommand{\dpmy}{D_{+y}D_{-y}}
\newcommand{\dpmze}{D_{+\zeta}D_{-\zeta}}
\newcommand{\dzze}{D_{0\zeta}}
\newcommand{\dze}{\Delta\zeta}
\newcommand{\dpmt}{D_{+t}D_{-t}}
\newcommand{\pt}{\partial_t}
\newcommand{\px}{\partial_x}
\newcommand{\py}{\partial_y}
\newcommand{\abs}[1]{\left\vert #1 \right\vert}

\newcommand{\dt}{\Delta t}
\newcommand{\dx}{\Delta x}
\newcommand{\dy}{\Delta y}

\newcommand{\dr}{{\Delta r}}

\newcommand{\defeq}{\overset{{\rm def}}{=}}
\newcommand{\eqdef}{\overset{{\rm def}}{=}}

\newcommand{\bv}{\mathbf{ b}}

\newcommand{\dv}{\mathbf{ d}}

\newcommand{\iv}{\mathbf{ i}}
\newcommand{\jv}{\mathbf{ j}}
\newcommand{\kv}{\mathbf{ k}}

\newcommand{\rv}{\mathbf{ r}}

\newcommand{\vv}{\mathbf{ v}}

\newcommand{\xv}{\mathbf{ x}}

\newcommand{\Gv}{\mathbf{ G}}

\newcommand{\Rv}{\mathbf{ R}}
\newcommand{\Sv}{\mathbf{ S}}

\newcommand{\Uv}{\mathbf{ U}}

\newcommand{\half}{{1\over2}}

\newcommand{\Real}{{\mathbb R}}

\newcommand{\zerov}{\mathbf{0}}

\newcommand{\Ac}{{\mathcal A}}
\newcommand{\Bc}{{\mathcal B}}
\newcommand{\Cc}{{\mathcal C}}
\newcommand{\Dc}{{\mathcal D}}

\newcommand{\Lc}{{\mathcal L}}

\newcommand{\Oc}{{\mathcal O}}

\newcommand{\grad}{\nabla}

\newcommand{\ds}{{\Delta s}}

\newcommand{\Dpx}{D_{+x}}
\newcommand{\Dmx}{D_{-x}}
\newcommand{\Dzx}{D_{0x}}
\newcommand{\Dpy}{D_{+y}}
\newcommand{\Dmy}{D_{-y}}

\newcommand{\Vhat}{\hat{V}}

\newcommand{\kHat}{\hat{k}}
\newcommand{\lambdaHat}{\hat{\lambda}}

\newcommand{\Qtilde}{\tilde{Q}}

\newcommand{\kHatzx}{\kHat_{0x}}
\newcommand{\kHatzy}{\kHat_{0y}}
\newcommand{\alphaMin}{\alpha_{\rm min}}
\newcommand{\betaMax}{\beta_{\rm max}}
\newcommand{\lambdaMax}{\lambda_{\rm max}}

\newcommand{\xt}{\tilde{x}}
\newcommand{\yt}{\tilde{y}}
\newcommand{\ut}{\tilde{u}}

%% file: tex/trimFig.tex
\newlength{\tfwidth}
\newlength{\tfheight}
\newlength{\tfxa}
\newlength{\tfxb}
\newlength{\tfya}
\newlength{\tfyb}
\newcommand{\trimFigWithBox}[6]{%
\setlength\fboxsep{0pt}%
\setlength\fboxrule{1.0pt}
\fbox{\includegraphics[width=#2, clip, trim=#3 #4 #5 #6]{#1}}%
}
\newcommand{\trimFigNoBox}[6]{%
\setlength\fboxsep{1pt}
\setlength\fboxrule{0.0pt}
\fbox{\includegraphics[width=#2, clip, trim=#3 #4 #5 #6]{#1}}%
}
\newcommand{\trimFigHeightWithBox}[6]{%
\setlength\fboxsep{0pt}%
\setlength\fboxrule{1.0pt}
\fbox{\includegraphics[height=#2, clip, trim=#3 #4 #5 #6]{#1}}%
}
\newcommand{\trimFigHeightNoBox}[6]{%
\setlength\fboxsep{1pt}
\setlength\fboxrule{0.0pt}
\fbox{\includegraphics[height=#2, clip, trim=#3 #4 #5 #6]{#1}}%
}

\newsavebox\figBox

\newcommand{\trimw}[6]{%
\sbox\figBox{\includegraphics{#1}}
\setlength{\tfwidth}{\the\wd\figBox}
\setlength{\tfheight}{\the\ht\figBox}
\setlength{\tfxa}{\tfwidth*\real{#3}}%
\setlength{\tfxb}{\tfwidth*\real{#4}}%
\setlength{\tfya}{\tfheight*\real{#5}}%
\setlength{\tfyb}{\tfheight*\real{#6}}%
\trimFigNoBox{#1}{#2}{\tfxa}{\tfya}{\tfxb}{\tfyb}%
}

\newcommand{\trimwb}[6]{%

\sbox\figBox{\includegraphics{#1}}
\setlength{\tfwidth}{\the\wd\figBox}
\setlength{\tfheight}{\the\ht\figBox}
\setlength{\tfxa}{\tfwidth*\real{#3}}%
\setlength{\tfxb}{\tfwidth*\real{#4}}%
\setlength{\tfya}{\tfheight*\real{#5}}%
\setlength{\tfyb}{\tfheight*\real{#6}}%
\trimFigWithBox{#1}{#2}{\tfxa}{\tfya}{\tfxb}{\tfyb}%
}

\newcommand{\trimh}[6]{%
\sbox\figBox{\includegraphics{#1}}
\setlength{\tfwidth}{\the\wd\figBox}
\setlength{\tfheight}{\the\ht\figBox}
\setlength{\tfxa}{\tfwidth*\real{#3}}%
\setlength{\tfxb}{\tfwidth*\real{#4}}%
\setlength{\tfya}{\tfheight*\real{#5}}%
\setlength{\tfyb}{\tfheight*\real{#6}}%
\trimFigHeightNoBox{#1}{#2}{\tfxa}{\tfya}{\tfxb}{\tfyb}%
}

\newcommand{\trimhb}[6]{%

\sbox\figBox{\includegraphics{#1}}
\setlength{\tfwidth}{\the\wd\figBox}
\setlength{\tfheight}{\the\ht\figBox}
\setlength{\tfxa}{\tfwidth*\real{#3}}%
\setlength{\tfxb}{\tfwidth*\real{#4}}%
\setlength{\tfya}{\tfheight*\real{#5}}%
\setlength{\tfyb}{\tfheight*\real{#6}}%
\trimFigHeightWithBox{#1}{#2}{\tfxa}{\tfya}{\tfxb}{\tfyb}%
}

%% file: tex/plotFigureMacros.tex
\newcommandx{\figByHeight}[9][5=0, 6=0, 7=0, 8=0,9=]{
\draw (#1,#2) node[anchor=south west,xshift=-16pt,yshift=-4pt] {\trimh{#3}{#4}{#5}{#6}{#7}{#8}};}
\newcommandx{\figByHeightb}[9][5=0, 6=0, 7=0, 8=0,9=]{
\draw (#1,#2) node[anchor=south west,xshift=-16pt,yshift=-4pt] {\trimhb{#3}{#4}{#5}{#6}{#7}{#8}};}

\newcommandx{\figByHeightWithLabel}[9][5=0, 6=0, 7=0, 8=0,9=]{
\draw (#1,#2) node[anchor=south west,xshift=-16pt,yshift=-4pt] {\trimh{#3}{#4}{#5}{#6}{#7}{#8}} node[draw=white,fill=white,inner sep=1pt,anchor=south west] {#9};}
\newcommandx{\figByHeightWithLabelb}[9][5=0, 6=0, 7=0, 8=0,9=]{
\draw (#1,#2) node[anchor=south west,xshift=-16pt,yshift=-4pt] {\trimhb{#3}{#4}{#5}{#6}{#7}{#8}} node[draw=white,fill=white,inner sep=1pt,anchor=south west] {#9};}

\newcommandx{\figByWidth}[9][5=0, 6=0, 7=0, 8=0,9=]{
\draw (#1,#2) node[anchor=south west,xshift=-16pt,yshift=-4pt] {\trimw{#3}{#4}{#5}{#6}{#7}{#8}};}
\newcommandx{\figByWidthb}[9][5=0, 6=0, 7=0, 8=0,9=]{
\draw (#1,#2) node[anchor=south west,xshift=-16pt,yshift=-4pt] {\trimwb{#3}{#4}{#5}{#6}{#7}{#8}};}

\newcommandx{\figByWidthWithLabel}[9][5=0, 6=0, 7=0, 8=0,9=]{
\draw (#1,#2) node[anchor=south west,xshift=-16pt,yshift=-4pt] {\trimw{#3}{#4}{#5}{#6}{#7}{#8}} node[draw=white,fill=white,inner sep=1pt,anchor=south west] {#9};}
\newcommandx{\figByWidthWithLabelb}[9][5=0, 6=0, 7=0, 8=0,9=]{
\draw (#1,#2) node[anchor=south west,xshift=-16pt,yshift=-4pt] {\trimwb{#3}{#4}{#5}{#6}{#7}{#8}} node[draw=white,fill=white,inner sep=1pt,anchor=south west] {#9};}

%% file: tex/intro.tex
\section{Introduction} \label{sec:introduction}

We describe a new approach for constructing discrete boundary conditions for high-order accurate numerical approximations to partial differential equations (PDEs).  The approach, called the Local Compatibility Boundary Condition (LCBC) method, combines the given physical boundary conditions (BCs) with additional compatibility boundary conditions (CBCs) formed from the PDE and its derivatives.  Our focus here is on finite-difference (and finite-volume) methods for both time-dependent and steady PDEs in second-order form with physical BCs of Dirichlet or Neumann type. 
A high-order accurate centered finite-difference approximation of the spatial operator of the PDE involves a wide stencil which then requires some special treatment to handle the approximation at grid points near the boundary.  Unlike a typical approach involving one-sided approximations of the PDE near the boundary and one-sided approximations of Neumann-type BCs, the LCBC approach results in {\it fully centered} approximations.  These centered approximations are generally more accurate than one-sided approximations, and for the case of time-dependent PDEs they are more stable and less stiff (i.e.~do not decrease the stable explicit time-step).  Furthermore, the new LCBC approach improves upon a more traditional derivation of discrete CBCs by defining local conditions that are not coupled to neighboring grids along the boundary in tangential directions.  As a result, there is no need to solve a system of equations along the boundary which is a significant advantage for explicit time-stepping schemes.  In the case of implicit time-stepping methods, and for approximations of steady (elliptic) PDEs, where the solution of large linear systems is required, this tangential decoupling can also be useful for iterative schemes, such as multigrid and Krylov methods.

The development of LCBCs is motivated by our interest in high-order accurate approximations of PDEs in complex domains using overset grids, although the applicability of LCBCs is broader.  As shown in Figure~\ref{fig:introFig}, an overset grid consists of multiple overlapping structured component grids used to cover a complex, and perhaps moving, problem domain.  A mapping is defined for each component grid from physical space to a unit square (or cube) in a computational (index) space, and the mapped PDE is discretized in the computational space.  Solution values on the grid are interpolated at internal boundaries where two component grids overlap, and the given physical BCs are applied on external grid boundaries.  We have developed second-order accurate schemes for the equations of linear and nonlinear elasticity~\cite{smog2012,flunsi2016}, and up to fourth-order accurate schemes for the incompressible Navier-Stokes equations~\cite{ICNS,fins2020} and Maxwell's equations~\cite{max2006b,adegdm2019,adegdmi2020} using overset grids, among other applications.  We generally use the physical BCs, along with CBCs, to define discrete centered boundary conditions at external boundaries (with the aid of ghost points), but this approach would become increasingly difficult as the order of the approximation increases.  The difficulty stems from the algebraic complexity associated with taking higher and higher derivatives of the spatial operator of the mapped PDE and working out its attendant discrete approximations (with tangential couplings).  An associated difficulty involves the special treatments required at corners of the problem domain where separate BCs along sides meet.  The LCBC approach overcomes these difficulties by introducing a polynomial interpolant of the solution about each point on the boundary.  The polynomial degree is determined by the desired order of accuracy of the approximation, and the coefficients of the polynomial are specified by imposing constraints involving known solution values at grid points interior to the boundary, the physical BCs and CBCs.  This approach only requires CBCs defined at a continuous level, and these  conditions can be applied to the polynomial interpolant recursively thus easing the aforementioned algebraic complexity.  Once defined, the polynomial interpolant can be used to specify solution values at ghost points normal to the boundary (or in corner ghost points for the case of a domain corner) without tangential couplings.

\input tex/introFig.tex

The aim of the present paper is to describe the LCBC approach in detail for a general class of PDEs in second-order form and to investigate the properties of the resulting discretizations.  For example, in the case of a straight boundary and where the spatial operator is the Laplacian, it is well known that for Dirichlet (Neumann) boundary conditions the solution has odd (even) symmetry at the boundary.  This leads to simple numerical reflection conditions, and we show that the LCBC approach naturally results in these same reflection conditions (while one-sided approximations would not in general).  Beyond this special case, we show that the LCBC approach leads to accurate discretizations of the PDEs, and their BCs, for all orders of accuracy tested (up to sixth order).  Further, we show that there is no additional time-step restriction for stability for the case of explicit time-stepping schemes.  We focus here on linear PDEs, but the approach should be extendible to nonlinear problems as well.  In this article we focus on scalar PDEs, but the approach is also applicable to problems with vector PDEs (e.g.~the equations of linear elasticity and Maxwell's equations) and to problems with material interfaces.  Our ultimate goal is to automate the construction of CBC conditions for any order of accuracy and for a wide range of PDEs.  We believe that by using the LCBC approach that this goal is achievable.  This construction includes the development of LCBC conditions at grid faces as well as at grid corners for two-dimensional domains and at grid edges and vertices for three-dimensional domains.

Compatibility boundary conditions have been used with finite-difference methods for many years\footnote{For example, CBCs were known to Professor H.-O. Kreiss and his students at least since the 1980's.}, although it appears that the approach is not widely known.
In our work, we have used CBCs for second-order and fourth-order accurate approximations of the heat equation~\cite{BCNS} and the incompressible Navier-Stokes equations~\cite{ICNS,fins2020}.  For wave equations, we have described the use of CBCs for the compressible Euler equations~\cite{mog2006} and linear elasticity~\cite{smog2012}, and for high-order accurate approximations to Maxwell's equations~\cite{max2006b,adegdm2019}.  CBCs are also useful for problems involving material interfaces, such as conjugate heat transfer~\cite{th2009} and electromagnetics~\cite{max2006b,adegdmi2020}.  In recent work, we have developed Added-Mass Partitioned (AMP) schemes for a wide range of fluid-structure interaction (FSI) problems, including schemes for incompressible flows coupled to rigid bodies~\cite{rbins2017,rbinsmp2017,rbins3d2018} and elastic solids~\cite{fibr2019,fibrmp2019}.  These strongly-partitioned schemes incorporate AMP interface conditions derived using CBCs and the physical matching conditions at fluid-solid interfaces in order to overcome added-mass instabilities that can occur for the case of light bodies~\cite{fib2014,fis2014}.  In related work, we have also used CBCs in the CHAMP scheme~\cite{champ2016} to form discrete interface conditions for a partitioned approach to the solution of conjugate heat transfer problems.

In other work, CBCs are used in the book by Gustafsson 
on high-order difference methods~\cite{GustafssonHighOrderDifferenceMethodsBook2008}.
CBCs have also been used to derive stable and accurate embedded boundary\footnote{By embedded boundary we mean a boundary curve (or boundary surface in three dimensions) that passes through a grid in an irregular fashion (as opposed to a boundary-conforming grid).}
 approximations~\cite{EBWaveOrder2,EBelasticWaveOrder2,EBWaveOrder4}.
CBCs have been incorporated into summation-by-parts schemes by Sj\"ogreen and Petersson for the equations of elasticity~\cite{SjogreenPetersson2012}.
CBCs have been used by LeVeque and Li with their immersed interface method
to develop accurate approximations at embedded interfaces~\cite{IMElliptic,IMStokesFlow,IMFluidComplexGeometry}.
Shu and collaborators have used CBCs in their inverse-Lax-Wendroff approach for
hyperbolic equations and conservation laws~\cite{ILWConservationLaws,ILWConservationLawsSequel,ILWHyperbolicChangingWind,ILWHyperbolic}
as well as for parabolic and advection-diffusion equations~\cite{ILWStabilityDiffusion,ILWStabilityParabolic,ILWConvectionDiffusion}.

In this article we focus on high-order accurate finite-difference schemes. We note, however, that CBCs could also be useful for Galerkin schemes.
Typical high-order accurate FEM or DG schemes that use polynomial approximations over an element effectively use one-sided approximations near boundaries.
This can result in time-step restrictions that force the time-step to decrease rather significantly as the order-of-accuracy increases~\cite{hughes00,strang08,hesthaven08}. Similarly for B-Spline FEM, as commonly used in isogeometric analysis, one-sided operators occurring near boundaries results in spurious large eigenvalues, so-called outlier eigenvalues~\cite{hughesNURBS_2014}.  Banks et~al.~\cite{banks16_GDCont,banks18_GD2D,jacangelo20DSpline}, however, have shown that when CBCs are used with their Galerkin-Difference method, a class of FEM schemes, the spectrum of the operator is near-optimal, and the time-step restriction for explicit integration gives approximately the maximal CFL-one stability.

The remaining sections of the paper are organized as follows.  In Section~\ref{sec:MotivationalProblem}, we consider two sample problems to introduce the LCBC method and to compare it to the standard CBC approach.  A second-order PDE initial-boundary-problem is introduced in Section~\ref{sec:ModelEqn}, and this model problem is used along with its physical boundary conditions and derived CBCs as a basis to describe the LCBC approach for discretizations of the equations up to sixth-order accuracy.  In Section~\ref{sec:LCBCmethod}, we detail the steps of the LCBC approach for boundary conditions of Dirichlet and Neumann type, as well as the implementation of the method at domain corners where boundary conditions of various types meet.  We provide algorithms for the efficient application of the LCBC method for both domain sides and corners, and we discuss the conditioning of the matrices generated by the LCBC procedure.  Various elements of the LCBC method, including solvability, symmetry properties, and stability for wave equations, are analyzed in Section~\ref{sec:analysis}. The accuracy and stability of the method is illustrated in Section~\ref{sec:NumericalResults} by considering the numerical results for a variety of initial-boundary-value problems.   Concluding remarks and a discussion of future directions of the LCBC approach are given in Section~\ref{sec:Conclusion}.

%% file: tex/introFig.tex
{
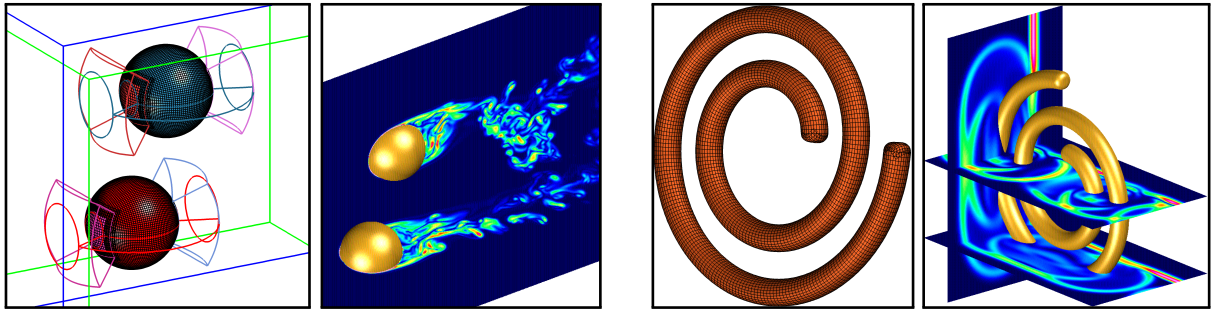
\begin{figure}[htb]
\begin{center}
\begin{tikzpicture}[scale=1]
  \useasboundingbox (0.25,0.45) rectangle (16,4.25);  

  \figByHeightb{0.0}{0}{fig/twoSpheresInAChannelGridZoom}{4cm}[.0][.0][.0][.0]
  \figByHeightb{4.2}{0}{fig/twoSpheresO4G8t13}{4cm}[.15][.3][.2][.2]

  \figByHeightb{8.6}{0}{fig/spiralWireGrid}{4cm}[.0][.0][.0][.0]

  \figByHeightb{12.2}{0}{fig/spiralWireEfieldNorm}{4cm}[.0][.0][.0][.0]
 
\end{tikzpicture}
\end{center}
  \caption{Some target applications for the new LCBC approach. Left: overset grid for two spherical bodies and computed incompressible flow (vorticity).
      Right: overset grid for a spiral wire and computed electromagnetic scattering.
   }
  \label{fig:introFig}
\end{figure}
}

%% file: tex/MotivationalProblem.tex
\section{Two sample problems}\label{sec:MotivationalProblem}

We begin by discussing CBCs and the LCBC approach for two relatively simple sample problems, the first involving the heat (diffusion) equation and the second involving the wave equation.  This is done to establish basic ideas and to introduce useful notation.

\subsection{Heat equation in one-dimension}

\newcommand{\Dcoeff}{\kappa}
Let $u(x,t)$ solve the following initial-boundary-value problem (IBVP) for the heat equation on the interval $x\in[a,b]$, with
a Dirichlet boundary condition on the left and a Neumann boundary condition on the right,
\begin{equation}\label{eq:heat1d}
	\begin{cases}
		\p_t u = \Dcoeff\, \p_x^2 u, & x\in(a,b), \quad t>0,\\
		u(a,t)=g_a(t), & t>0,\\
		\p_x u(b,t)=g_b(t), & t>0,\\
		u(x,0) = u_0(x) & x\in[a,b],
	\end{cases}
\end{equation}
where $\Dcoeff>0$ is a constant diffusivity, and $g_a(t)$, $g_b(t)$ and $u_0(x)$ are given smooth functions.  The CBCs we use for this problem are found by taking $\nu$ time derivatives of the boundary conditions in~\eqref{eq:heat1d}, and then using the PDE to replace time derivatives with space derivatives. This leads to a sequence of conditions given by 
\begin{equation}\label{eq:CBCs1d}   
\Dcoeff^\nu \p_x^{2\nu} u(a,t) = \p_t^\nu g_a(t), \qquad 
    \Dcoeff^\nu \p_x^{2\nu+1} u(b,t) = \p_t^\nu g_b(t ), \qquad \nu=1,2,\ldots
\end{equation}
The conditions in~\eqref{eq:CBCs1d}, denoted by ${\rm CBC}_s[\nu]$, $s=a$ and $b$, can be used in a finite difference scheme for the IBVP in~\eqref{eq:heat1d} as additional {\it numerical boundary conditions}, in place of extrapolation conditions or one-sided approximations.  For example, suppose the PDE and the Neumann boundary condition are discretized to order of accuracy $2p$ using central difference approximations with a stencil width $2p+1$, and suppose CBCs are used to generate additional numerical boundary conditions.  The number of derivatives in ${\rm CBC}_s[\nu]$ increases with $\nu$, but it can be shown that the CBCs can be approximated to lower-order accuracy as $\nu$ increases so that the stencil width remains $2p+1$.

To illustrate the use of CBCs for a specific case, consider a fourth-order accurate ($p=2$) approximation of~\eqref{eq:heat1d}.  Our focus is on the spatial approximations, especially near the boundaries, and so we consider a method-of-lines approach.  Let $U_j(t) \approx u(x_j,t)$, $j=-2,-1,0,\ldots,\Nx+2$, where $x_j=a+ j h$ are points on a uniform grid with grid spacing $h=(b-a)/\Nx$ as shown in Figure~\ref{fig:grid1d}.  Note that the grid includes {\it ghost points} at each boundary, and these are appended to the grid covering the problem domain $x\in[a,b]$ to facilitate the approximations of the spatial derivatives in the Neumann boundary condition at $x=b$ and the CBCs applied at both boundaries.

The fourth-order accurate approximation to~\eqref{eq:heat1d}, using a CBC at each boundary, is given by
\begin{equation}\label{eq:heat1dOrder4}
	\begin{cases}
		\p_t U_j(t) = \Dcoeff \, D_{4xx} U_j, & j=1,2,\ldots,\Nx, \quad t>0,\\
		U_0(t) = g_a(t), & t>0,\\
		\Dcoeff D_{4xx}  U_0(t) = \p_t g_a(t), & t>0,\qquad {\rm CBC}_a[1],\\
		\Dcoeff^2 (\Dpx\Dmx)^2 U_0(t) = \p_t^2 g_a(t), & t>0,\qquad {\rm CBC}_a[2],\\
		D_{4x} U_{\Nx}(t) = g_b(t), & t>0,\\
		\Dcoeff\Dzx\Dpx\Dmx U_{\Nx}(t) = \p_t g_b(t), & t>0,\qquad {\rm CBC}_b[1].\\
	\end{cases}
\end{equation}
Here, $D_{4x}$ and $D_{4xx}$ are fourth-order accurate approximations to $\p_x$ and $\p_x^2$, respectively, defined by
\begin{equation}\label{eq:Dops4}
D_{4x} \eqdef \Dzx\left( I - \frac{h^2}{6} \Dpx\Dmx \right),\qquad D_{4xx} \eqdef  \Dpx\Dmx\left( I - \frac{h^2}{12} \Dpx\Dmx \right),
\end{equation}
where $\Dpx$, $\Dmx$ and $\Dzx$ are the usual divided difference operators defined by
\begin{equation}\label{eq:Dops}
\Dpx U_j \eqdef \frac{U_{j+1}-U_j}{h}, \qquad \Dmx U_j \eqdef \frac{U_{j}-U_{j-1}}{h},\qquad \Dzx U_j \eqdef \frac{U_{j+1}-U_{j-1}}{2h}.
\end{equation}
Note that the discretization of the PDE is applied on all interior grid points and on the right boundary ($j=\Nx$) as a CBC for the Neumann boundary condition.  Note also that ${\rm CBC}_a[2]$ on the left and ${\rm CBC}_b[1]$ on the right are approximated only to second-order accuracy.  Despite this, the scheme is fourth-order accurate in space.  Initial conditions have been omitted in~\eqref{eq:heat1dOrder4} as these are not essential for the present discussion.

\input tex/gridFig1d

\newcommand{\uTilde}{\tilde{u}}
We now describe the LCBC method as an alternative approach to derive discrete boundary conditions in~\eqref{eq:heat1dOrder4}.  Let us first consider boundary conditions at $x=a$.  In the LCBC approach we approximate the solution near the boundary at $x=a$ as a polynomial $\uTilde(x)$ of degree $2p=4$.  The polynomial approximation changes in time as the solution evolves, but this dependence is suppressed here for notational convenience.  In terms of a monomial basis\footnote{In practice we use a Lagrange polynomial basis as it leads to better conditioned equations.}, for example, we have
\ba
   \uTilde(x) = \sum_{n=0}^4 c_n\, (x-a)^n . \label{eq:LCBCpoly1da}
\ea
The five coefficients $c_n$, $n=0,1,\ldots,4$, in~\eqref{eq:LCBCpoly1da} are found by requiring $\uTilde$ to match $U_j$ on the boundary and two interior points, as well as to satisfy the two compatibility conditions ${\rm CBC}_a[1]$ and ${\rm CBC}_a[2]$, i.e.
\bse
\label{eq:Aconditions}
\bat
     \uTilde(a) &= g_a(t),                     \quad&&  \\
     \uTilde(x_j) &= U_j(t),                     \quad&& j=1,2, \\
     \Dcoeff^\nu \p_x^{2\nu} \uTilde(a) &= \p_t^\nu g_a(t), \quad&& \nu=1,2,
\eat
\ese
for a fixed time $t$.  The five constraints in~\eqref{eq:Aconditions} imply a local system of linear equations for the coefficients $\{c_n\}$ in~\eqref{eq:LCBCpoly1da}.  Once the system is solved, ghost values of $U_j$ are then determined by setting
\ba
  U_j = \uTilde(x_j), \qquad j=-2,-1.
\ea
For the right boundary, a similar polynomial approximation to that in~\eqref{eq:LCBCpoly1da} is defined about $x=b$, and its coefficients are found using the constraints
\bse
\label{eq:rightconditions}
\bat
     \p_x\uTilde(b) &= g_b(t),                     \quad&&  \\
     \uTilde(x_j) &= U_j(t),                     \quad&& j=\Nx-2,\Nx-1,\Nx, \\
     \Dcoeff^\nu \p_x^{2\nu+1} \uTilde(b) &= \p_t^\nu g_b(t), \quad&& \nu=1.
\eat
\ese
The polynomial approximation about $x=b$ is then used to determine ghost values at $x_{\Nx+1}$ and $x_{\Nx+2}$.  We note that the LCBC approach uses CBCs at a continuous level to define local polynomial approximations instead of their discrete approximations as in~\eqref{eq:heat1dOrder4}.  However, for this simple one-dimensional problem, it can be shown that the LCBC approach yields an equivalent approximation to that in~\eqref{eq:heat1dOrder4}.  

\subsection{The wave equation in two dimensions} \label{sec:wave2d}

We now consider an IBVP for the wave equation (with unit wave speed) on a unit square domain, $\xv=(x,y)\in\Omega = (0,1)^2$.  For this problem, we let $u(\xv,t)$ solve
\begin{equation}\label{eq:wave2d}
	\begin{cases}
		\p_{t}^2 u = \Delta u, & \xv\in\Omega, \quad t>0,\\
		u = g_\ell(y,t)& \xv \in \p\Omega_\ell, \quad t>0,\\
		u = g_r(y,t)& \xv \in \p\Omega_r, \quad t>0,\\
		\p_yu = g_b(x,t) & \xv\in\p\Omega_b, \quad t>0,\\
		\p_yu = g_\tau(x,t) & \xv\in\p\Omega_\tau, \quad t>0,\\
		u = u_0(\xv),\quad \p_tu=u_1(\xv), & \xv \in\overline{\Omega},\quad t=0,
	\end{cases}
\end{equation}
where $\Delta$ is the Laplacian operator.  Dirichlet conditions are imposed on the left ($x=0$) and right ($x=1$) boundaries, denoted by $\p\Omega_\ell$ and $\p\Omega_r$, respectively, with given smooth functions $g_\ell(y,t)$ and $g_r(y,t)$.  Similarly, Neumann conditions are specified on the bottom ($y=0$) and top ($y=1$) boundaries, denoted by $\p\Omega_b$ and $\p\Omega_\tau$, respectively, with given smooth functions $g_b(y,t)$ and $g_\tau(y,t)$.  Initial conditions are determined by the given functions $u_0(\xv)$ and $u_1(\xv)$ defined for $\overline{\Omega}=\Omega\cup\p\Omega_s$, $s=\ell,r,b,\tau$.

As before, we consider CBCs for the problem by taking time derivatives of the boundary conditions and then replacing even time derivatives of $u(\xv,t)$ in favor of spatial derivatives using the PDE.  For example, $2\nu$ time derivatives of the Dirichlet conditions at $s=\ell$ gives
\begin{equation}
\Delta\sp\nu u(0,y,t)=\p_t\sp{2\nu}g_\ell(y,t),\qquad \nu=1,2,\ldots, 
\end{equation}
which implies a sequence of CBCs, denoted by ${\rm CBC}_\ell[\nu]$, $\nu=1,2,\ldots$, whose first two conditions are given by
\bse
\label{eq:cbcl}
\begin{align}
\p_x\sp2u(\xv,t)&=\p_t\sp2g_\ell(y,t) - \p_y\sp2g_\ell(y,t),  &\xv\in\p\Omega_\ell, \\ 
\p_x\sp4u(\xv,t)&=\p_t\sp4g_\ell(y,t) - 2\p_t\sp2\p_y\sp2g_\ell(y,t) + \p_y\sp4g_\ell(y,t), &\xv\in\p\Omega_\ell.
\end{align}
\ese
We can also write down CBCs corresponding to the Neumann boundary conditions by taking time derivatives.  At the bottom boundary, for example, this leads to ${\rm CBC}_b[\nu]$, $\nu=1,2,\ldots$, whose first two conditions are given by
\bse
\label{eq:cbcb}
\begin{align}
\p_y\sp3u(\xv,t)&=\p_t\sp2g_b(x,t) - \p_x\sp2g_b(x,t),  &\xv\in\p\Omega_b, \\ 
\p_y\sp5u(\xv,t)&=\p_t\sp4g_b(x,t) - 2\p_t\sp2\p_x\sp2g_b(x,t) + \p_x\sp4g_b(x,t), &\xv\in\p\Omega_b.
\end{align}
\ese
The CBCs along the boundaries with $s=r$ and $s=\tau$ are similar to those given in~\eqref{eq:cbcl} and~\eqref{eq:cbcb}, respectively.

A standard use of CBCs to define numerical boundary conditions for a high-order accurate discretization of the IBVP in~\eqref{eq:wave2d} for the wave equation in two space dimensions follows similar steps as were used above for the heat equation in one space dimension.  We first introduce a Cartesian grid for $\overline{\Omega}$ defined by
\begin{equation} \label{eq:grid}
\overline{\Omega}_h \eqdef \left\{ \xv_{\iv} = (x_i,y_j) = (i\dx,j\dy), \ i = 0,\dots, \Nx, \ j = 0,\ldots, \Ny \right\},
\end{equation}
where $\Nx$ and $\Ny$ determine the number of grid lines in the $x$ and $y$ directions, respectively, $\dx = 1/\Nx$ and $\dy = 1/\Ny$ are grid spacings, and $\iv = (i,j)$ is a multi-index, see the left plot of Figure~\ref{fig:grid2d}.  Let $\p\Omega_{s,h}$, $s=\ell,r,b,\tau$, represent the grid points on the left, right, bottom and top boundaries of $\overline{\Omega}_h$ respectively. The set $\p\Omega_h=\bigcup_{s = \ell,r,\tau,b}\p\Omega_{s,h}$ represents all the boundary grid points of $\overline{\Omega}_h$ and $\Omega_{h} = \overline{\Omega}_h\backslash\p\Omega_{h}$ denotes the interior grid points.  The grid includes ghost points, as shown in the figure, and these are used to facilitate the approximation of the Neumann boundary conditions and the CBCs.  Let $U_\iv(t) \approx u(x_i,y_j,t)$, and consider a fourth-order accurate semi-discrete scheme for~\eqref{eq:wave2d} given by
\begin{equation}\label{eq:wave2dOrder4}
	\begin{cases}
		\pt^2U_{\iv}(t) = \Delta_{4h}U_{\iv}, & \xv_{\iv} \in \Omega_{h}\cup\p\Omega_{b,h}\cup\p\Omega_{\tau,h},\quad t>0, \\
		U_{\iv}(t) = g_\ell(y_j,t), & \xv_{\iv} \in\p\Omega_{\ell,h},\quad t>0, \\
		U_{\iv}(t) = g_r(y_j,t), & \xv_{\iv} \in\p\Omega_{r,h},\quad t>0, \\
		D_{4y}U_{\iv}(t) = g_b(x_i,t), & \xv_{\iv} \in\p\Omega_{b,h},\quad t>0, \\
		D_{4y}U_{\iv}(t) = g_\tau(x_i,t), & \xv_{\iv} \in\p\Omega_{\tau,h},\quad t>0,
	\end{cases}
\end{equation}
where $\Delta_{4h}$ and $D_{4y}$ are fourth-order centered approximations of $\Delta$ and $\p_y$, respectively, following the definitions in~\eqref{eq:Dops4} and~\eqref{eq:Dops}.  Discrete approximations of ${\rm CBC}_\ell[1]$ and ${\rm CBC}_\ell[2]$ in~\eqref{eq:cbcl} given by
\bse
\label{eq:cbclapprox}
\begin{align}
D_{4xx}U_\iv(t)&=\p_t\sp2g_\ell(y_j,t) - \p_y\sp2g_\ell(y_j,t),  &\xv_\iv\in\p\Omega_{\ell,h}, \\ 
(D_{+x}D_{-x})\sp2U_\iv(t)&=\p_t\sp4g_\ell(y_j,t) - 2\p_y\sp2\p_y\sp2g_\ell(y_j,t) + \p_y\sp4g_\ell(y_j,t), &\xv_\iv\in\p\Omega_{\ell,h},
\end{align}
\ese
respectively, are added to~\eqref{eq:wave2dOrder4}.  Note that the approximation of ${\rm CBC}_\ell[1]$ is fourth-order accurate, while a second-order accurate approximation is used for ${\rm CBC}_\ell[2]$, following the scheme in~\eqref{eq:heat1dOrder4} for the CBCs about its Dirichlet boundary.  Similar approximations of ${\rm CBC}_r[1]$ and ${\rm CBC}_r[2]$ are also added to~\eqref{eq:wave2dOrder4}.  For the Neumann conditions, we use second-order accurate approximations of ${\rm CBC}_b[1]$ and ${\rm CBC}_\tau[1]$.  For example, the approximation of ${\rm CBC}_b[1]$ in~\eqref{eq:cbcb} given by
\begin{align}
D_{0y}D_{+y}D_{-y}U_\iv(t)&=\p_t\sp2g_b(x_i,t) - \p_x\sp2g_b(x_j,t), &\xv_\iv\in\p\Omega_{b,h},
\end{align}
is added to~\eqref{eq:wave2dOrder4}.  As before, the specifications of the initial conditions have been omitted in~\eqref{eq:wave2dOrder4}.

\input tex/gridFig2d

\newcommand{\ip}{{\hat{m}}}
\newcommand{\jp}{{\hat{n}}}
\newcommand{\ih}{{\hat{i}}}
\newcommand{\jh}{{\hat{j}}}
\newcommand{\ib}{{\tilde{i}}}
\newcommand{\jb}{{\tilde{j}}}

We now consider the LCBC approach for the fourth-order accurate approximation of the wave equation in~\eqref{eq:wave2dOrder4}.  For this two-dimensional problem, it is convenient to define an interpolating polynomial $\tilde{u}(x,y)$, centered about $(\tilde{x},\tilde{y})$, as
\begin{equation}\label{eqn:InterpolatingPolynomial}
	\tilde{u}(x,y) \eqdef \sum_{\jp = -p}^p\sum_{\ip = -p}^pd_{\ip,\jp}L_\ip\left(\dfrac{x - \tilde{x}}{\dx}\right)L_\jp\left(\dfrac{y - \tilde{y}}{\dy}\right), \qquad p\in\mathbb{N},
\end{equation}
where $d_{\ip,\jp}$ are coefficients and $L_k(z)$ is a Lagrange basis function given by
\begin{equation}
	L_k(z) = \prod_{\substack{l = -p\\ l\neq k}}^p{(z - l)\over(k-l)}, \qquad k = -p,\dots, p.
\end{equation}
Note that $\tilde{u}$ has the property 
\begin{equation}
	\tilde{u}\bigl(\tilde{x}+\ih\dx,\tilde{y} + \jh\dy\bigr) = d_{\ih,\jh}, \qquad \ih,\jh = -p,\dots, p.
\end{equation}
We set $p=2$ in~\eqref{eqn:InterpolatingPolynomial} for the fourth-order accurate scheme, and use known data at grid points in the interior and on the boundary, along with the physical boundary conditions and CBCs, to determine the coefficients of $\tilde{u}(x,y)$ centered at grid points along the boundaries.  For example, along the left Dirichlet boundary at a fixed time~$t$, let $(\tilde{x},\tilde{y}) = \xv_{0,\jb}$, $\jb\in\{2,3, \dots, \Ny-2\}$, and impose the constraints
\bse
\label{eq:leftconditions4}
\bat
    \tilde{u}\bigl(\tilde{x},\tilde{y} + \jh\dy\bigr) &= g_\ell\bigl(\tilde{y} + \jh\dy,t\bigr),                     && \quad \jh = -2,\dots, 2, \\
    \tilde{u}\bigl(\tilde{x}+\ih\dx,\tilde{y} + \jh\dy\bigr)  &= U_{\ih,\tilde{j}+\jh}(t),                    \quad& \ih = 1,2,& \quad \jh = -2,\dots, 2, \\
    \py^{\mu}\Delta^{\nu}\tilde{u}(\tilde{x},\tilde{y}) &= \py^{\mu}\pt^{2\nu}g_\ell(\tilde{y},t), \qquad& \nu=1,2,& \quad \mu = 0,\dots, 4.
\eat
\ese
The 25 constraints in~\eqref{eq:leftconditions4} lead to a linear system for the 25 coefficients $d_{\ip,\jp}$, $\ip,\jp = -2,\dots,2$, in~\eqref{eqn:InterpolatingPolynomial} for each point $\xv_{0,\jb}$ along the boundary.  Once solved, the coefficients determine the ghost points as $U_{\ih,\tilde{j}} = d_{\ih,0}$, $\ih=-2,-1$.  Similarly, along the bottom Neumann boundary, let $(\tilde{x},\tilde{y}) = \xv_{\ib,0}$, $\ib\in\{2,3, \dots, \Nx-2\}$, and impose the constraints
\bse
\label{eq:bottomconditions4}
\bat
    \tilde{u}\bigl(\tilde{x}+\ih\dx,\tilde{y} + \jh\dy\bigr)  &= U_{\ib+\ih,\jh}(t),                    \quad&\jh = 0,1,2, & \quad \ih = -2,\dots, 2, \\
    \px^{\mu}\Delta^{\nu}\p_y\tilde{u}(\tilde{x},\tilde{y}) &= \px^{\mu}\pt^{2\nu}g_b(\tilde{x},t), \qquad&\nu=0,1, & \quad \mu = 0,\dots, 4. \label{eq:cbcN}
\eat
\ese

Note that the Neumann boundary condition appears in~\eqref{eq:cbcN} for the case $\nu=\mu=0$, and tangential derivatives of this boundary condition are applied for $\mu>0$.  As before, the 25 constraints in~\eqref{eq:bottomconditions4} determine the 25 coefficients in~\eqref{eqn:InterpolatingPolynomial}, which then specify ghost points as $U_{\ib,\jh} = d_{0,\jh}$, $\jh=-2,-1$, for each point $\xv_{0,\jb}$ along the boundary.

The application of the LCBC approach at corners requires special treatment.  For example, let us consider the corner at $\tilde{\xv}=(0,0)$ where the Dirichlet boundary on the left meets the Neumann boundary along the bottom (see the right plot of Figure~\ref{fig:grid2d}).  Here, we choose a linearly independent set of constraints from those implied by~\eqref{eq:leftconditions4} and~\eqref{eq:bottomconditions4}.  We begin with the data available in the interior and on the Neumann boundary,
\bse
\label{eq:cornerconditions4}
\begin{equation}
\tilde{u}\bigl(\ih\dx,\jh\dy\bigr) = U_{\ih,\jh}(t),\qquad \ih = 1,2, \quad \jh = 0,1,2.
\end{equation}
Next, we impose tangential derivatives of the primary boundary conditions
\bat
    \py\sp{\mu_\ell}\tilde{u}(0,0) &= \py\sp{\mu_\ell}g_\ell(0,t),                     & \qquad \mu_\ell = 0,2,3,4, \label{eq:cornerDl}\\
    \px\sp{\mu_b}\py\tilde{u}(0,0) &= \px\sp{\mu_b}g_b(0,t),                     & \qquad \mu_b = 1,2,3,4, \label{eq:cornerNb}
\eat
along with
\begin{equation}
\py\tilde{u}(0,0) = {1\over2}\bigl(\py g_\ell(0,t)+g_b(0,t)\bigr).
\end{equation}
The latter is an average of the linearly dependent constraints in~\eqref{eq:cornerDl} and~\eqref{eq:cornerNb} corresponding to $\mu_\ell=1$ and $\mu_b=0$.  The remaining constraints and taken from compatibility conditions, and they are
\bat
    \py^{\mu_\ell}\Delta^{\nu}\tilde{u}(0,0) &= \py^{\mu_\ell}\pt^{2\nu}g_\ell(0,t), \qquad& \nu=1,2,& \quad \mu_\ell = 0,2,4, \label{eq:cornerDCl}\\
    \px^{\mu_b}\Delta\p_y\tilde{u}(0,0) &= \px^{\mu_b}\pt^{2}g_b(0,t), \qquad&& \quad \mu_b = 1,3,4, \label{eq:cornerNCb}
\eat
and
\begin{equation}
\px\sp2\py\sp3\tilde{u}(0,0) = {\dy\sp2\py\sp3\pt\sp2g_\ell(0,t)+\dx\sp2\px\sp2\pt\sp2g_b(0,t)\over\dy\sp2+\dx\sp2}. \label{eq:cornerAvg}
\end{equation}
\ese
Again, an average is used when the conditions obtained from the two boundaries are linearly dependent.  This occurs for the constraints in~\eqref{eq:cornerDCl} and~\eqref{eq:cornerNCb} corresponding to $(\nu,\mu_\ell)=(1,3)$ and $\mu_b=2$, respectively.  Note that we use a weighted average in~\eqref{eq:cornerAvg} to balance the tangential derivatives taken in the $x$ and $y$ directions, and that certain high derivatives of $\tilde u$ implied by the compatibility conditions have been omitted since they vanish for the chosen degree of the interpolating polynomial in~\eqref{eqn:InterpolatingPolynomial}.  There are 25 constraints in~\eqref{eq:cornerconditions4} that determine the 25 coefficients of $\tilde u(\xv)$ about $\tilde{\xv}=(0,0)$, and this interpolating polynomial is evaluated at $\xv_{\iv}$, for $\iv=[-2,-1]\times[-2,1]$ and $\iv=[0,1]\times[-2,-1]$, to determine the ghost values there.

Note that we have chosen to use the corner LCBC polynomial to assign solution values in corner ghost points, e.g.~$U_\iv$ for $\iv=[-2,1]\times[-2,-1]$, as well as values in nearby ghost points along the adjacent faces, e.g.~$U_\iv$ for $\iv=[-2,-1]\times[0,1]$ and $\iv=[0,1]\times[-2,-1]$, see Figure~\ref{fig:grid2d}.  These latter values could instead have been determined by special LCBC polynomials that used appropriate compatibility conditions from the primary face and the secondary adjacent face. For higher-order approximations this would result in many more special cases as the corner is approached.  To avoid these complications, we have chosen to employ the simpler approach of using the corner LCBC polynomial at all adjacent face ghost points where the standard LCBC polynomial for the face does not apply.

The application of the LCBC method just outlined for the two sample problems describes the basic approach. In the subsequent sections, we elaborate on the LCBC method by considering initial-boundary-value problems involving a broader class of linear PDEs and corresponding spatial discretizations of higher orders of accuracy.

%% file: tex/gridFig1d.tex
\newcommand{\labelFont}{\small}
\begin{figure}[hbt!]
  \centering
  \begin{tikzpicture}
    \useasboundingbox (0,0) rectangle (12,2);

    \begin{scope}[yshift=.25cm]
    \draw[ultra thick,black] (2,0) -- (10,0);
      \foreach \x in {2,10}
        \draw[ultra thick,black] (\x,-4pt) -- (\x,+4pt);  
  
      \draw ( 2,0) node[anchor=north,yshift=-3pt] {\labelFont$x=a$};         
      \draw (10,0) node[anchor=north,yshift=-3pt] {\labelFont$x=b$}; 
    \end{scope}        

    \draw[ultra thick,blue] (0,1) -- (12,1);
    \foreach \x in {0,1,2,3,4,5,8,9,10,11,12}
      \draw[ultra thick,blue] (\x,1cm-4pt) -- (\x,1cm +4pt); 
    \draw ( 0,1) node[anchor=north,yshift=-3pt] {\labelFont$x_{-2}$} node[anchor=south,yshift=3pt]{\labelFont$V_{-2}$};
    \draw ( 1,1) node[anchor=north,yshift=-3pt] {\labelFont$x_{-1}$} node[anchor=south,yshift=3pt]{\labelFont$V_{-1}$};
    \draw ( 2,1) node[anchor=north,yshift=-3pt] {\labelFont$x_{0}$}  node[anchor=south,yshift=3pt]{\labelFont$V_{0}$};
    \draw ( 3,1) node[anchor=north,yshift=-3pt] {\labelFont$x_{1}$}  node[anchor=south,yshift=3pt]{\labelFont$V_{1}$};
    \draw ( 4,1) node[anchor=north,yshift=-3pt] {\labelFont$x_{2}$}  node[anchor=south,yshift=3pt]{\labelFont$V_{2}$};

    \draw ( 5,1) node[anchor=north,yshift=-3pt] {$\cdots$}  node[anchor=south,yshift=5pt]{$\cdots$};

    \draw ( 9,1) node[anchor=north,yshift=-3pt] {\labelFont$x_{\Nx-1}$}  node[anchor=south,yshift=3pt]{\labelFont$V_{\Nx-1}$};
    \draw (10,1) node[anchor=north,yshift=-3pt] {\labelFont$x_{\Nx}$}  node[anchor=south,yshift=3pt]{\labelFont$V_{\Nx}$};
    \draw (11,1) node[anchor=north,yshift=-3pt] {\labelFont$x_{\Nx+1}$}  node[anchor=south,yshift=3pt]{\labelFont$V_{\Nx+1}$};
    \draw (12,1) node[anchor=north,yshift=-3pt] {\labelFont$x_{\Nx+2}$}  node[anchor=south,yshift=3pt]{\labelFont$V_{\Nx+2}$};

%
%
  \end{tikzpicture}
  \caption{Grid in one space dimension, with ghost points, for the fourth-order accurate approximation of the heat equation in~\eqref{eq:heat1dOrder4}.}
  \label{fig:grid1d}
\end{figure}
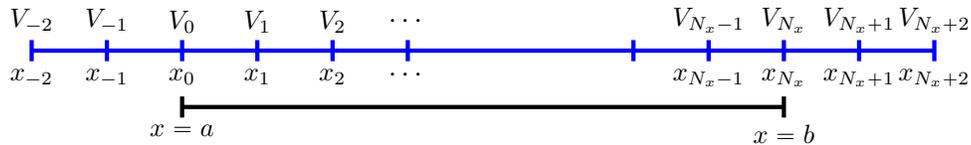

%% file: tex/gridFig2d.tex
\usetikzlibrary{shapes,positioning,intersections,quotes}
\usetikzlibrary{shapes.geometric}
\usetikzlibrary{arrows.meta}
\tikzset{
	ring shading/.code args={from #1 at #2 to #3 at #4}{
		\def\colin{#1}
		\def\radin{#2}
		\def\colout{#3}
		\def\radout{#4}
		\pgfmathsetmacro{\proportion}{\radin/\radout}
		\pgfmathsetmacro{\outer}{.8818cm}
		\pgfmathsetmacro{\inner}{.8818cm*\proportion}
		\pgfmathsetmacro{\innerlow}{\inner-0.01pt}
		\pgfdeclareradialshading{ring}{\pgfpoint{0cm}{0cm}}%
		{
			color(0pt)=(white);
			color(\innerlow)=(white);
			color(\inner)=(#1);
			color(\outer)=(#3)
		}
		\pgfkeysalso{/tikz/shading=ring}
	},
}
\definecolor{pinegreen}{rgb}{0.0, 0.47, 0.44}
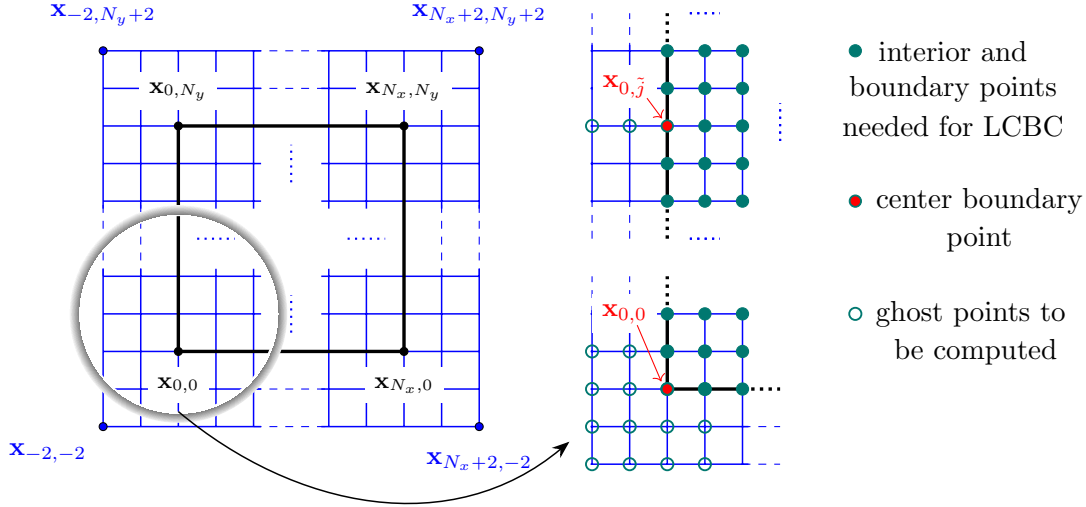
\begin{figure}[H]
	\centering
	\begin{tikzpicture}
		\draw[step=.5cm,blue,line width=0.2mm] (-1,1.9) grid (1.1,4);
		\draw[step=.5cm,blue,line width=0.2mm] (-1,-1) grid (1.1,1.1);
		\draw[step=.5cm,blue,line width=0.2mm] (1.9,-1) grid (4,1.1);
		\draw[step=.5cm,blue,line width=0.2mm] (1.9,1.9) grid (4,4);
		
		\foreach \x in {0,1,9,10}
		\draw[dashed,blue] (-1+.5*\x,1) -- (-1+.5*\x,2);
		\foreach \y in {0,1,9,10}
		\draw[dashed,blue] (1,-1+.5*\y) -- (2,-1+.5*\y);
		\foreach \x in {0,2}
		\draw[dotted,thick,blue] (\x+0.25,1.5) -- (\x+0.75,1.5);
		\foreach \y in {0,2}
		\draw[dotted,thick,blue] (1.5,\y+0.25) -- (1.5,\y+0.75);
		
		\draw[line width=0.45mm,-] (0,0) -- (3,0);
		\draw[line width=0.45mm,-] (0,0) -- (0,3);
		\draw[line width=0.45mm,-] (3,0) -- (3,3);
		\draw[line width=0.45mm,-] (0,3) -- (3,3);
		
		\filldraw (0,0) circle[radius=1.5pt];
		\node[below=4pt of {(0,0)}, outer sep=2pt,fill=white] {\footnotesize{$\xv_{0,0}$}};
		\filldraw (3,0) circle[radius=1.5pt];
		\node[below=4pt of {(3,0)}, outer sep=2pt,fill=white] {\footnotesize{$\xv_{\Nx,0}$}};
		\filldraw (3,3) circle[radius=1.5pt];
		\node[above=4pt of {(3,3)}, outer sep=2pt,fill=white] {\footnotesize{$\xv_{\Nx,\Ny}$}};
		\filldraw (0,3) circle[radius=1.5pt];
		\node[above=4pt of {(0,3)}, outer sep=2pt,fill=white] {\footnotesize{$\xv_{0,\Ny}$}};
		
		\draw (-1,-1) circle[radius=1.5pt] [fill = blue];
		\node[below left=1pt of {(-1,-1)}, outer sep=2pt,fill=white] {\small{\textcolor{blue}{$\xv_{-2,-2}$}}};
		\draw (4,-1) circle[radius=1.5pt] [fill = blue];
		\node[below=4pt of {(4,-1)}, outer sep=2pt,fill=white] {\small{\textcolor{blue}{$\xv_{\Nx+2,-2}$}}}; 
		\draw (-1,4) circle[radius=1.5pt] [fill = blue];
		\node[above=4pt of {(-1,4)}, outer sep=2pt,fill=white] {\small{\textcolor{blue}{$\xv_{-2,\Ny+2}$}}};
		\draw (4,4) circle[radius=1.5pt] [fill = blue];
		\node[above=4pt of {(4,4)}, outer sep=2pt,fill=white] {\small{\textcolor{blue}{$\xv_{\Nx +2,\Ny+2}$}}};
		
		\draw[step=.5cm,blue,line width=0.2mm] (5.48,1.98) grid (7.58,4.08); 
		\foreach \x in {0,1}
		\draw[dashed,blue] (5.5+.5*\x,4) -- (5.5+.5*\x,4.5);
		\foreach \x in {0,1}
		\draw[dashed,blue] (5.5+.5*\x,1.5) -- (5.5+.5*\x,2.0);
		\draw[line width=0.45mm,-] (6.5,2) -- (6.5,4); \draw[dotted,line width=0.45mm] (6.5,4) -- (6.5,4.5);\draw[dotted,line width=0.45mm] (6.5,1.5) -- (6.5,2);
		\draw[dotted,thick,blue] (8,2.8) -- (8,3.3); \draw[dotted,thick,blue] (6.8,4.5) -- (7.2,4.5); \draw[dotted,thick,blue] (6.8,1.5) -- (7.2,1.5);
		\foreach \x in {0,1,2}{
			\foreach \y in {0,1,2,3,4}{
				\draw[pinegreen] (6.5+0.5*\x,2+0.5*\y) circle[radius=2.3pt][fill = pinegreen];
		}}
		\draw[red] (6.5,3) circle[radius=1.6pt] [fill = red];
		\node[above left=2.5pt of {(6.5,3.1) }, outer sep=2.3pt,fill=white] {\small{\textcolor{red}{$\xv_{0,\tilde{j}}$}}};
		\draw[<-,red] (6.43,3.1) -- (6.23,3.37); 
		\foreach \x in {0,1}
		{\draw[line width=0.25mm,pinegreen] (5.5+0.5*\x,3) circle[radius=2.3pt];}
		
		\draw[step=.5cm,blue,line width=0.2mm] (5.48,-1.5) grid (7.58,0.55); 
		\foreach \x in {0,1}
		\draw[dashed,blue] (5.5+.5*\x,0.5) -- (5.5+.5*\x,1);
		\foreach \y in {0,1}
		\draw[dashed,blue] (7.5,-1.5+.5*\y) -- (8,-1.5+.5*\y);
		\draw[line width=0.45mm,-] (6.5,-.5) -- (6.5,0.5); 
		\draw[dotted,line width=0.45mm] (6.5,0.5) -- (6.5,1);
		\draw[line width=0.45mm,-] (6.5,-.5) -- (7.5,-.5); 
		\draw[dotted,line width=0.45mm] (7.5,-.5) -- (8,-.5);
		\draw[dotted,thick,blue] (8,2.8) -- (8,3.3); \draw[dotted,thick,blue] (6.8,4.5) -- (7.2,4.5); \draw[dotted,thick,blue] (6.8,1.5) -- (7.2,1.5);
		\foreach \x in {0,1,2}{
			\foreach \y in {0,1,2}{
				\draw[pinegreen] (6.5+0.5*\x,-0.5+0.5*\y) circle[radius=2.3pt][fill = pinegreen];
		}}
		\draw[red] (6.5,-.5) circle[radius=1.6pt] [fill = red];
		\node[above left=2.5pt of {(6.5,0.1) }, outer sep=2.3pt,fill=white] {\small{\textcolor{red}{$\xv_{0,0}$}}};
		\draw[<-,red] (6.42,-0.4) -- (6.1,0.35); 
		\foreach \y in{0,1,2,3}{
		\foreach \x in {0,1,2,3}{
		\draw[line width=0.25mm,pinegreen] (5.5+0.5*\x,-1.5+0.5*\y) circle[radius=2.3pt];
	}}
		
		\shade[even odd rule,ring shading={from gray at 1.5 to white at 1.65}]
		(0,0.5) circle (1.45) circle (1.3);
		\draw [-{Stealth[scale=1.2]},line width=0.2mm](0,-0.8) to [out=320,in=220] (5.2,-1.20);
		
		\draw[pinegreen] (9,4) circle[radius=2.3pt][fill = pinegreen] (10.3,4) node{\textcolor{black}{interior and}}
		(10.3,3.5) node{\textcolor{black}{boundary points}}
		(10.3,3) node{\textcolor{black}{needed for LCBC}};
		\draw[pinegreen] (9,2) circle[radius=2.3pt][fill = pinegreen];
		\draw[red] (9,2) circle[radius=1.6pt] [fill = red] (10.65,2) node{\textcolor{black}{center boundary}}
		(10.65,1.5) node{\textcolor{black}{point}};
		\draw[line width=0.25mm,pinegreen] (9,0.5) circle[radius=2.3pt] (10.6,0.5) node{\textcolor{black}{ghost points to }}
		(10.6,0) node{\textcolor{black}{be computed}};
	\end{tikzpicture}
\caption{Grid in two space dimension, with ghost points, for the fourth-order accurate approximation of the wave equation in~\eqref{eq:wave2dOrder4}.}  \label{fig:grid2d}
\end{figure}

%% file: tex/ModelEqn.tex
\section{Second-order PDE initial-boundary-value problems and discretizations}\label{sec:ModelEqn}

In this section, we consider initial-boundary-value problems for a general scalar second-order PDE and corresponding high-order accurate finite-difference approximations as a basis for a full description of the LCBC approach which follows in the next section.

\subsection{Second-order PDE initial-boundary-value problems}\label{subsec:MainEqn}

Consider the initial-boundary-value problem on $[0,T]\times \Omega$, $T>0$, given by
\begin{equation}\label{eqn:MainEqn}
	\begin{cases}
		\Lc_q u = Qu + f(\xv,t), & \xv \in \Omega, \quad t\in(0,T], \quad q = 0,1,2,\\
		\Bc u(\xv,t) = g(\xv,t), & \xv \in\p\Omega, \quad t\in[0,T], \\
		\pt^{\alpha - 1}u(\xv,0) = u_{\alpha - 1}(\xv), & \xv\in \overline{\Omega}, \quad \alpha = 1, \dots, q, \quad q  = 1,2.
	\end{cases}
\end{equation}
Here, $\Omega \subset \mathbb{R}^2$ is a general domain, $\p\Omega$ denotes the boundary of $\Omega$, and $\overline{\Omega} = \Omega\cup\p\Omega$ as before.  We define
\begin{equation}
	\Lc_qu \eqdef \begin{cases}
		0, & q = 0,\\
		\pt^q u, & q = 1,2,
	\end{cases}
\end{equation}
and
\begin{equation}\label{eqn:OperatorQ}
	Qu \eqdef c_{11}(\xv)\p_x\sp2u + 2c_{12}(\xv)\p_x\p_yu + c_{22}(\xv)\p_y\sp2u + c_{1}(\xv)\p_xu + c_{2}(\xv)\p_yu + c_{0}(\xv)u, \qquad \xv \in \Omega.
\end{equation}
We assume that the coefficient functions $c_{11}(\xv)$, $c_{12}(\xv)$, etc., are smooth, and they are chosen, together with the boundary and initial conditions, so that the problem is well posed.  
\wdh{For example, for $q=1,2$,} necessary conditions are that $c_{11}(\xv)>0$, $c_{22}(\xv)>0$ and 
\begin{equation}
	c_{11}(\xv)c_{22}(\xv) - c_{12}^2(\xv) \ge\delta>0, \qquad \forall \xv\in\Omega.
\end{equation}
We note that~\eqref{eqn:OperatorQ} is taken in non-conservative form for the purposes of this article; LCBC methods for problems in conservative form are left to future work.

The governing equation in~\eqref{eqn:MainEqn}, with given forcing function $f(\xv,t)$, takes the form of a general elliptic ($q=0$), parabolic ($q=1$) or hyperbolic ($q=2$) PDE in second-order form depending on the choice of the index~$q$.  The boundary conditions in~\eqref{eqn:MainEqn}, with given forcing function $g(\xv,t)$, are written in terms of the boundary operator given by
\begin{equation} \label{eqn:OperatorB}
	\Bc u \eqdef b_{1}(\xv)u + b_{2}(\xv)\p_n u, \qquad \xv \in \partial\Omega,
\end{equation}
where $\p_n$ is the outward normal derivative and the coefficient functions satisfy $|b_1(\xv)|+|b_2(\xv)|\ne0$, $\forall \xv\in\p\Omega$.  Initial conditions are specified for the cases, $q=1$ and~$2$.

We are motivated by the application of the LCBC method for high-order accurate discretizations of the model problem in~\eqref{eqn:MainEqn} on mapped grids.  For such discretizations, we consider a smooth mapping from the unit square to $\Omega$. The form of the model problem remains unchanged in the mapped domain, so it suffices to study the governing equations in~\eqref{eqn:MainEqn} over the domain $\Omega = (0,1)^2$ as a model problem.

\subsection{Semi-discrete approximations}
Let
\begin{equation*}
	U_{\iv} \approx	\begin{cases}
		u(\xv_{\iv}), & q = 0, \\
		u(\xv_{\iv},t), & q = 1,2,
	\end{cases}
\end{equation*} 
represent the numerical approximation of the exact solution of \eqref{eqn:MainEqn} at discrete points $\xv_{\iv}$ on the Cartesian grid $\overline{\Omega}_h$ given in~\eqref{eq:grid} and at a fixed time $t\in[0,T]$.  Our principal focus is on discretizations of~\eqref{eqn:MainEqn} to fourth and sixth-order accuracy, although we also consider second-order accurate approximations as a baseline.  A second-order accurate discretization of~\eqref{eqn:MainEqn} employs standard centered differences given by
\begin{equation}\label{eqn:opsTwo}
D_{2,\zeta} \eqdef \dzze,\qquad D_{2,\zeta\zeta} \eqdef \dpmze,\qquad \zeta = x,y.
\end{equation}
Fourth-order accurate operators, $D_{4,\zeta}$ and $D_{4,\zeta\zeta}$, $\zeta=x,y$, were introduced in the previous section, and we define sixth-order accurate centered operators as
\begin{equation}\label{eqn:opsSix}
   \left.
   \begin{array}{l}
   \displaystyle{
	D_{6,\zeta} \eqdef \dzze \left(I - \frac{\dze^2}{6}\dpmze+\dfrac{\dze^4}{30}\left(\dpmze\right)^2\right),
	}\medskip\\
   \displaystyle{
	D_{6,\zeta\zeta} \eqdef \dpmze\left(I - \dfrac{\dze^2}{12}\dpmze+\dfrac{\dze^4}{90}\left(\dpmze\right)^2\right),
	}
	\end{array}\quad\right\}\qquad \zeta = x,y.
\end{equation}
Thus, 
\begin{equation}\label{eqn:Qh}
	Q_{d,h} \eqdef c_{11}(\xv_{\iv})D_{d,xx} + 2c_{12}(\xv_{\iv})D_{d,x}D_{d,y} + c_{22}(\xv_{\iv})D_{d,yy} + c_1(\xv_{\iv})D_{d,x} + c_{2}(\xv_{\iv})D_{d,y} + c_0(\xv_{\iv})I,
\end{equation}
for $\xv_{\iv}\in \Omega_h$.  Set $\Bc_{d,h}$ to be the $d\sp{{\rm th}}$-order accurate centered discretization of the boundary operator~$\Bc$.

In view of the finite difference operators defined above, we note that~\eqref{eqn:Qh} can be written as
\begin{equation}\label{eqn:Qhexp}
\begin{array}{l}
\displaystyle{
Q_{d,h}V_\iv=c_{11}(\xv_{\iv})D_{+x}D_{-x}\left\{V_\iv+\sum_{n=1}\sp{k-1}a_{n}\dx\sp{2n}\bigl(D_{+x}D_{-x}\bigr)\sp{n}V_\iv\right\}
}\smallskip\\
\qquad\qquad\displaystyle{
+2c_{12}(\xv_{\iv})D_{0x}D_{0y}\left\{V_\iv+\sum_{n=1}\sp{k-1}\sum_{l=0}\sp{n}b_{l}b_{n-l}\dx\sp{2l}\dy\sp{2(n-l)}\bigl(D_{+x}D_{-x}\bigr)\sp{l}\bigl(D_{+y}D_{-y}\bigr)\sp{n-l}V_\iv\right\}
}\smallskip\\
\qquad\qquad\;\;\displaystyle{
+c_{22}(\xv_{\iv})D_{+y}D_{-y}\left\{V_\iv+\sum_{n=1}\sp{k-1}a_{n}\dy\sp{2n}\bigl(D_{+y}D_{-y}\bigr)\sp{n}V_\iv\right\}
}\smallskip\\
\qquad\qquad\;\;\;\;\displaystyle{
+c_{1}(\xv_{\iv})D_{0x}\left\{V_\iv+\sum_{n=1}\sp{k-1}b_{n}\dx\sp{2n}\bigl(D_{+x}D_{-x}\bigr)\sp{n}V_\iv\right\}
}\smallskip\\
\qquad\qquad\;\;\;\;\;\;\displaystyle{
+c_{2}(\xv_{\iv})D_{0y}\left\{V_\iv+\sum_{n=1}\sp{k-1}b_{n}\dy\sp{2n}\bigl(D_{+y}D_{-y}\bigr)\sp{n}V_\iv\right\} + c_0(\xv_\iv)V_\iv,
}
\end{array}
\end{equation}
for some grid function $V_\iv$ and for $d=2k$, $k=1,2,3$.  The higher-order correction terms given by the sums in~\eqref{eqn:Qhexp} are omitted for the second-order accurate case with $k=1$.  The coefficients $a_n$ and $b_n$ of the correction terms involving powers of $D_{+x}D_{-x}$ and $D_{+y}D_{-y}$ are given by
\begin{equation}\label{eq:QhexpCoeffs}
a_1=-{1\over12},\quad a_2={1\over90},\qquad b_0=1,\quad b_1=-{1\over6},\quad b_2={1\over30}.
\end{equation}
The form of the discrete operator in~\eqref{eqn:Qhexp} is particularly useful for the calculation of $(Q_{d,h})\sp\nu V_\iv$, $\nu=1,2,\ldots$, which is required for the approximation of compatibility boundary conditions as is discussed in the next section.

The semi-discrete model problem takes the form 
\begin{equation}\label{eqn:SemiDiscreteEqn}
	\begin{cases}
		\Lc_q U_{\iv}(t) = Q_{d,h}U_{\iv}(t) + f(\xv_{\iv},t), & \xv_{\iv} \in \Omega_h, \quad t\in(0,T], \quad q = 0,1,2,\\
		\Bc_{d,h} U_{\iv}(t) = g(\xv_{\iv},t), & \xv_{\iv} \in\p\Omega_h, \quad t\in[0,T], \\
		\pt^{\alpha - 1}U_{\iv}(0) = u_{\alpha - 1}(\xv_{\iv}), & \xv_{\iv}\in \overline{\Omega}_h, \quad \alpha = 1, \dots, q, \quad q  = 1,2.
	\end{cases}
\end{equation}
Grid points along ghost lines at each boundary of $\Omega_h$ are introduced to accommodate the stencil of the discrete spatial operators near the boundaries, and these are included in the extended grid defined by
\begin{equation} \label{eqn:extDomain}
	\Omega_{h}^e \eqdef \left\{\xv_{\iv} \ \vert\ \iv = (i,j), \ i = -p, \dots, \Nx+p, \ j = -p,\dots, \Ny+p\right\},
\end{equation}
where $p=d/2$.  We evaluate the solution at the ghost points using the LCBC method.

\subsection{Compatibility boundary conditions}\label{subsec:CBC}
The LCBC method uses compatibility boundary conditions obtained from the primary boundary conditions and the governing PDE (and its derivatives) applied on the boundary. The steps taken to obtain the CBCs depend on whether the PDE is time-dependent or not.  We begin with the time-dependent cases for which $q=1$ or~$2$.  Here, we take $q$ time derivatives of the primary boundary condition in~\eqref{eqn:MainEqn} to give
\begin{equation}\label{eqn:CBCstep1}
	\Bc\pt^{q} u(\xv,t) = \pt^{q}g(\xv,t), \qquad \xv \in \p\Omega,
\end{equation}
at a fixed time $t\in[0,T]$.  Now apply the PDE from \eqref{eqn:MainEqn} to obtain
\begin{equation}\label{eqn:CBC1}
	\Bc Qu(\xv,t) = \pt^qg(\xv,t) - \Bc f(\xv,t), \qquad \xv\in \p\Omega. 
\end{equation}
Repeat the process $\nu$ times to find the $\nu\sp{{\rm th}}$ compatibility condition, denoted by ${\rm CBC}_{\Bc,q}[\nu]$, as
\begin{equation}\label{eqn:CBCq}
	\Bc Q^{\nu}u(\xv,t) = \pt^{q\nu}g(\xv,t) - \Bc\Psi_{\nu}f(\xv,t), \qquad \xv\in \p\Omega, \quad\nu=1,2,\ldots,
\end{equation}
where $\Psi_{\nu}$ is a differential operator defined by
\begin{equation} \label{eq:PsiDef}
	\Psi_{\nu} f(\xv,t) \eqdef \sum_{k = 1}^{\nu} Q^{k-1}\pt^{q(\nu - k)}f(\xv,t), \qquad \xv\in \p\Omega, \quad\nu=1,2,\ldots.
\end{equation} 

For the time-independent case with $q = 0$, we derive the $\nu\sp{{\rm th}}$ compatibility condition by first applying the elliptic operator $Q$ to the governing PDE $\nu -1$ times to obtain
\begin{equation}\label{eqn:ellipticComp}
	Q^{\nu}u(\xv) = - Q^{\nu-1}f(\xv), \qquad \xv \in \overline\Omega, \quad \nu=1,2,\ldots.
\end{equation}
The boundary operator in~\eqref{eqn:MainEqn} is then applied to \eqref{eqn:ellipticComp} to give
\begin{equation}\label{eqn:CBC0}
	\Bc Q^{\nu}u(\xv) = -\Bc Q^{\nu - 1}f(\xv), \qquad \xv \in \p\Omega, \quad\nu=1,2,\ldots,
\end{equation}
which we denote by ${\rm CBC}_{\Bc,0}[\nu]$.  We observe that the form of the CBCs in~\eqref{eqn:CBCq} for the cases $q=1$ and $2$ and in~\eqref{eqn:CBC0} for the case $q=0$ are similar, and these conditions provide the additional constraints needed to complete the specification of the local interpolants of the solution on the boundary for the LCBC method.

%% file: tex/LCBCapproach.tex
\section{LCBC method}\label{sec:LCBCmethod}

We now provide a full description of the LCBC method for the semi-discrete model problem in~\eqref{eqn:SemiDiscreteEqn}.  We first consider a coordinate boundary away from corners where two coordinate boundaries meet. We choose a Dirichlet-type boundary condition and introduce the LCBC method using a \textit{direct approach}. For a more efficient implementation, we improve upon this direct approach by adopting a stencil representation of the solution at the ghost points; we call this improved method the \textit{stencil approach}. After this is done, we present an example where the boundary condition is of the Neumann type. Finally, we describe the treatment near the corner and discuss the conditioning of linear systems that arise in the LCBC method.

\newcommand{\mt}{{\tilde{m}}}

\subsection{Dirichlet boundary} \label{sec:LCBCside}

As an example of the LCBC method for Dirichlet boundary conditions, let us consider the left boundary, $x=0$ with $y\in[0,1]$, and assume that the boundary operator in~\eqref{eqn:MainEqn} becomes
\begin{equation}
  u(\xv,t) = g_\ell(y,t), \qquad \xv\in \p\Omega_\ell,
\end{equation}
for a fixed time $t\in[0,T]$.  Following the approach described in Section~\ref{sec:wave2d} for the wave equation defined on a unit square domain, we consider the interpolating polynomial $\tilde u(\xv)$ in~\eqref{eqn:InterpolatingPolynomial} centered about a point $\tilde{\xv}\in\p\Omega_\ell$, and then specify its $\mt=(2p+1)\sp2$ coefficients $d_{\ip,\jp}$, $\ip,\jp=-p,\ldots,p$, by enforcing the constraints
\bse
\label{eq:leftconditions}
\bat
    \tilde{u}\bigl(0,\tilde{y} + \jh\dy\bigr) &= g_\ell\bigl(\tilde{y} + \jh\dy,t\bigr),                     && \quad \jh = -p,\dots, p, \label{eq:leftbc} \\
    \tilde{u}\bigl(\ih\dx,\tilde{y} + \jh\dy\bigr)  &= U_{\ih,\tilde{j}+\jh}(t),                    \quad& \ih = 1,\dots,p,& \quad \jh = -p,\dots, p,\label{eq:leftinterior} \\
    \py^{\mu}Q^{\nu}\tilde{u}(0,\tilde{y}) &= \py^{\mu}R_{\ell,\nu}(\tilde{y},t), \qquad& \nu=1,\ldots,p,& \quad \mu = 0,\dots,2p, \label{eq:leftcbc}
\eat
\ese
where
\begin{equation} \label{eq:Rdef}
R_{\ell,\nu}(y,t)\eqdef
\begin{cases}
-Q^{\nu - 1}f(0,y), & q = 0,  \\[5pt]
\pt^{q\nu}g_\ell(y,t) - \Psi_{\nu}f(0,y,t), & q=1,2.
\end{cases}
\end{equation}
The constraints in~\eqref{eq:leftbc} are the Dirichlet boundary condition applied at $2p+1$ grid points about the boundary point $(0,\tilde y)$, while~\eqref{eq:leftinterior} sets $\tilde{u}$ equal to $U_\iv$ at $p(2p+1)$ grid points interior to the boundary point.  The last constraints in~\eqref{eq:leftcbc} require that $\tilde{u}$ satisfy $(2p+1)$ tangential derivatives of the compatibility boundary conditions, ${\rm CBC}_{\ell,q}[\nu]$, $\nu=1,\ldots,p$, evaluated at the boundary point~$(0,\tilde y)$.  Together, the constraints in~\eqref{eq:leftconditions} imply $\mt=(2p+1)\sp2$ linear equations for the $\mt$ coefficients in~$\tilde{u}$ for each point $(0,\tilde y)\in \p\tilde{\Omega}_{\ell,h}$, where
\begin{equation}\label{eq:discreteBoundary}
  \p\tilde{\Omega}_{\ell,h} \eqdef \bigl\{\,\xv_{\iv}\ \vert \ i = 0, \ j = p,p+1,\ldots,\Ny-p\,\bigr\},
\end{equation}
is the set of grid points along the left boundary $x=0$ sufficiently separated from the corners at $y=0$ and~$1$. 

The $\mt\times\mt$ linear system implied by~\eqref{eq:leftconditions} has the form
\begin{equation}\label{eq:leftSys}
A\dv = \bv,
\end{equation}
where $A\in\mathbb{R}^{\mt\times\mt}$ is a coefficient matrix, $\bv\in\mathbb{R}^{\mt}$ is a right-hand side vector and $\dv\in\mathbb{R}^{\mt}$ is a vector containing the coefficients of the interpolating polynomial in~\eqref{eqn:InterpolatingPolynomial} organized as
\begin{equation}\label{eqn:vectord}
  \dv = \left[d_{-p,-p}, \dots, d_{-p,p}\,\vert\, d_{-p+1,-p},\dots, d_{-p+1,p}\,\vert\cdots\vert\, d_{p,-p},\dots,d_{p,p}\right]^{T}.
\end{equation}
The matrix $A$, as constructed in Algorithm~\ref{alg:formA} for a point $\tilde{\xv}$ on the boundary, has the $2\times2$ block structure
\begin{equation}\label{eqn:matrixA}
A=\left[\begin{array}{cc}
A_{11} & A_{12} \smallskip\\
0      & I
\end{array}\right].
\end{equation}
The elements in the matrices $A_{11}\in\mathbb{R}\sp{\mt_1\times\mt_1}$ and $A_{12}\in\mathbb{R}\sp{\mt_1\times\mt_2}$, with $\mt_1=p(2p+1)$ and \linebreak$\mt_2=(p+1)(2p+1)$, are obtained from derivatives of the interpolating polynomial $\tilde  u$ implied by the conditions in~\eqref{eq:leftcbc}.  The $\mt_2\times\mt_2$ identity in the lower-right block of $A$ is implied by the conditions in~\eqref{eq:leftbc} and~\eqref{eq:leftinterior}.  The matrix $A$ is nonsingular provided that the coefficient function $c_{11}(\xv)$ associated with the highest $x$-derivative in the differential operator~Q does not vanish (see Theorem~\ref{theorem:faceSolvability} discussed later in Section~\ref{sec:solvability}).  Algorithm~\ref{alg:formb} shows the construction of the right-hand side vector $\bv$ which follows similar steps to that used to build~$A$.  The solution of~\eqref{eq:leftSys} yields the coefficients $d_{\ip,\jp}$ of the interpolating polynomial, and in particular
\begin{equation} \label{eq:ghost}
U_{\ih,\tilde{j}}=d_{\ih,0},\qquad \ih=-p,\ldots,-1,
\end{equation}
which sets the values of $U_{\iv}$ in the $p$ ghost points corresponding to the boundary point $\tilde\xv$.

\input tex/FormAbAlg.tex

\subsubsection{LCBC method: Direct approach}

In the direct approach to the LCBC method, the matrix $A$ and vector $\bv$ in~\eqref{eq:leftSys} are constructed for each point on the boundary, and then the system is solved to determine ghost points following the assignments in~\eqref{eq:ghost} for example.  Algorithms~\ref{alg:formA} and~\ref{alg:formb} provide the steps for a point $\tilde\xv=(\tilde x,\tilde y)$ along the boundary $\tilde x=0$ for the case of a Dirichlet boundary condition, while the case of a Neumann boundary condition is described in Section~\ref{sec:LCBCneumann} below.  Points on the boundary near corners require special treatment, and this is discussed in Section~\ref{sec:LCBCcorner}.  

An important element of the direct approach, and the stencil approach discussed next, is an efficient calculation of the matrix $A$.  The main step in this calculation appears in line~8 of Algorithm~\ref{alg:formA}, which is independent of time $t$ and need only be performed once for a given problem.  This step involves applying repeated $y$-derivatives and powers of the operator $Q$ on the product of Lagrange basis functions $L_\ip$ and $L_\jp$, and then evaluating the result at a point $\tilde{\xv}$ on the boundary.  While this calculation can be carried out analytically, the form of $Q$ in~\eqref{eqn:OperatorQ} involving general coefficient functions, $c_{11}(\xv)$, $c_{12}(\xv)$, etc., makes this calculation increasingly messy as the order of accuracy determined by $p$ increases.  Also, it is desirable to avoid having to specify derivatives of the coefficient functions.  With these issues in mind, a more practical approach is described in Algorithm~\ref{alg:elementA} which computes suitable approximations of these elements, denoted by  $Z_{\ip,\jp}[\mu,\nu]$, in a particular column of~$A$ determined by given values of $\ip,\jp\in\{-p,\dots,p\}$ defining the basis functions.  The row entries are determined by the integers $\mu$ and $\nu$, and we note in advance that the algorithm only requires evaluations of the coefficient functions at points on the grid.

\input tex/elementAalg.tex

The first collection of steps in the algorithm results in the calculation of the grid function $V_{\hat{\iv}}[\nu+1,k]$ in line~16 defined by
\begin{equation}
V_{\hat{\iv}}[\nu,k]\eqdef(Q_{d,h})\sp\nu L_{\ip}\bigl(\ih\bigr)L_{\jp}\bigl(\jh\bigr),\qquad \nu=1,\ldots,p,
\end{equation}
where the indices $(\ip,\jp)$ are fixed and the order of accuracy of the approximation is $d=2k$, $k=1,\ldots,p+1-\nu$.  Note that the highest order of accuracy, given by $2(p+1-\nu)$, decreases as $\nu$ increases.  The calculation of $V_{\hat{\iv}}[\nu+1,k]$, determined by the function \textcolor{blue}{${\rm applyQh}$}, follows from the form of the discrete operator $Q_{d,h}$ given in~\eqref{eqn:Qhexp}.  Using the correction terms denoted by $W_{\hat{\iv}}\sp{(m,n)}$ in line~11, the function sets
\begin{equation}\label{eqn:applyQh}
\begin{array}{l}
\displaystyle{
V_{\hat{\iv}}[\nu+1,k]=c_{11}(\xv_{\hat{\iv}})D_{+x}D_{-x}\left\{V_{\hat{\iv}}[\nu,k]+\sum_{n=1}\sp{k-1}a_{n}\dx\sp{2n}W_{\hat{\iv}}\sp{(n,0)}[\nu,k-n]\right\}
}\smallskip\\
\qquad\qquad\qquad\displaystyle{
+2c_{12}(\xv_{\hat{\iv}})D_{0x}D_{0y}\left\{V_{\hat{\iv}}[\nu,k]+\sum_{n=1}\sp{k-1}\sum_{l=0}\sp{n}b_{l}b_{n-l}\dx\sp{2l}\dy\sp{2(n-l)}W_{\hat{\iv}}\sp{(l,n-l)}[\nu,k-n]\right\}
}\smallskip\\
\qquad\qquad\qquad\;\;\displaystyle{
+c_{22}(\xv_{\hat{\iv}})D_{+y}D_{-y}\left\{V_{\hat{\iv}}[\nu,k]+\sum_{n=1}\sp{k-1}a_{n}\dy\sp{2n}W_{\hat{\iv}}\sp{(0,n)}[\nu,k-n]\right\}
}\smallskip\\
\qquad\qquad\qquad\;\;\;\;\displaystyle{
+c_{1}(\xv_{\hat{\iv}})D_{0x}\left\{V_{\hat{\iv}}[\nu,k]+\sum_{n=1}\sp{k-1}b_{n}\dx\sp{2n}W_{\hat{\iv}}\sp{(n,0)}[\nu,k-n]\right\}
}\smallskip\\
\qquad\qquad\qquad\;\;\;\;\;\;\displaystyle{
+c_{2}(\xv_{\hat{\iv}})D_{0y}\left\{V_{\hat{\iv}}[\nu,k]+\sum_{n=1}\sp{k-1}b_{n}\dy\sp{2n}W_{\hat{\iv}}\sp{(0,n)}[\nu,k-n]\right\} + c_0(\xv_{\hat{\iv}})V_{\hat{\iv}}[\nu,k],
}
\end{array}
\end{equation}
where the coefficients $(a_n,b_n)$ are given in~\eqref{eq:QhexpCoeffs}.  The domain for the local index $\hat{\iv}$, denoted by $\hat{\Omega}_h[\nu,k]$, for each calculation is defined by
\begin{equation}\label{eq:indexDomain}
\hat{\Omega}_h[\nu,k]\eqdef[-w_x,w_x]\times[-w_y,w_y],\qquad w_x=p-(\nu+k-1),\quad w_y=w_x+p,
\end{equation}
and this gives the minimum stencil width required for the subsequent calculation of the discrete $y$-derivatives of $V_{\hat{\iv}}[\nu,k]$ performed in the second collection of steps starting at line~20.  Here, the main step involves the function \textcolor{blue}{${\rm applyDy}$} in line~25 which computes the odd/even derivative pair $Z_{\ip,\jp}[\mu-1,\nu]$ and $Z_{\ip,\jp}[\mu,\nu]$ using standard centered finite differences in the $y$-direction to order of accuracy $d=2k=2(p+1-\nu)$.

The elements of the right-hand side vector $\bv$ in~\eqref{eq:leftSys} are specified by Algorithm~\ref{alg:formb} for the case of a Dirichlet boundary along $\tilde x=0$.  The difficult step appears in line~5 and it involves the calculation of successive $y$-derivatives of $R_{\ell,\nu}(\tilde y,t)$ defined in~\eqref{eq:Rdef}.  The calculation of $R_{\ell,\nu}(y,t)$, in turn, requires powers of the operator $Q$ applied to the forcing function $f(\xv,t)$.
As before, we use a practical approach in which the various derivatives, both in space and time, are performed approximately to appropriate orders of accuracy.  At present we have considered only a spatial discretization in the semi-discrete model in~\eqref{eqn:SemiDiscreteEqn} and so we assume the time derivatives in $R_{\ell,\nu}(\tilde y,t)$ are exact for now.  In terms of the spatial approximations, a key step involves applying powers of the discrete operator $Q_{d,h}$ onto $f(\xv,t)$ evaluated at grid points about $\tilde\xv$, and this can be done efficiently following steps similar to those described in Algorithm~\ref{alg:elementA}.  Discrete $y$-derivatives are then applied to the result, again following the previous algorithm.  The principal details involve the approximations of $R_{\ell,\nu}(\tilde y,t)$ and these are given in Algorithm~\ref{alg:elementB} for the time-dependent cases $q=1,2$ (with straightforward simplifications for the steady case $q=0$). For Algorithm~\ref{alg:elementB}, we redefine the domain for the local index $\hat{\iv}$ in \eqref{eq:indexDomain} such that $w_x = p - (\nu + k)$.

\input tex/elementBalg.tex

It is worth noting that for the time-dependent cases, the elements of $\bv$ must be calculated at each time step.  Also, the approximation of $\p_y\sp\mu R_{\ell,\nu}(\tilde y,t)$ uses values of $R_{\ell,\hat{j}}[\nu,t]$ about $\tilde y$, computed in Algorithm~\ref{alg:elementB}, and these can be used by the approximations at neighboring values along the boundary.  This observation suggests a possible savings in computational cost that is explored with the stencil approach discussed next.

\subsubsection{LCBC method: Stencil approach}

The aim of the stencil approach is to manipulate the linear system in~\eqref{eq:leftSys} so that the values in the ghost points in~\eqref{eq:ghost} corresponding to a point $\tilde\xv$ on the boundary can be computed using the stencil formula
\begin{align}
  U_{\ih,\tilde{j}} = & \sum_{\nu = 1}^p\sum_{j = \tilde{j}-p}^{\tilde{j}+p}\alpha_{\ih,\tilde j}^{(\nu,j)}R_{\ell,j}[\nu,t] + \sum_{i = 0}^{p}\sum_{j = \tilde{j}-p}^{\tilde{j}+p}\beta_{\ih,\tilde j}^{(i,j)}U_{i,j}(t), \qquad \ih = -p, \dots, -1, \label{eqn:StencilRep}
\end{align}
where $\alpha_{\ih,\tilde j}^{(\nu,j)}$ and $\beta_{\ih,\tilde j}^{(i,j)}$ are coefficients belonging to the left boundary centered at $\xv_{0,\jb}$.  A central point is that the coefficients in~\eqref{eqn:StencilRep} do not depend on time $t$ and can be computed from the matrix~$A$ in~\eqref{eqn:matrixA}.  Thus, the values in the ghost points can be computed efficiently via a fixed linear combination of the relevant time-dependent data, assuming $q=1,2$, given by $R_{\ell,j}[\nu,t]$ and the grid data given by $U_{i,j}(t)$.  This grid data includes values at interior points close to the boundary for $i=1,\ldots,p$ and Dirichlet boundary data, $U_{0,j}(t)=g_\ell(y_j,t)$.  Note that Algorithm~\ref{alg:elementB} computes $R_{\ell,\jh}[\nu,t]$ for values of the local index $\jh$ about $\jb$, but the range of the $y$-index can be extended readily to cover the whole left boundary (sufficiently separated from the corners).

To compute the coefficients in~\eqref{eqn:StencilRep}, we consider the linear system in~\eqref{eq:leftSys} in the form
\begin{equation} \label{eq:blocksys}
\left[\begin{array}{cc}
A_{11} & A_{12} \smallskip\\
0      & I
\end{array}\right]
\left[\begin{array}{c}
\dv_1 \smallskip\\
\dv_2
\end{array}\right]=
\left[\begin{array}{c}
D_y\Rv(t) \smallskip\\
\Uv(t)
\end{array}\right],
\end{equation}
where $\dv=\bigl[\dv_1,\dv_2]\sp{T}$ holds the coefficients of the interpolating polynomial, $\Rv(t)\in\mathbb{R}\sp{\mt_1}$ is a vector containing $R_{\ell,j}[\nu,t]$, $\Uv(t)\in\mathbb{R}\sp{\mt_2}$ is a vector containing $U_{i,j}(t)$, and $D_y\in\mathbb{R}\sp{\mt_1\times\mt_1}$ is the matrix operator representing the discrete $y$-derivatives of $R_{\ell,j}[\nu,t]$.  We are mainly interested in the elements of $\dv_1$ which give the ghost values in~\eqref{eq:ghost}.  The lower set of $\mt_2$ equations in~\eqref{eq:blocksys} implies $\dv_2=\Uv(t)$ so that the upper set of $\mt_1$ equations becomes
\begin{equation} \label{eq:stencilEqn}
A_{11}\dv_1=D_y\Rv(t)-A_{12}\Uv(t).
\end{equation}
Let $C_\alpha\in\mathbb{R}\sp{\mt_1\times\mt_1}$ and $C_\beta\in\mathbb{R}\sp{\mt_1\times\mt_2}$ solve the matrix systems
\begin{equation} \label{eq:stencilCoeffs}
A_{11}C_\alpha=D_y,\qquad A_{11}C_\beta=-A_{12},
\end{equation}
so that~\eqref{eq:stencilEqn} reduces to
\begin{equation} \label{eqn:StencilMatrixRep}
\dv_1=C_\alpha\Rv(t)+C_\beta\Uv(t).
\end{equation}
The sets of coefficients, $\bigl\{\alpha_{\ih,\tilde j}^{(\nu,j)}\bigr\}$ and $\bigl\{\beta_{\ih,\tilde j}^{(i,j)}\bigr\}$, in the stencil formula in~\eqref{eqn:StencilRep} are given by the elements along selected rows of $C_\alpha$ and $C_\beta$, respectively, corresponding to the desired ghost values in $\dv_1$.  We note also that the linear systems in~\eqref{eq:stencilCoeffs} are dense but not very large, e.g.~$A_{11}$ is $21\times21$ for $p=3$, and they can be solved readily using standard linear algebra software.

\subsection{Neumann boundary}\label{sec:LCBCneumann}

We now turn our attention to the case of a Neumann boundary.  Let us again treat the left boundary, $x=0$ with $y\in[0,1]$, and assume a primary boundary condition given by
\begin{equation}
  \p_xu(\xv,t) = g_{\ell}(y,t), \qquad \xv\in \p\Omega_\ell,
\label{eq:nbc}
\end{equation}
for a fixed time $t\in[0,T]$.  We consider the interpolating polynomial $\tilde u(\xv)$ in~\eqref{eqn:InterpolatingPolynomial} centered about a point $\tilde{\xv}\in\p\Omega_\ell$, and then specify its coefficients by enforcing the constraints
\bse
\label{eq:leftconditionsN}
\bat
    \tilde{u}\bigl(\ih\dx,\tilde{y} + \jh\dy\bigr)  &= U_{\ih,\tilde{j}+\jh}(t),                    \qquad& \ih = 0,\dots,p, & \quad \jh = -p,\dots, p,\label{eq:leftinteriorN} \\
    \py^{\mu}\p_xQ^{\nu}\tilde{u}(0,\tilde{y}) &= \py^{\mu}S_{\ell,\nu}(\tilde{y},t), \qquad& \nu=0,\ldots,p-1, & \quad \mu = 0,\dots,2p, \label{eq:leftcbcN}
\eat
\ese
where
\begin{equation} \label{eq:Sdef}
S_{\ell,\nu}(y,t)\eqdef
\begin{cases}
-\p_xQ^{\nu - 1}f(0,y), & q = 0,  \\[5pt]
\pt^{q\nu}g_\ell(y,t) - \p_x\Psi_{\nu}f(0,y,t), & q=1,2.
\end{cases}
\end{equation}
The interpolating polynomial $\tilde u$ is set equal to the data $U_\iv(t)$ in~\eqref{eq:leftinteriorN} at $2p+1$ points on the boundary and at $p(2p+1)$ points in the interior adjacent to the boundary.  The compatibility conditions in~\eqref{eq:leftcbcN} include tangential derivatives of the Neumann boundary condition in~\eqref{eq:nbc} for the case $\nu=0$.  Note that we use only $p-1$ CBCs here since the boundary condition involves a normal derivative, whereas $p$ CBCs are used for the case of a Dirichlet boundary.  As before, the constraints in~\eqref{eq:leftconditionsN} imply $\mt=(2p+1)\sp2$ linear equations for the $\mt$ coefficients in~$\tilde{u}$ for each point $(0,\tilde y)\in \p\tilde{\Omega}_{\ell,h}$ given in~\eqref{eq:discreteBoundary}, and these equations can be written in a matrix form following~\eqref{eq:leftSys} with a coefficient matrix~$A$ and right-hand side vector~$\bv$ generated in Algorithms~\ref{alg:formAN} and~\ref{alg:formbN}, respectively.  The solution of the linear system gives~$\dv$, the coefficients of~$\tilde u(\xv)$, and thus the ghost point values in~\eqref{eq:ghost} associated with the point $\tilde{\xv}$ on the Neumann boundary.

\input tex/FormAbNAlg.tex

The elements in line~8 of Algorithm~\ref{alg:formAN}, involving repeated applications of the operator $Q$ on the $x$-derivative of the Lagrange basis functions followed by successive $y$-derivatives, are obtained approximately.  The details of these approximations are similar to those given in Algorithm~\ref{alg:elementA} for the Dirichlet case.  The main difference for the case of a Neumann boundary is the addition of a discrete $x$-derivative applied to $V_{\hat{\iv}}[\nu,k]$, to an appropriate order of accuracy, prior to applications of the discrete $y$-derivatives.  The $y$-derivatives of $S_{\ell,\nu}(\tilde y,t)$ in line~5 of Algorithm~\ref{alg:formbN} are also computed approximately.

Instead of solving the linear system in~\eqref{eq:leftSys} at each time step (for $q=1,2$), we prefer the stencil approach.  Here, the stencil formula for the ghost points associated with a point $\xv_{0,\jb}$ along the left boundary becomes
\begin{align}
  U_{\ih,\tilde{j}} = & \sum_{\nu = 0}^{p-1}\sum_{j = \tilde{j}-p}^{\tilde{j}+p}\bar\alpha_{\ih,\tilde j}^{(\nu,j)}S_{\ell,j}[\nu,t] + \sum_{i = 0}^{p}\sum_{j = \tilde{j}-p}^{\tilde{j}+p}\bar\beta_{\ih,\tilde j}^{(i,j)}U_{i,j}(t), \qquad \ih = -p, \dots, -1, \label{eqn:StencilRepN}
\end{align}
where $S_{\ell,j}[\nu,t]\approx S_{\ell,\nu}(y_j,t)$ for $q>0$ and $\nu=0,\ldots,p-1$.
The coefficients in~\eqref{eqn:StencilRepN} for each boundary point are obtained by manipulating the system
\begin{equation} \label{eq:blocksysN}
\left[\begin{array}{cc}
\bar{A}_{11} & \bar{A}_{12} \smallskip\\
0      & I
\end{array}\right]
\left[\begin{array}{c}
\dv_1 \smallskip\\
\dv_2
\end{array}\right]=
\left[\begin{array}{c}
\bar{D}_y\Sv(t) \smallskip\\
\Uv(t)
\end{array}\right].
\end{equation}
For the Neumann case, the block matrices $\bar A_{11}\in\mathbb{R}\sp{\mt_1\times\mt_1}$ and $\bar A_{12}\in\mathbb{R}\sp{\mt_1\times\mt_2}$ are obtained from the approximations of the elements in line~8 of Algorithm~\ref{alg:formAN}.  The vectors $\Sv(t)\in\mathbb{R}\sp{\mt_1}$ and $\Uv(t)\in\mathbb{R}\sp{\mt_2}$ contain $S_{\ell,j}[\nu,t]$ and $U_{i,j}(t)$, respectively, and $\bar D_y\in\mathbb{R}\sp{\mt_1\times\mt_1}$ contains the coefficients associated with the discrete $y$-derivatives of $S_{\ell,j}[\nu,t]$, $\nu=0,\ldots,p-1$.  Let $\bar C_{\bar\alpha}\in\mathbb{R}\sp{\mt_1\times\mt_1}$ and $\bar C_{\bar\beta}\in\mathbb{R}\sp{\mt_1\times\mt_2}$ solve the systems
\begin{equation} \label{eq:stencilCoeffsN}
\bar A_{11}\bar C_{\bar\alpha}=\bar D_y,\qquad \bar A_{11}\bar C_{\bar\beta}=-\bar A_{12},
\end{equation}
so that
\begin{equation} \label{eqn:StencilMatrixRepN}
\dv_1=\bar C_{\bar\alpha}\Sv(t)+\bar C_{\bar\beta}\Uv(t).
\end{equation}
As before, the sets of coefficients of the stencil formula in~\eqref{eqn:StencilRepN} are given by the elements along selected rows of $\bar C_{\bar\alpha}$ and $\bar C_{\bar\beta}$ corresponding to the desired ghost values in $\dv_1$.

\subsection{LCBC conditions at a corner}\label{sec:LCBCcorner}

As a representative case involving the conditions at a corner, let us consider the bottom-left corner, $\tilde x=(0,0)$, where two Dirichlet boundaries meet.  The cases of a Neumann-Neumann corner and a Dirichlet-Neumann corner are discussed in~\ref{sec:LCBCcornerGeneral}.  The physical (primary) boundary are taken to be
\bse
\begin{align}
u(\xv,t)&=g_\ell(y,t), \qquad\xv\in\partial\Omega_\ell, \\
u(\xv,t)&=g_b(x,t), \qquad \xv\in\partial\Omega_b,
\end{align}
\ese
for a some fixed time $t$.  We start by specifying the interpolating polynomial $\tilde u(\xv)$ at known interior data given by
\bse \label{eq:cornerconditions}
\begin{equation}
\tilde{u}\bigl(\ih\dx,\jh\dy\bigr) = U_{\ih,\jh}(t),\qquad \ih = 1,\dots,p, \quad \jh = 1,\dots, p.
\label{eq:interiorC}
\end{equation}
Next, we apply tangential derivatives of the primary boundary conditions and compatibility conditions given by
\begin{equation}
\left.
\begin{array}{l}
\py^{\mu}\tilde{u}(0,0) = \py^{\mu}g_{\ell}(0,t) \smallskip\\
\px^{\mu}\tilde{u}(0,0) = \px^{\mu}g_{b}(0,t)
\end{array}\;\right\}
\qquad \mu\in{\cal M}_0,
\label{eq:bcC}
\end{equation}
and
\begin{equation}
\left.
\begin{array}{l}
\py^{\mu}Q^{\nu}\tilde{u}(0,0) = \py^{\mu}R_{\ell,\nu}(0,t) \smallskip\\
\px^{\mu}Q^{\nu}\tilde{u}(0,0) = \px^{\mu}R_{b,\nu}(0,t)
\end{array}\;\right\}
\qquad \nu=1,\ldots,p, \quad \mu\in{\cal M}_\nu,
\label{eq:cbcCC}
\end{equation}
\ese
respectively, where $R_{\ell,\nu}(y,t)$ is defined in~\eqref{eq:Rdef} and $R_{b,\nu}(x,t)$ is defined by
\begin{equation} \label{eq:Rbdef}
R_{b,\nu}(x,t)\eqdef
\begin{cases}
-Q^{\nu - 1}f(x,0), & q = 0,  \\[5pt]
\pt^{q\nu}g_b(x,t) - \Psi_{\nu}f(x,0,t), & q=1,2.
\end{cases}
\end{equation}
The sets ${\cal M}_\nu$, $\nu=0,\ldots,p$, chosen to eliminate redundant constraints, are given by
\begin{equation}
  {\cal M}_\nu=\left\{
  \begin{array}{ll}
  0,1,2,3,\ldots,2p-1,2p, &\hbox{if $\nu=0$, with an average for $\mu=0$,}\smallskip\\
  \quad 1,2,3,4,\ldots,2p-1,2p, &\hbox{if $\nu=1$, with an average for $\mu=2$,}\smallskip\\
  \quad 1,\quad 3,4,5,\ldots,2p-1,2p, &\hbox{if $\nu=2$, with an average for $\mu=4$,}\smallskip\\
  \quad\quad\quad\vdots       &\quad\quad\vdots \smallskip\\
  \quad 1,\quad 3,\quad 5,\ldots,2p-1,2p, &\hbox{if $\nu=p$, with an average for $\mu=2p$.}
  \end{array}
  \right.  \label{eq:MnuDirichletDirichletCorner}
\end{equation}
Note that there is one value for $\mu$ in each set ${\cal M}_\nu$ where the pairs in~\eqref{eq:bcC} and~\eqref{eq:cbcCC} are averaged to resolve linearly dependent constraints (and to balance the constraints on the left and bottom boundaries).  The weights for the averages are $\dy\sp\mu$ and $\dx\sp\mu$ for the CBCs arising from the left and bottom boundaries, respectively, to balance the tangential derivatives taken in the $y$ and $x$ directions (following the previous discussion in Section~\ref{sec:wave2d}).  The $\mt=(2p+1)\sp2$ constraints in~\eqref{eq:cornerconditions} for the $\mt$ coefficients $d_{\ip,\jp}$ of $\tilde u(\xv)$ are a generalization of the approach for the bottom-left corner given in~\eqref{eq:cornerconditions4} for the wave equation.  Ghost points near the corner can be obtained from the solution of the linear system implied by~\eqref{eq:cornerconditions} following a direct approach, or these ghost points can be written in terms of the stencil formula
\begin{align}
  U_{\ih,\jh} = & \sum_{\nu = 0}^p\sum_{j = -p}^{p}\tilde\alpha_{\ih,\jh}^{(\nu,j)}R_{\ell,j}[\nu,t] + \sum_{\nu = 0}^{p}\sum_{i =-p}^{p}\tilde\beta_{\ih,\jh}^{(\nu,i)}R_{b,i}[\nu,t] + \sum_{i = 1}^{p}\sum_{j = 1}^{p}\tilde\gamma_{\ih,\jh}^{(i,j)}U_{i,j}(t), \qquad \hat{\iv}\in\hat\Omega_c,   \label{eqn:StencilRepC}
\end{align}
where
\begin{equation}
\hat\Omega_c\eqdef\bigl\{\,\hat{\iv}=(\ih,\jh)\,\vert\,-p\le(\ih,\jh)<p\;\backslash\; 1\le(\ih,\jh)<p\bigr\}
\end{equation}
defines the set of local indices for the ghost-point values in~\eqref{eqn:StencilRepC}.  The time-dependent data $R_{\ell,j}[\nu,t]$ and $R_{b,i}[\nu,t]$ in~\eqref{eqn:StencilRepC}, assuming $q=1,2$, are discrete approximations of $R_{\ell,\nu}(j\dy,t)$ and $R_{b,\nu}(i\dx,t)$, respectively, for $\nu=1,\ldots,p$.  The boundary conditions are specified in~\eqref{eqn:StencilRepC} by setting
\bse
\begin{align}
R_{\ell,j}[0,t]&=g_\ell(j\dy,t), \qquad j=-p,\ldots,p, \\
R_{b,i}[0,t]&=g_b(i\dx,t), \qquad i=-p,\ldots,p,
\end{align}
\ese
similar to previous specifications.  The coefficients in the stencil formula are derived from the $\mt\times\mt$ linear system implied by~\eqref{eq:cornerconditions} following the analysis described for the Dirichlet and Neumann boundaries.

Our choice for the constraints in~\eqref{eq:cornerconditions} is guided by the case when $Q$ in~\eqref{eqn:OperatorQ} is the Laplacian operator.  For this case, the constraints are linearly independent.  For the more general operator $Q$ with variable coefficients, the constraints remain linearly independent provided $c_{11}(\xv)>0$, $c_{22}(\xv)>0$ and $\max\left\{\abs{c_{12}(\xv)}/\sqrt{c_{11}(\xv)c_{22}(\xv)}\right\}$ is small, with $\xv$ evaluated at the corner.  Should these conditions be violated, the $\mt\times\mt$ matrix $A$ implied by~\eqref{eq:cornerconditions} may become singular or badly conditioned.  For example, if $c_{11}>0$, $c_{22}>0$ and $c_{12}$ are constants, and if $c_{1} = c_{2} = c_{0} = 0$, then the determinant of~$A$ for the case $p=1$ ($d=2$) has the form
\begin{equation*}
  \det(A) = - D \dx \dy (c_{11} + c_{22})(c_{11}c_{22} - 4c_{12}^2), \qquad D=\hbox{constant}>0.
\end{equation*}
Thus, $A$ becomes singular when $\abs{c_{12}} = \sqrt{c_{11}c_{22}}/2$.  Another case for which $A$ is rank deficient occurs when $c_{11} = c_{22} = 1$, $c_{12} = 1/2$, $c_{1} = c_{2} = c_{0} = 0$ and $\dx = \dy$, and for any value of~$p$.

As noted earlier, we are motivated by high-order accurate discretizations of the model problem in~\eqref{subsec:MainEqn}.  For many problems of interest, this model problem is obtained by an orthogonal, or near-orthogonal, mapping of a PDE in physical space involving the Laplacian operator.  The resulting mapped problem, of the form in section~\eqref{subsec:MainEqn}, would have $\vert c_{12}(\xv)\vert$ small relative to $c_{11}(\xv)$ and $c_{22}(\xv)$ resulting in a nonsingular matrix $A$ implied by the constraints in~\eqref{eq:cornerconditions} for a Dirichlet-Dirichlet corner.  The matrices for the Neumann-Neumann and Dirichlet-Neumann corners, discussed in the appendix, are also nonsingular under these conditions.
\subsection{Conditioning of LCBC matrices}
We are interested in the sensitivity of the computational procedures involving LCBC matrices to data perturbations and round-off errors.  We first note that the rows of the LCBC matrices for both sides and corners vary in magnitude significantly.  Therefore, we choose to use the \texttt{equilibrate} function in MATLAB to optimize the conditioning of the LCBC matrices by permuting the rows and then applying row and column scalings. The diagonal entries of each permuted and scaled matrix are all one in magnitude while the off-diagonal elements are less than or equal to one. This decreases the number of pivots needed for numerical stability in Gaussian elimination (see~\cite{equilibrateDoc}). We shall then adopt the scaling proposed by \texttt{equilibrate} to generate the condition numbers for the numerical results presented below.  

Let $\kappa(A)=\| A \| \| A^{-1}\| $ be the 2-norm condition number of an LCBC matrix~$A$ generated for a point on a boundary or a corner, and let $\kappa_{\text{max}}$ be the maximum of $\kappa(A)$ for all points along the four boundaries or for the four corners.  We compute $\kappa_{\text{max}}$ for the following test cases:
\begin{enumerate}
	\item We set $\Omega = (0,1)^2$, $Q = \Delta$ and take Dirichlet boundary conditions at $x = 0$ and~$1$, $y\in[0,1]$, and Neumann boundary conditions at $y=0$ and~$1$, $x\in[0,1]$. This case corresponds to the numerical example presented later in Section~\ref{NumRes:example1}.
	\item We set $\Omega = \left\{\xv=(\rho\cos\theta,\rho\sin\theta)\;|\;1<\rho<2,\,0<\theta<\pi/2\right\}$ and $Q = \Delta$, and take Dirichlet boundary conditions at $\rho=1$, $\theta\in[0,\pi/2]$ and at $\theta = 0$, $\rho\in[1,2]$ and Neumann boundary conditions at $\rho=2$, $\theta\in[0,\pi/2]$ and at $\theta=\pi/2$, $\rho\in[1,2]$. This case corresponds to the numerical example presented in Section~\ref{NumRes:example3}.
	\item We set $\Omega$ to be the wavy channel domain defined in~\eqref{eqn:flagregion} for the problem described in Section~\ref{NumRes:HeatFlowWavyChannel}. The spacial operator is $Qu = \Dc \Delta u - \vv\cdot \grad{u}  + \gamma u$,  where $\Dc = 0.2$, $\vv = (0.5,0.3)$ and $\gamma = 1$, and we use Dirichlet conditions for $\xv\in\partial\Omega$.
\end{enumerate}

Figure~\ref{fig:conditionNumbers} shows the computed values of $\kappa_{\text{max}}$ for the side and corner LCBC matrices as functions of the grid spacing $h$ used for each test case and for orders of accuracy $d=2$, $4$ and~$6$.  For the first two test cases, we observe that $\kappa_{\text{max}}$ for the side and corner LCBC matrices are approximately equal, and they are nearly independent of the grid spacing $h$, although $\kappa_{\text{max}}$ shows a moderate increase with increasing values of the order of approximation $d$ for $Q_{d,h}$.  The behavior of $\kappa_{\text{max}}$ for the third test case shows more variation with $h$ for the corner matrices when $d=4$ and~$6$.  The plots in Figure~\ref{fig:conditionNumbersSpecial} show more details of the behavior of $\kappa_{\text{max}}$ as a function of $h$ for $d=4$.  There is more variation for coarser grids, which we attribute to the affects of the lower-order terms in $Q$, while the behavior of $\kappa_{\text{max}}$ appears to settle to a constant value as the grid spacing gets smaller (and affects of the lower-order terms become less significant).
{
	\newcommand{\figWidth}{5.7cm}
	\newcommand{\trimfig}[2]{\trimw{#1}{#2}{0.}{0.}{.0}{.0}}  
	\newcommand{\trimfigb}[2]{\trimw{#1}{#2}{.0}{.05}{.1}{.1}} 
	\begin{figure}[H]
		\begin{center}
			\begin{tikzpicture}[scale=1]
				\useasboundingbox (0,0) rectangle (16,6.3);  
				\draw (-1.17,-0.5) node[anchor=south west] {\trimfig{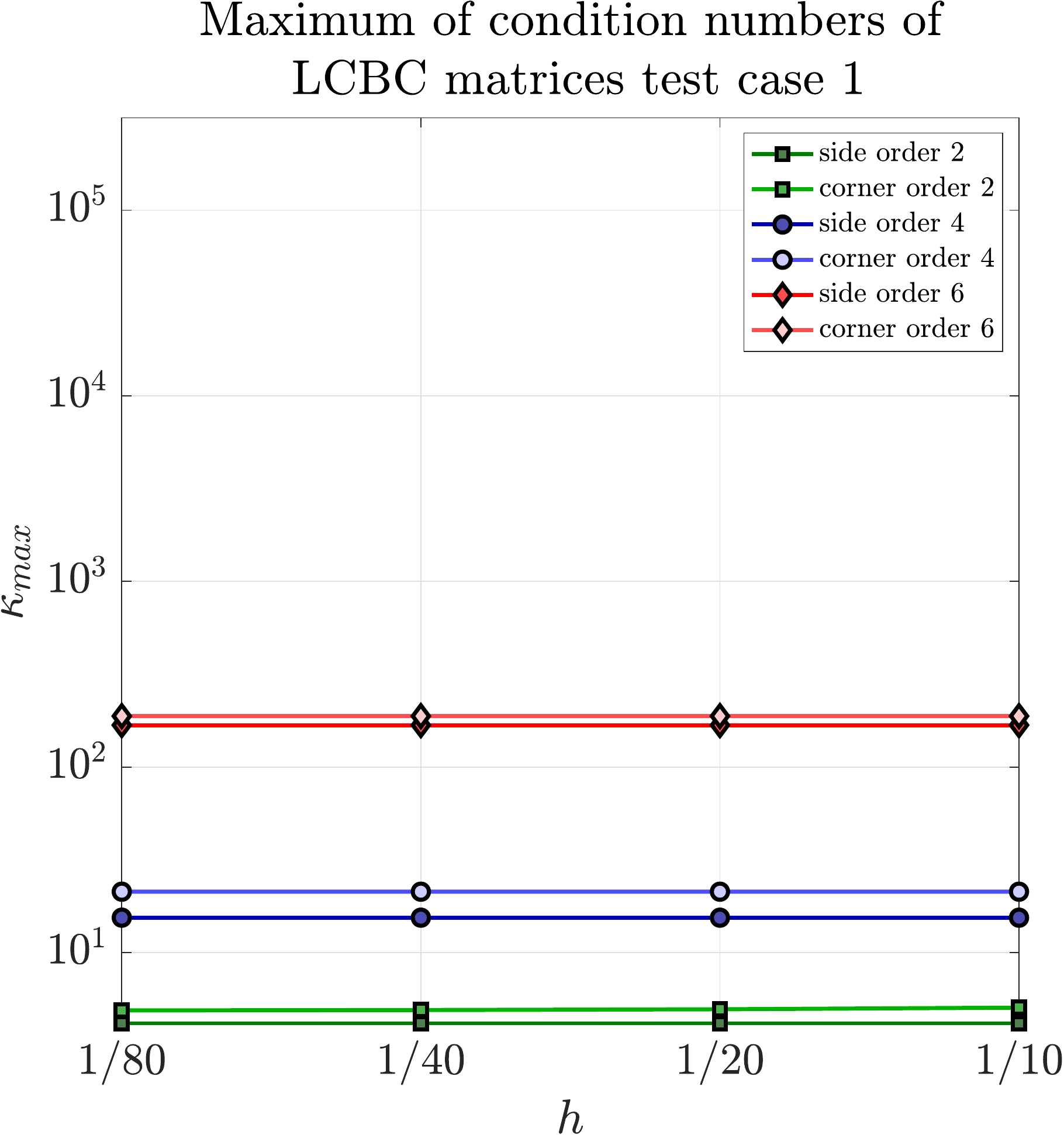}{\figWidth}};
				\draw (4.52,-0.5) node[anchor=south west] {\trimfig{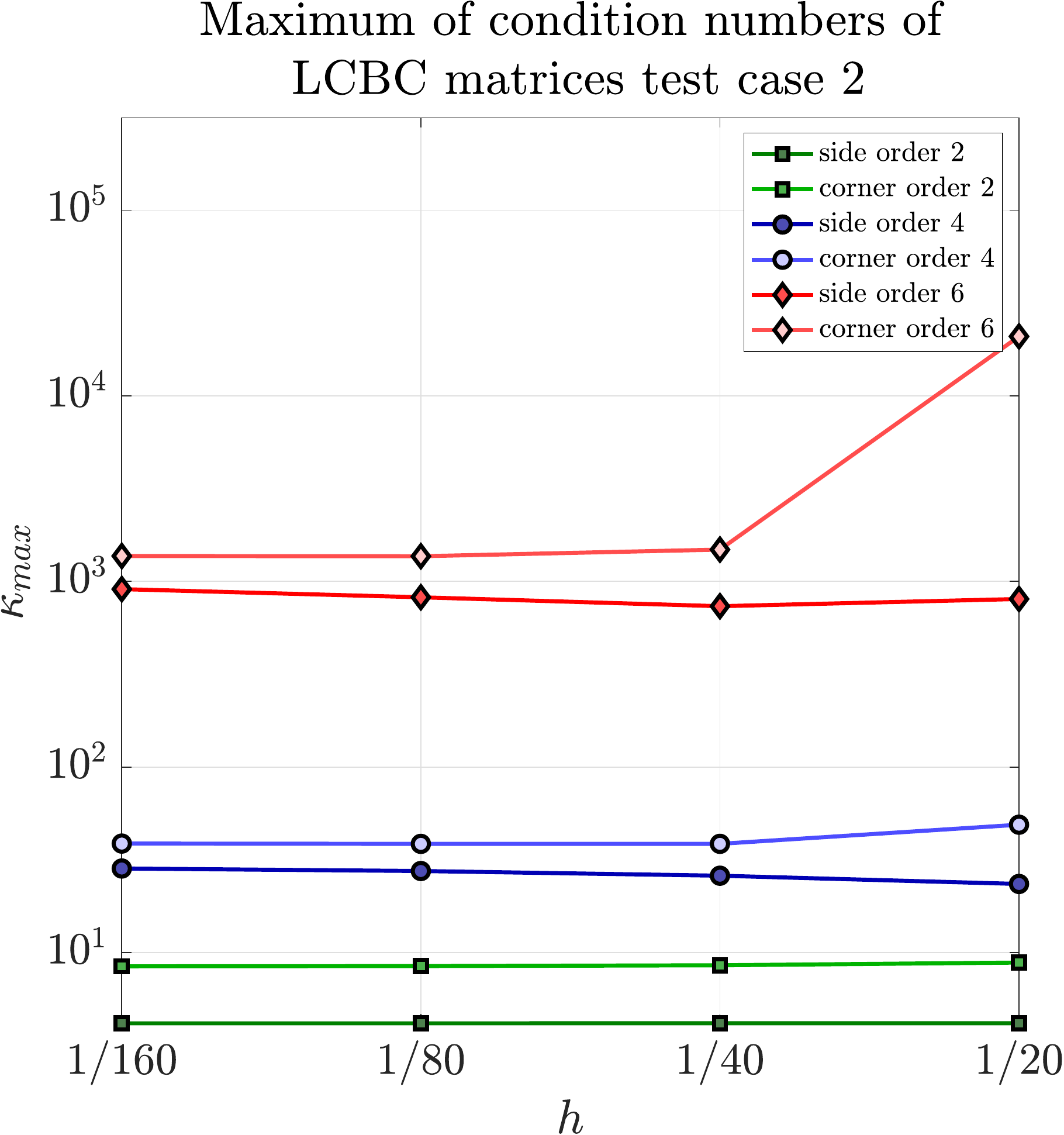}{\figWidth}};
				\draw (10.21,-0.5) node[anchor=south west] {\trimfig{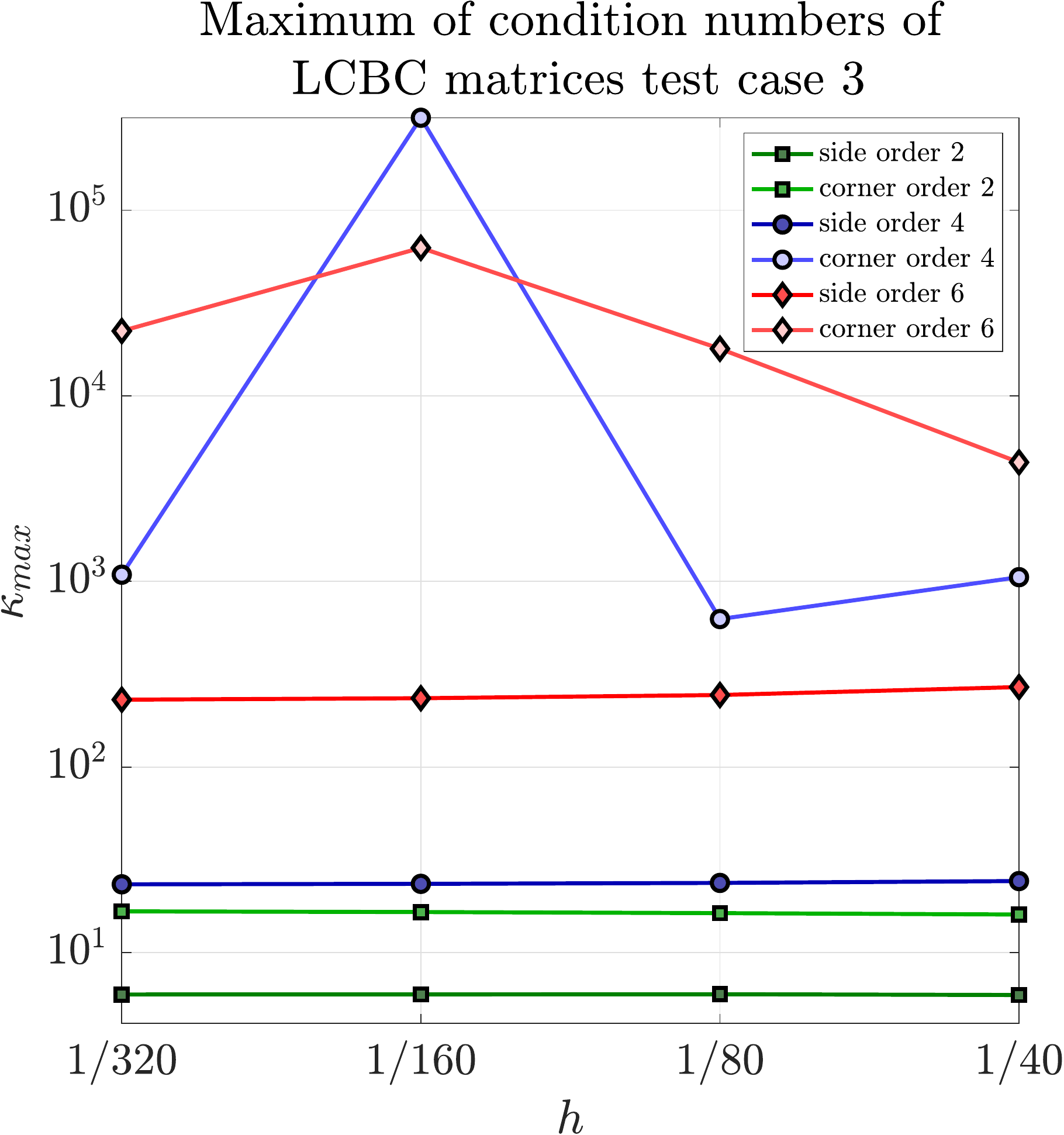}{\figWidth}};
			\end{tikzpicture}
		\end{center}
		\caption{Condition numbers of LCBC matrices as functions of the grid spacing $h$ for test cases 1, 2 and 3.}
		\label{fig:conditionNumbers}
	\end{figure}
}
{
	\newcommand{\figWidth}{7cm}
	\newcommand{\trimfig}[2]{\trimw{#1}{#2}{0.}{0.}{.0}{.0}}  
	\newcommand{\trimfigb}[2]{\trimw{#1}{#2}{.0}{.05}{.1}{.1}} 
	\begin{figure}[H]
		\begin{center}
			\begin{tikzpicture}[scale=1]
				\useasboundingbox (0,0) rectangle (13,6.5);  
				\draw (-0.9,-.5) node[anchor=south west] {\trimfig{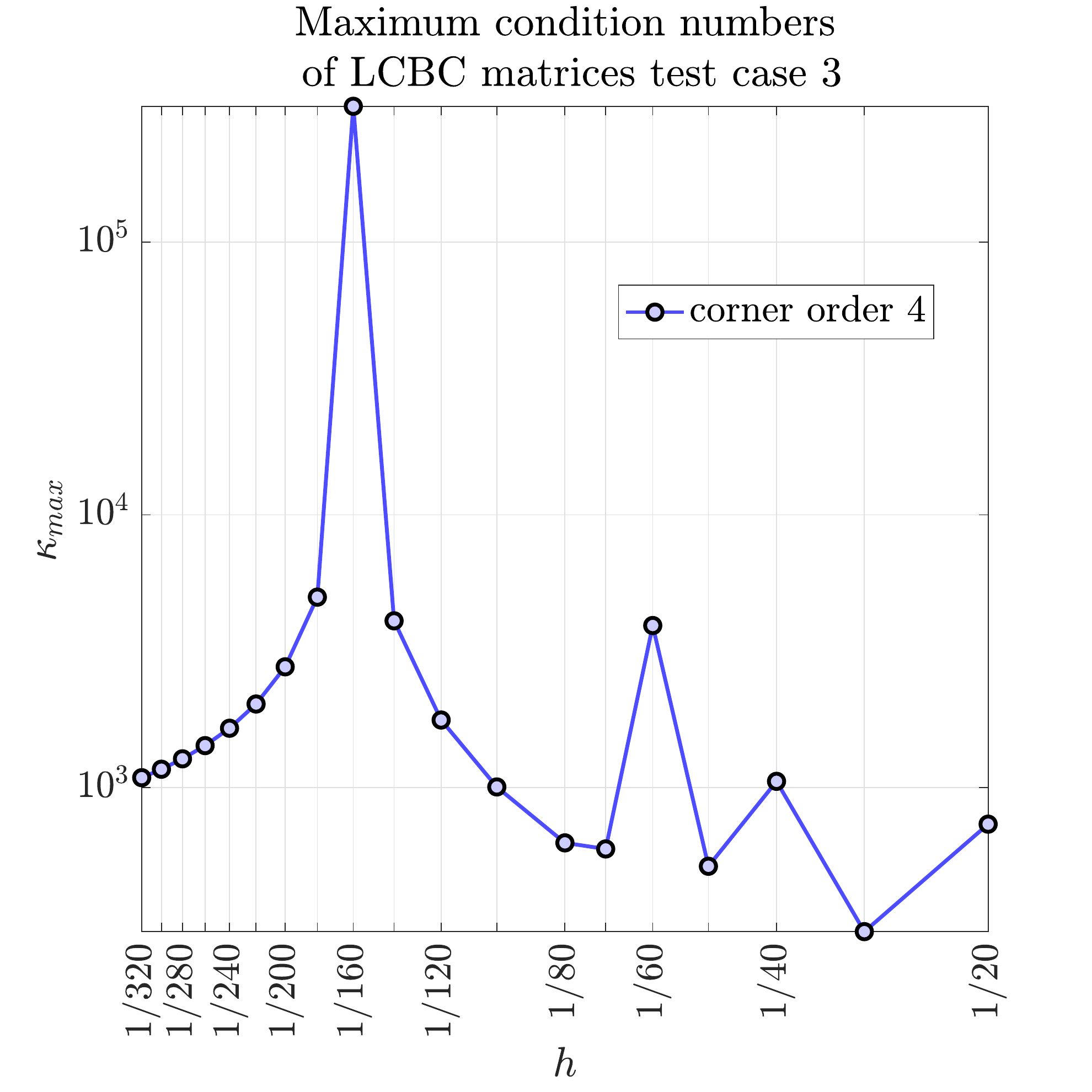}{\figWidth}};
				\draw (5.7,-.5) node[anchor=south west] {\trimfig{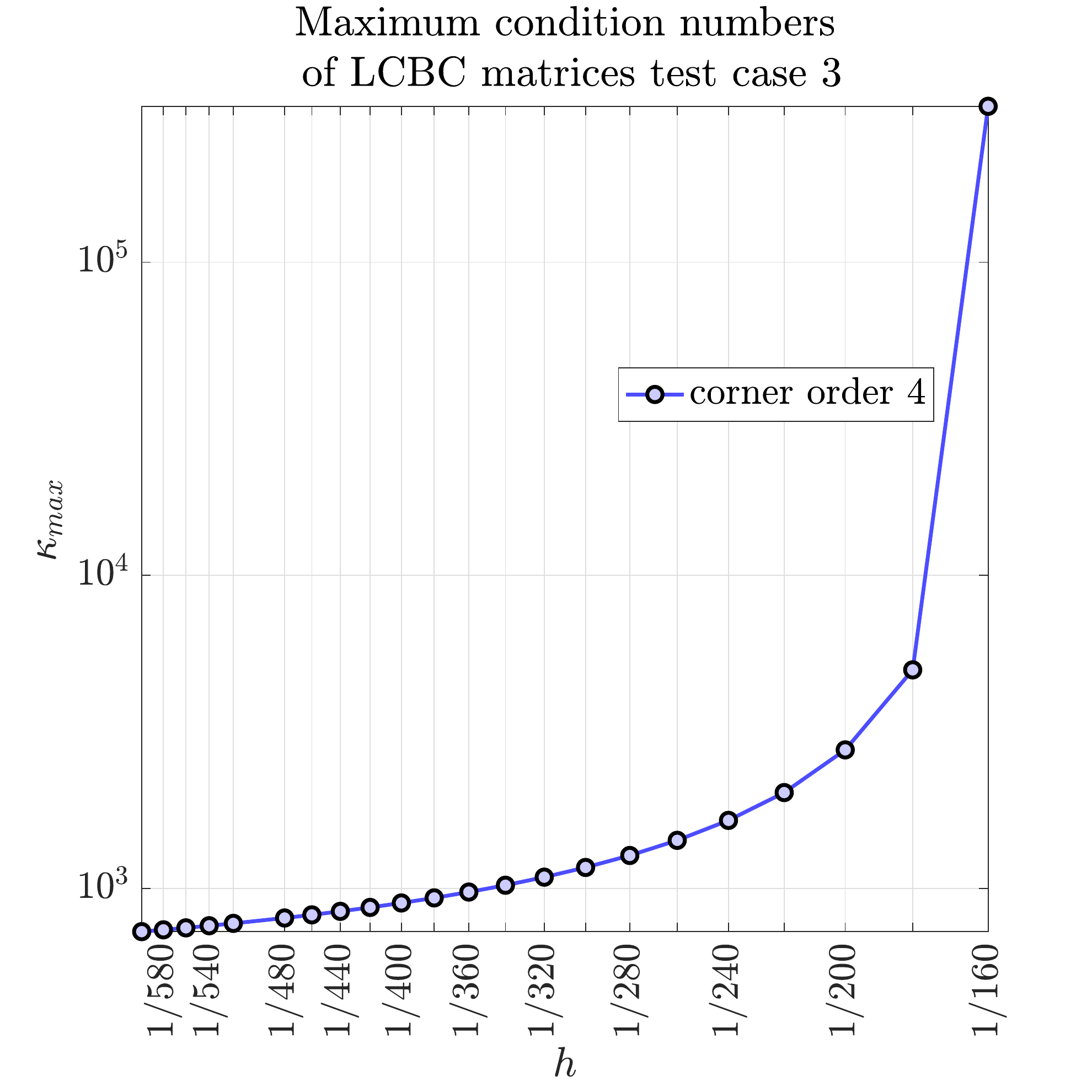}{\figWidth}};
			\end{tikzpicture}
		\end{center}
		\caption{Test case 3: $\kappa_{\text{max}}$ versus $h$ for corner matrices at order $d=4$.}
		\label{fig:conditionNumbersSpecial}
	\end{figure}
}

%% file: tex/FormAbAlg.tex
\begin{algorithm}[H]
	\algFontSize 
	\caption{~Construct the coefficient matrix $A$ for a Dirichlet boundary.}
	\begin{algorithmic}[1]
		\State $r = 0$;
		\For{$\nu = 1, \dots, p$}
		\For{$\mu = 0,\dots, 2p$}
		\State $r = r + 1$; 
		\For{$\ip = -p, \dots,p$}
		\For{$\jp = -p,\dots, p$}
		\State $c = (2p+1)(\ip + p) + \jp + p + 1$; 
		\State $A(r,c) = \py^{\mu}Q^{\nu}\bigl.L_{\ip}\bigl((x-\tilde{x})/\dx\bigr)L_{\jp}\bigl((y-\tilde{y})/\dy\bigr)\bigr\vert_{\xv=\tilde\xv}$\,; \Comment Elements of $A$ from~\eqref{eq:leftcbc}
		\EndFor
		\EndFor
		\EndFor
		\EndFor
		\For{$\ih = 0, \dots, p$}
		\For{$\jh = -p, \dots, p$}
		\State $r = r+1$; 
		\State $A(r,r) = 1$; \Comment Elements of $A$ from~\eqref{eq:leftbc} and~\eqref{eq:leftinterior}
		\EndFor
		\EndFor
	\end{algorithmic}
	\label{alg:formA}
\end{algorithm}

\begin{algorithm}[H]
	\algFontSize 
	\caption{~Construct the right-hand side vector $\bv$ for a Dirichlet boundary.}
	\begin{algorithmic}[1]
		\State $r = 0$; 
		\For{$\nu = 1, \dots,p$}
		\For{$\mu = 0,\dots, 2p$}
		\State $r = r + 1$; 
		\State $b(r) = \py^{\mu}R_{\ell,\nu}(\tilde y,t)$; \Comment Elements of $\bv$ from~\eqref{eq:leftcbc}
		\EndFor
		\EndFor
		\For{$\jh = -p, \dots, p$}
		\State $r = r + 1$; 
		\State $b(r) = g_\ell\bigl(\tilde y+\jh\dy,t\bigr)$;  \Comment Elements of $\bv$ from~\eqref{eq:leftbc}
		\EndFor
		\For{$\ih = 1, \dots, p$}
		\For{$\jh = -p, \dots, p$}
		\State $r = r + 1$; 
		\State $b(r) = U_{\ih,\tilde{j}+\jh}(t)$;  \Comment Elements of $\bv$ from~\eqref{eq:leftinterior}
		\EndFor
		\EndFor
	\end{algorithmic}
\label{alg:formb}	
\end{algorithm}

%% file: tex/elementAalg.tex
%

\begin{algorithm}[H]
	\algFontSize 
	\caption{~Compute $Z_{\ip,\jp}[\mu,\nu]\approx\py^{\mu}Q^{\nu}\bigl.L_{\ip}\bigl((x-\tilde{x})/\dx\bigr)L_{\jp}\bigl((y-\tilde{y})/\dy\bigr)\bigr\vert_{\xv=\tilde\xv}$.}
	\begin{algorithmic}[1]
		\For{$k=1,\ldots,p$}
		\For{$\hat{\iv}\in\hat{\Omega}_h[0,k]$}  \Comment Initialize $V_{\hat{\iv}}[0,k]=L_\ip(\ih)L_\jp(\jh)$
		\State $V_{\hat{\iv}}[0,k] = L_\ip(\ih)L_\jp(\jh)$;
		\EndFor
		\EndFor
		\For{$\nu = 0, \dots, p-1$}
		\For{$k=1,\ldots,p-\nu$}
		\For{$l=1,\ldots,k-1$}
		\For{$\hat{\iv}\in\hat{\Omega}_h[\nu,k]$}  \Comment Compute corrections $W_{\hat{\iv}}\sp{(m,n)}[\nu,l]$ involving $V_{\hat{\iv}}[\nu,l]$
		\For{$m=0,\ldots,k-l$}
        \State $W_{\hat{\iv}}\sp{(m,(k-l)-m)}[\nu,l]=(D_{+x}D_{-x})\sp{m}(D_{+y}D_{-y})\sp{(k-l)-m}V_{\hat{\iv}}[\nu,l]$;
		\EndFor
		\EndFor
		\EndFor
		\For{$\hat{\iv}\in\hat{\Omega}_h[\nu+1,k]$} \Comment Compute $V_{\hat{\iv}}[\nu+1,k]=\bigl(Q_{2k,h}\bigr)V_{\hat{\iv}}[\nu,k]$
        \State $V_{\hat{\iv}}[\nu+1,k]=\textcolor{blue}{{\rm applyQh}}\bigl\{V_{\hat{\iv}}[\nu,k],W_{\hat{\iv}}\sp{(m,n)}[\nu,k-1],\ldots,W_{\hat{\iv}}\sp{(m,n)}[\nu,1]\bigr\}$;
		\EndFor
		\EndFor
		\EndFor
		\For{$\nu=1,\ldots,p$} \Comment Compute $Z_{i,j}[\mu,\nu]$ using $V_{\hat{\iv}}[\nu,k]$, $k=1,2,\ldots,p+1-\nu$
		\State $k=p+1-\nu$;
		\State $Z_{\ip,\jp}[0,\nu]=V_{0,0}[\nu,k]$
		\For{$l=1\ldots,p$}
		\State $\mu=2l$;
		\State $\bigl\{Z_{\ip,\jp}[\mu-1,\nu],Z_{\ip,\jp}[\mu,\nu]\bigr\}=\textcolor{blue}{{\rm applyDy}}\bigl\{V_{\hat{\iv}}[\nu,1],\ldots,V_{\hat{\iv}}[\nu,k]\bigr\}$;
		\EndFor
		\EndFor

	\end{algorithmic}
	\label{alg:elementA}
\end{algorithm}

%% file: tex/elementBalg.tex
\begin{algorithm}[H]
	\algFontSize 
	\caption{~Compute $R_{\ell,\hat{j}}[\nu,t]\approx R_{\ell,\nu}\bigl(\tilde{y}+\hat{j}\dy,t\bigr)$ for $q>0$}
	\begin{algorithmic}[1]
		\For{$\nu=1,\ldots,p$}
		\For{$\hat{j}\in[-p,p]$}  \Comment Initialize $R_{\ell,\hat{j}}[\nu,t]=\p_t\sp{q\nu}g_{\ell}\bigl(\tilde y+\hat{j}\dy,t\bigr)$
		\State $R_{\ell,\hat{j}}[\nu,t]=\textcolor{blue}{{\rm applyDt}}\bigl\{g_{\ell,\hat{j}}(t),q\nu\bigr\}$;
		\EndFor
		\EndFor
		\For{$n=0,\ldots,p-1$}
		\For{$k=1,\ldots,p$}
		\For{$\hat{\iv}\in\hat{\Omega}_h[0,k]$}  \Comment Initialize $F_{\hat{\iv}}[0,k,t]=\p_t\sp{qn}f(\tilde\xv+\xv_{\hat{\iv}},t)$
		\State $F_{\hat{\iv}}[0,k,t] = \textcolor{blue}{{\rm applyDt}}\bigl\{f_{\hat{\iv}}(t),qn\bigr\}$;
		\EndFor
		\EndFor
		\For{$\bar\nu = 0, \dots, p-n-2$}
		\For{$k=1,\ldots,p-\bar\nu-1$}
		\For{$l=1,\ldots,k-1$}
		\For{$\hat{\iv}\in\hat{\Omega}_h[\bar\nu,k]$}  \Comment Compute corrections $W_{\hat{\iv}}\sp{(m,n)}[\bar\nu,l,t]$ involving $F_{\hat{\iv}}[\bar\nu,l,t]$
		\For{$m=0,\ldots,k-l$}
        \State $W_{\hat{\iv}}\sp{(m,(k-l)-m)}[\bar\nu,l,t]=(D_{+x}D_{-x})\sp{m}(D_{+y}D_{-y})\sp{(k-l)-m}F_{\hat{\iv}}[\bar\nu,l,t]$;
		\EndFor
		\EndFor
		\EndFor
		\For{$\hat{\iv}\in\hat{\Omega}_h[\bar\nu+1,k]$} \Comment Compute $F_{\hat{\iv}}[\bar\nu+1,k,t]=\bigl(Q_{2k,h}\bigr)F_{\hat{\iv}}[\bar\nu,k,t]$
        \State $F_{\hat{\iv}}[\bar\nu+1,k,t]=\textcolor{blue}{{\rm applyQh}}\bigl\{F_{\hat{\iv}}[\bar\nu,k,t],W_{\hat{\iv}}\sp{(m,n)}[\bar\nu,k-1,t],\ldots,W_{\hat{\iv}}\sp{(m,n)}[\bar\nu,1,t]\bigr\}$;
		\EndFor
		\EndFor
		\EndFor
		\For{$\nu = n+1, \dots, p$}
		\State $\bar\nu=\nu-n-1$;
		\State $k={\rm min}\{p-\bar\nu\,,\,p\}$;
		\For{$\hat{j}\in[-p,p]$}   \Comment Update $R_{\ell,\hat{j}}[\tilde\nu,t]$
		\State $R_{\ell,\hat{j}}[\nu,t]=R_{\ell,\hat{j}}[\nu,t]-F_{0,\hat{j}}[\bar\nu,k,t]$;
		\EndFor
		\EndFor
		\EndFor
	\end{algorithmic}
	\label{alg:elementB}
\end{algorithm}

%% file: tex/FormAbNAlg.tex
\begin{algorithm}[H]
	\algFontSize 
	\caption{~Construct the coefficient matrix $A$ for a Neumann boundary.}
	\begin{algorithmic}[1]
		\State $r = 0$;
		\For{$\nu = 0, \dots, p-1$}
		\For{$\mu = 0,\dots, 2p$}
		\State $r = r + 1$; 
		\For{$\ip = -p, \dots,p$}
		\For{$\jp = -p,\dots, p$}
		\State $c = (2p+1)(\ip + p) + \jp + p + 1$; 
		\State $A(r,c) = \py^{\mu}\p_xQ^{\nu}\bigl.L_{\ip}\bigl((x-\tilde{x})/\dx\bigr)L_{\jp}\bigl((y-\tilde{y})/\dy\bigr)\bigr\vert_{\xv=\tilde\xv}$\,; \Comment Elements of $A$ from~\eqref{eq:leftcbcN}
		\EndFor
		\EndFor
		\EndFor
		\EndFor
		\For{$i = 0, \dots, p$}
		\For{$j = -p, \dots, p$}
		\State $r = r+1$; 
		\State $A(r,r) = 1$; \Comment Elements of $A$ from~\eqref{eq:leftinteriorN}
		\EndFor
		\EndFor
	\end{algorithmic}
	\label{alg:formAN}
\end{algorithm}

\begin{algorithm}[H]
	\algFontSize 
	\caption{~Construct the right-hand side vector $\bv$ for a Neumann boundary.}
	\begin{algorithmic}[1]
		\State $r = 0$; 
		\For{$\nu = 0, \dots,p-1$}
		\For{$\mu = 0,\dots, 2p$}
		\State $r = r + 1$; 
		\State $b(r) = \py^{\mu}S_{\ell,\nu}(\tilde y,t)$; \Comment Elements of $\bv$ from~\eqref{eq:leftcbcN}
		\EndFor
		\EndFor
		\For{$\ih = 0, \dots, p$}
		\For{$\jh = -p, \dots, p$}
		\State $r = r + 1$; 
		\State $b(r) = U_{\ih,\tilde{j}+\jh}(t)$;  \Comment Elements of $\bv$ from~\eqref{eq:leftinteriorN}
		\EndFor
		\EndFor
	\end{algorithmic}
\label{alg:formbN}	
\end{algorithm}

%% file: tex/analysis.tex
\section{Analysis of the LCBC approach} \label{sec:analysis}

In this section, we provide some results of an analysis of the LCBC approach.  In particular, we consider the solvability of the matrix systems associated with the constraints implied by the LCBC method for points along a grid side and at a grid corner.  We then consider symmetry properties of the discrete approximations generated by the LCBC method for the case when the PDE involves the Laplacian operator.  Finally, we examine the stability of explicit time-stepping schemes for the wave equation with numerical boundary conditions given by the LCBC approach.

\subsection{Solvability of the LCBC matrix systems}  \label{sec:solvability}

We first consider conditions required for the LCBC matrix systems to be non-singular.  This is done for the case of a constant-coefficient operator $Q$ given by
\ba
 &  Q  = c_{11} \p_x^2 + 2 c_{12} \p_x\p_y  + c_{22} \p_y^2  + c_1 \p_x  + c_2 \p_y  + c_0 . \label{eq:QconstantCoeff}
\ea
For this operator, we have the following result:
\begin{theorem}[Solvability on a face] 
\label{theorem:faceSolvability}
The matrix resulting from the order $2p=2,4,6$ LCBC constraints for the constant-coefficient operator
$Q$ in~\eqref{eq:QconstantCoeff} with a Dirichlet or Neumann boundary condition on a grid face is non-singular
provided $c_{11} >0$ and $\dx$ is sufficiently small (left or right face) 
or  $c_{22}>0$ and $\dy$ is sufficiently small (bottom or top face). 
If $c_1=0$ (left face) or $c_2=0$ (right face), then the matrix is non-singular for any~$\dx$ and~$\dy$.
\end{theorem}
\begin{proof}
Let us focus on the left boundary, while similar arguments hold for the other boundaries.  The form of the LCBC matrix system is given in~\eqref{eq:leftSys}, where the coefficient matrix $A$ is obtained from the constraints in~\eqref{eq:leftconditions} for a Dirichlet boundary and in~\eqref{eq:leftconditionsN} for a Neumann boundary.  For either a Dirichlet or Neumann boundary, the determinant of $A$, for order of accuracy $2p=2,4,6$, has the form
\begin{equation}
\det(A)=K_p G_p(\xi),\qquad \xi={c_1\dx\over c_{11}},\qquad p=1,2,3,
\label{eq:sideDetCond}
\end{equation}
where $K_p$ is a non-zero constant depending on $\dx$, $\dy$ and $c_{11}$, and where $G_p(\xi)$ is a polynomial satisfying $G_p(0)=1$.  For the Dirichlet case, the polynomials are given by
\[
\begin{array}{l}
\displaystyle{
G_1(\xi)=\left(1-{\xi\over2}\right)\sp3,
} \medskip\\
\displaystyle{
G_2(\xi)=\left(1-{3\xi\over2}+{\xi\sp2\over2}-{\xi\sp3\over18}\right)\sp5,
} \medskip\\
\displaystyle{
G_3(\xi)=\left(1-3\xi+{11\xi\sp2\over4}-{1691\xi\sp3\over1440}+{121\xi\sp4\over480}-{11\xi\sp5\over400}+{\xi\sp6\over800}\right)\sp7,
}
\end{array}
\]
and for the Neumann case, we have
\[
\begin{array}{l}
\displaystyle{
G_1(\xi)=1,
} \medskip\\
\displaystyle{
G_2(\xi)=\left(1-{2\xi\over9}\right)\sp5,
} \medskip\\
\displaystyle{
G_3(\xi)=\left(1-{23\xi\over30}+{11\xi\sp2\over75}-{\xi\sp3\over100}\right)\sp7.
}
\end{array}
\]
The result of the theorem follows from the form of the determinant of $A$ in~\eqref{eq:sideDetCond}.  As expected, the lower order terms in~\eqref{eq:QconstantCoeff} become less important for the solvability of the system as the grid spacings tend to zero.

\noindent{\Huge$\square$}
\end{proof}

The solvability conditions at a corner are more complicated.  For this case, we focus on the constant-coefficient operator in~\eqref{eq:QconstantCoeff} with the coefficients of the lower-order terms set to zero, i.e.~$c_0=c_1=c_2=0$, and define the dimensionless parameters
\[
\gamma={c_{12}\over\sqrt{c_{11}c_{22}}},\qquad \sigma=\sqrt{{c_{11}/\dx\sp2\over c_{22}/\dy\sp2}},
\]
assuming $c_{11}>0$ and $c_{22}>0$.  Recall that when choosing the corner compatibility conditions we assumed that $\vert c_{12}\vert$ is small compared to $c_{11}$ and $c_{22}$, and this now corresponds to $\vert\gamma\vert$ small.  The following theorem describes the solvability of the LCBC matrix systems for the Dirichlet-Dirichlet (D-D), Neumann-Neumann (N-N) and Dirichlet-Neumann (D-N) corners.
\begin{theorem}[Solvability at a corner] 
\label{theorem:cornerSolvability}
	The matrices resulting from the LCBC constraints at D-D, N-N and D-N corners for the constant-coefficient operator $Q$ in~\eqref{eq:QconstantCoeff} with $c_{11}>0$, $c_{22}>0$, and $c_0=c_1=c_2=0$ are nonsingular provided any of the following
	conditions hold:
	\begin{enumerate}
	  \item $\gamma=0$ ($c_{12}=0$), for orders $2p=2,4,6$.  
	  \item $\vert\gamma\vert$ is sufficiently small, for orders $2p=2,4$.
	  \item $\gamma<0$ and $\vert\gamma\vert$ is sufficiently small, for order $2p=6$.
	  \item $\gamma>0$ and $(\sigma+1/\sigma)\gamma$ is sufficiently small, for order $2p=6$. 
	\end{enumerate}  
\end{theorem}
\begin{proof}
We consider the corner where the left and bottom boundaries meet, while similar arguments hold for the other corners.  The form of the LCBC matrix system is given in~\eqref{eq:leftSys}, where the coefficient matrix $A$ is obtained from the constraints in~\eqref{eq:cornerconditions}, \eqref{eq:NNcornerconditions} and~\eqref{eq:DNcornerconditions} for D-D, N-N and D-N corners, respectively.  For all three cases, the determinant of $A$, for order of accuracy $2p=2,4,6$, has the form
\begin{equation}
\det(A)=K_p H_p(\gamma) F_p(\gamma,\sigma),\qquad p=1,2,3,
\label{eq:cornerDetCond}
\end{equation}
where $K_p$ is a non-zero constant depending on $\dx$, $\dy$, $c_{11}$ and $c_{22}$, $H_p(\gamma)$ is a polynomial satisfying $H_p(0)=1$, and $F_p(\gamma,\sigma)$ is a polynomial in $\gamma$ with coefficients that depend on $\sigma$.  For a D-D corner, we have
\[
\begin{array}{l}
\displaystyle{
H_1(\gamma)=1-4\gamma\sp2,
} \medskip\\
\displaystyle{
H_2(\gamma)=\left(1-4\gamma\sp2\right)\sp2\left(1-28\gamma\sp2+208\gamma\sp4-256\gamma\sp6\right),
} \medskip\\
\displaystyle{
H_3(\gamma)=\left(1-4\gamma\sp2\right)\sp4\left(1-12\gamma\sp2+16\gamma\sp4\right)\sp2\left(1-104\gamma\sp2+3984\gamma\sp4-68480\gamma\sp6\right.
} \medskip\\
\quad\qquad\qquad\qquad\qquad\qquad\qquad\qquad\qquad\qquad\displaystyle{
\left.+509440\gamma\sp8-1278976\gamma\sp{10}+921600\gamma\sp{12}\right),
}
\end{array}
\]
and
\[
\begin{array}{l}
\displaystyle{
F_1(\gamma,\sigma)=1,
} \medskip\\
\displaystyle{
F_2(\gamma,\sigma)=3\left(\sigma+{1\over\sigma}\right)-4\gamma,
} \medskip\\
\displaystyle{
F_3(\gamma,\sigma)=7200\left(\sigma\sp3+\sigma+{1\over\sigma}+{1\over\sigma\sp3}\right)-\gamma\left[3960\left(\sigma\sp4+{1\over\sigma\sp4}\right)+28070\left(\sigma\sp2+{1\over\sigma\sp2}\right)+26620\right]
} \medskip\\
\qquad\qquad\qquad\displaystyle{
+\gamma\sp2\left[13423\left(\sigma\sp3+{1\over\sigma\sp3}\right)+39483\left(\sigma+{1\over\sigma}\right)\right]
-\gamma\sp3\left[14399\left(\sigma\sp2+{1\over\sigma\sp2}\right)+28798\right]
} \medskip\\
\qquad\qquad\qquad\qquad\displaystyle{
+\gamma\sp4\left[5940\left(\sigma+{1\over\sigma}\right)\right].
}
\end{array}
\]
For a N-N corner, we have
\begin{align*}
&H_1(\gamma) = 1, \\
&H_2(\gamma) = 1 - 4\gamma^2, \\
&H_3(\gamma) = (1-4\gamma^2)^2\left(1 - 28\gamma^2 + 208\gamma^4 - 256\gamma^6\right),
\end{align*}
and 
\begin{align*}
& F_1(\gamma,\sigma) = 1, \\
& F_2(\gamma,\sigma) = 27\left(\sigma + \frac{1}{\sigma}\right)-32\gamma,\\
& F_3(\gamma,\sigma) = 56250\left(\sigma^3 + \sigma + \f{1}{\sigma} + \f{1}{\sigma^3}\right) -\gamma\left[ 24750\left(\sigma^4+\f{1}{\sigma^4}\right) + 194625\left(\sigma^2 + \f{1}{\sigma^2}\right)+ 189750\right]\\
&\qquad \qquad\qquad  + \gamma^2\left[240310\left(\sigma + \f{1}{\sigma}\right) + 25450\left(\sigma^3 + \f{1}{\sigma^3}\right)\right] - \gamma^3\left[71564\left(\sigma^2 + \f{1}{\sigma^2}\right)+150040\right]\\
&\qquad \qquad \qquad \qquad + \gamma^4\left[25344\left(\sigma + \f{1}{\sigma}\right)\right].
\end{align*}
Finally, for a D-N corner, we have
\begin{align*}
	&H_1(\gamma) = 1, \\
	&H_2(\gamma) = 1-12\gamma^2 + 32\gamma^4, \\
	&H_3(\gamma) = (1  -4\gamma^2)^2\left(1 - 64 \gamma^2 + 1504\gamma^4 - 16128\gamma^6 + 80640\gamma^8 - 171008\gamma^{10} + 122880\gamma^{12}\right), 
\end{align*}
and 
\begin{align*}
	& F_1(\gamma,\sigma) = 1,\\
	& F_2(\gamma,\sigma) = 3\left(\sigma + \f{1}{\sigma}\right) - 4\gamma,\\
	& F_3(\gamma,\sigma) = 135000\left(\sigma^3 + \sigma + \f{1}{\sigma} + \f{1}{\sigma^3}\right) - 75\gamma\left(660\sigma^4 + 7353\sigma^2 + 6523 + 6418\f{1}{\sigma^{2}} + 1188\f{1}{\sigma^{4}}\right)\\
	& \qquad \qquad\qquad + 5\gamma^2\left(35985\sigma^3 + 158555\sigma + 126287\f{1}{\sigma} + 54738\f{1}{\sigma^{3}}\right)\\
	& \qquad \qquad \qquad\qquad  -\gamma^3\left(215985\sigma^2 + 520784 + 239532\f{1}{\sigma^{2}}\right) + 35640\gamma^4\left(3\sigma + 2\f{1}{\sigma}\right). 
\end{align*}
Note that when $\gamma<0$, the functions $F_p$ are always positive and bounded away from zero.  The result of the theorem follows from the form of the determinant of $A$ in~\eqref{eq:cornerDetCond}. 

\noindent{\Huge$\square$}
\end{proof}

We note that a good quality grid usually aims to have $\sigma \approx 1$.  One way to see this is to note that if $c_{11} \ll c_{22}$ then there could be boundary layers near $x=0$ or $x=1$, which would require a small value for $\dx$ to resolve the solution there.  We also note that for order $2p=6$ when $\gamma>0$ ($c_{12}>0$), we require not just $\gamma$ to be small but also $\gamma\sigma$ and $\gamma/\sigma$ to be small. Thus the corner LCBC matrix could be
poorly conditioned if $\sigma$ becomes large or small when $c_{12}>0$.  This could occur, for example,
if one only refined the grid in the $x$-direction.

\subsection{Symmetry properties of the LCBC conditions} \label{sec:symmetry}

The next two theorems concern symmetry properties of the numerical boundary conditions generated by the LCBC method for a boundary face and corner.  These symmetry conditions pertain to the case when $Q$ is the Laplacian operator and the domain is represented by a Cartesian grid.  The first theorem considers the symmetry for a boundary face.

\begin{theorem}[Symmetry on a face] When applied to the operator $Q=\Delta$ on a Cartesian grid, the LCBC approach on a face, at any order $2p=2,4,6,\ldots$, results in numerical boundary conditions with odd symmetry for the case of homogeneous Dirichlet boundary conditions and with even symmetry for the case of homogeneous Neumann boundary conditions, for example, 
\begin{subequations}
\label{eq:symmtryCond}
\bat
 &  U_{i-\hat{i},j} = - U_{i+\hat{i},j}, \quad&& \hat{i}=1,\ldots,p ,  \qquad&& \text{Dirichlet BC at $i=0$ or $i=N_x$} ,\\
 &  U_{i-\hat{i},j} =  \;\; U_{i+\hat{i},j}, \quad&& \hat{i}=1,\ldots,p ,  \qquad&& \text{Neumann BC at $i=0$ or $i=N_x$} .
\eat
\end{subequations}
\end{theorem}
\begin{proof}

\mni
First consider the case of a homogeneous Dirichlet boundary condition on the left side, $i=0$, away from the corner. 
Without loss of generality we may take $\xt=0$ and $\yt=0$, and then the polynomial interpolant $\ut$ can be written
as
\ba
   \ut(x,y) = \sum_{n=0}^{2p} \sum_{m=0}^{2p} a_{n,m} \, x^n \, y^m, \qquad p=1,2,\ldots.
   \label{eq:uTildeMonomial}
\ea
We wish to show that $\ut(x,0)$ is an odd function in $x$, so that
\ba
   \ut(x,0) = a_{1,0} x + a_{3,0} x^3 + \ldots a_{2p-1,0} x^{2p-1}, \label{eq:DirichletOddPoly}
\ea
for then we have $\ut(-x,0)= - \ut(x,0)$ and the desired result follows.  The CBCs in~\eqref{eq:leftcbc} reduce to
\ba
   \p_y^\mu \Delta^\nu \ut(\zerov) = 0 , \qquad \nu=0,\ldots,p, \quad \mu=0,\ldots,2p,  \quad {\rm CBC}[\mu,\nu],
\label{eq:DirichletCBCs}
\ea
where the case $\nu=0$ follows since $U_{0,j}=0$ from the homogeneous boundary condition.  For the purposes of the proof, we have labeled the conditions in~\eqref{eq:DirichletCBCs} as ${\rm CBC}[\mu,\nu]$.  We will show that~\eqref{eq:DirichletCBCs} implies that all even $x$-derivatives of $\ut$ at $\xv=\zerov$ are zero, 
\ba
   \p_x^{2\nu} \ut(\zerov) = 0 ,  \qquad~\nu=0,\ldots,p, \label{eq:DirichletCBCevenDerivs}
\ea
which implies~\eqref{eq:DirichletOddPoly}.  The conditions in~\eqref{eq:DirichletCBCevenDerivs} can be shown as follows.  We have $\p_y^\mu \ut(\zerov) = 0$, 
for $\mu=0,1,\ldots$, since the Dirichlet conditions are homogeneous and since $\ut$ is a polynomial of finite degree.  Then, from ${\rm CBC}[0,1]$, we see that~\eqref{eq:DirichletCBCevenDerivs} holds for $\nu=1$ since 
\ba
    \p_x^2 \ut(\zerov) = - \p_y^2 \ut(\zerov) = 0 , 
\ea
and from ${\rm CBC}[\mu,1]$ we also find 
\ba
   \p_y^\mu \p_x^2 \ut(\zerov) = - \p_y^{\mu+2} \ut(\zerov)  = 0 , \qquad \mu=0,1,\ldots. 
\ea
Now from ${\rm CBC}[0,2]$, we find that~\eqref{eq:DirichletCBCevenDerivs} holds for $\nu=2$, since
\ba
   \p_x^4 \ut(\zerov) = (- 2\p_x^2\p_y^2 - \p_y^4)\ut(\zerov) = 0,
\ea
and from ${\rm CBC}[\mu,2]$ we also find 
\ba 
 \p_y^\mu \p_x^4 \ut(\zerov) = 0 , \qquad \mu=0,1,\ldots. 
\ea
The process can be repeated to show~\eqref{eq:DirichletCBCevenDerivs}.

\mni
The argument is similar for the case of a homogeneous Neumann boundary condition except that in this case
it can be shown that all odd $x$-derivatives are zero, 
\ba
   \p_x^{2\nu+1} \ut(\zerov) = 0 ,  \qquad~\nu=0,\ldots,p, \label{eq:NeumannCBCoddDerivs}
\ea
so that $\ut(-x,0)=\ut(x,0)$. 

\noindent{\Huge $\square$}
\end{proof}

We now consider the symmetry at a corner.  For this case, note that the LCBC conditions are used to obtain values in ghost points in the corner of the extended grid and also at nearby ghost points belonging to the adjacent faces, see Figure~\ref{fig:grid2d} for the case $p=2$ for example.

\begin{theorem}[Symmetry at a corner]
When applied to the operator $Q=\Delta$ on a Cartesian grid, the LCBC approach applied at any corner and at any order $2p=2,4,6,\ldots$, results in numerical boundary conditions on the adjacent faces with odd symmetry for the case of homogeneous Dirichlet boundary conditions and with even symmetry for the case of homogeneous Neumann boundary conditions.  At a left boundary, for example, the symmetries are given in~\eqref{eq:symmtryCond}.  Values at the corner ghost points have even symmetry for Dirichlet-Dirichlet (D-D) or Neumann-Neumann (N-N) corners and odd symmetry for Dirichlet-Neumann (D-N) corners.  At a bottom-left corner, for example, the values satisfy
\begin{subequations}
\bat
 &  U_{i-\hat{i},j-\hat{j}} =   \;\;U_{i+\hat{i},j+\hat{j}}, \quad&& \hat{i},\hat{j}=1,\ldots,p ,  \qquad&& \text{D-D or N-N corners}, \\
 &  U_{i-\hat{i},j-\hat{j}} =  -U_{i+\hat{i},j+\hat{j}}, \quad&& \hat{i},\hat{j}=1,\ldots,p ,  \qquad&& \text{D-N corner}.
\eat
\end{subequations}
\end{theorem}
\begin{proof}
Consider the case of homogeneous Dirichlet boundary conditions on the left side, $i = 0$, and the bottom side, $j = 0$, so that we have a D-D corner at $\xv = (0,0)$ and grid index $\iv = (0,0)$. With $\tilde{u}$ given in \eqref{eq:uTildeMonomial}
we show that 
\begin{equation} \label{eq:Csym}
\tilde{u}(-x,y) = -\tilde{u}(x,y),\qquad \tilde{u}(x,-y) = -\tilde{u}(x,y), 
\end{equation}
and thus 
\begin{equation}
\tilde{u}(-x,-y) = \tilde{u}(x,y). 
\end{equation}
To show \eqref{eq:Csym}, we show
\begin{subequations}
\begin{align}
\px^{m_1}\py^{m_2}\tilde{u}(\zerov) = 0, \qquad m_1=2k, \quad m_2 = 0,1,\dots,2p,\label{eq:Csym1}\\
\px^{m_1}\py^{m_2}\tilde{u}(\zerov) = 0, \qquad m_1=0,1,\dots,2p,\quad  m_2 =2k,\label{eq:Csym2}
\end{align}
\end{subequations}
where $k = 0,1,\dots,p$.

Recall that $\ut$ satisfies the boundary conditions in~\eqref{eq:bcC} and the compatibility conditions in~\eqref{eq:cbcCC} with homogeneous boundary data, so that
\begin{equation}
  \left.
  \begin{array}{l}
    \py^{\mu}\Delta^{\nu}\tilde{u}(\zerov) = 0 \smallskip\\
    \px^{\mu}\Delta^{\nu}\tilde{u}(\zerov) = 0
  \end{array}\;\right\}
  \qquad \nu=0,1,\ldots,p, \quad \mu\in{\cal M}_\nu,
  \label{eq:cbcC}
\end{equation}
where $\cal M_{\nu}$ is defined in \eqref{eq:MnuDirichletDirichletCorner}.
Using mathematical induction, we find that \eqref{eq:cbcC} implies 
\begin{equation}
  \left.
  \begin{array}{l}
    \py^{\mu}\px^{2\nu}\tilde{u}(\zerov) =0 \smallskip\\
    \px^{\mu}\py^{2\nu}\tilde{u}(\zerov) = 0
  \end{array}\;\right\}
  \qquad \nu=0,1,\ldots,p, \quad \mu\in{\cal M}_\nu. 
  \label{eq:cbcChom}
\end{equation}
Set $m_1 = 2k$ for $k = 0,1,\dots, p$. The first set of conditions in \eqref{eq:cbcChom} implies that
\begin{equation}
	\px^{m_1}\py^{m_2}\tilde{u}(\zerov) = 0, \qquad \text{ for } m_2= 1,3,5,\dots,2k-1,2k,2k+1,\dots, 2p,
\end{equation} 
while the second set of conditions in \eqref{eq:cbcChom} gives 
\begin{equation}
	\px^{m_1}\py^{m_2}\tilde{u}(\zerov) = 0, \qquad \text{ for } m_2= 0,2,4,\dots, 2k. 
\end{equation}
Hence, for $m_1 = 2k$, we have 
\begin{equation}
	\px^{m_1}\py^{m_2}\tilde{u}(\zerov) = 0, \qquad \text{ for } m_2= 0,1,2,3,\dots, 2p,
\end{equation}
for any $k = 0,1,\dots, p$. The result in \eqref{eq:Csym2} follows using a symmetric argument.  Therefore, we have odd symmetry on the Dirichlet sides near the corner and even symmetry at the D-D corner.  The results for N-N and D-N corners follow using similar arguments.

\noindent{\Huge $\square$}
\end{proof}

\subsection{Stability of LCBC approximations for the wave equation}  \label{sec:stabilityWaveEquation}

We now consider the stability of an explicit modified equation (ME) time-stepping algorithm for the wave equation on a Cartesian grid using the LCBC approach at the boundary.  For the present model problem in~\eqref{eqn:MainEqn}, the standard wave equation corresponds to the case $q=2$ and $Q=c\sp2\Delta$ for a wave speed $c>0$, and the ME time-stepping schemes are given in~\eqref{eqn:ModifiedEqn2}, \eqref{eqn:ModifiedEqn4} and~\eqref{eqn:ModifiedEqn6} for orders of accuracy $2p=2$, $4$ and~$6$, respectively.  In~\cite{AnneJolyHuyTran2000} it was shown that an ME scheme for the wave equation in one space dimension is stable at any order of accuracy, $2p=2,4,6,\ldots$, under the condition $c\dt/\dx<1$, where $\dt$ is the time-step.  In two dimensions (or three dimensions), the time-step condition depends on whether selected terms are dropped to retain a stencil width of $2p+1$ or not.  For example, at sixth-order, the term $\Delta_{4,h}^2 U_\iv^n$ appears, and it has a term proportional to $\dx^4 (\Dpx\Dmx)^4 U_\iv^n$ which can be dropped (since it is also multiplied by $\dt^2$).  If appropriate terms are dropped so that the stencil width of the ME scheme is $2p+1$, then the time-step restriction for two-dimensional problems is
\ba
     c^2 \dt^2 \Big( \f{1}{\dx^2} + \f{1}{\dy^2} \Big) < 1  , \label{eq:waveEquationTimeStep}
\ea
for orders of accuracy $2p=2,4,6$, as given by Theorem~\ref{th:stabilityWaveEquation} discussed below.  We call this version the \textit{compact} ME scheme, and we conjecture that the condition in~\eqref{eq:waveEquationTimeStep} holds at any even order $2p=2,4,6,\ldots$ (with a similar result holding for three-dimensional problems).

The compact ME scheme with LCBC conditions thus has some nice properties.  It achieves high-order accuracy in space and time in a single step.  In addition, the time-step restriction does not change as the order of accuracy increases, in contrast to some other high-order accurate schemes (e.g.~explicit multi-step methods)  where the stable time-step decreases significantly as the order of accuracy increases.

\begin{theorem}[Stability of approximations for the wave equation]  \label{th:stabilityWaveEquation}
The IBVP in~\eqref{eqn:MainEqn} for the wave equation with $q=2$ and $Q=c\sp2\Delta$ discretized  to orders $2p=2,4,6$ with the compact ME time-stepping scheme and the LCBC method on a Cartesian grid with Dirichlet or Neumann boundary conditions is stable under the time-step restriction given in~\eqref{eq:waveEquationTimeStep}.
\end{theorem}
\begin{proof}
  Let the domain be $\Omega = [0,L_x]\times[0,L_y]$, i.e.~a physical domain with lengths $L_x$ and $L_y$.
  We consider the case of Dirichlet boundary conditions on the left and
  right faces and Neumann boundary conditions on the top and bottom. The proof for other combinations
  of boundary conditions follow in a similar way.
  Let us look for normal mode solutions of the form
  \ba
      W_\iv^n = A^n \, \kappa_x^i \, \kappa_y^j ,  \label{eq:normalModeWave}
  \ea
  where $A$ is an amplification factor, $(\kappa_x,\kappa_y)$ are constants and $\iv=(i,j)$.  Since the LCBC approach leads to discrete boundary conditions that enforce even and odd symmetry, we can
  look for normal-mode solutions in space that satisfy these symmetry conditions.
  In this case we find that the normal modes are
  \ba
     W_\iv^n = A_{\pm}^n \, \sin\Big( \f{\pi k_x}{L_x} x_i\Big) \, 
               \cos\Big( \f{\pi k_y}{L_y} y_j\Big), \qquad k_x=1,2,\ldots,N_x-1, \quad k_y=0,1,2,\ldots,N_y,
  \ea
  where $A_{\pm}$ are two possible values for the amplification factor (see below).
  Any grid function $V_\iv$ satisfying the boundary conditions can then be represented as a sum of normal modes,
  \ba
     V_\iv^n = \sum_{k_x=1}^{N_x-1} \sum_{k_y=0}^{N_y}  \Vhat_{\kv}\, \sin\Big( \f{\pi k_x}{L_x} x_i\Big) \, 
               \cos\Big( \f{\pi k_y}{L_y} y_j\Big), 
  \ea
  for some coefficients $\Vhat_{\kv}$, where $\kv=(k_x,k_y)$.
  The general solution to the IBVP takes the form,
  \ba
     U_\iv^n = \sum_{k_x=1}^{N_x-1} \sum_{k_y=0}^{N_y}  
        \Big( A_{+,\kv}^n \Vhat_{+,\kv}\, + A_{-,\kv}^n \Vhat_{-,\kv} \Big) \, \sin\Big( \f{\pi k_x}{L_x} x_i\Big) \, 
               \cos\Big( \f{\pi k_y}{L_y} y_j\Big), 
  \ea
  where the coefficients $\Vhat_{\pm,\kv}$ are determined from the two initial conditions.

  For stability we choose $\dt$ so that $|A_{\pm,\kv}|\le 1$ for all valid $k_x$
  and $k_y$. 
  It is straightforward to find the symbols of $\Dpx\Dmx$ and $\Dpy\Dmy$,
  \begin{subequations}
  \ba
    &  \Dpx\Dmx \sin\Big( \f{\pi k_x}{L_x} x_i\Big)  = -\kHat_x ^2 \,\sin\Big( \f{\pi k_x}{L_x} x_i\Big), \\
    &  \Dpy\Dmy \cos\Big( \f{\pi k_y}{L_y} y_j\Big)  = -\kHat_x ^2 \,\cos\Big( \f{\pi k_y}{L_y} y_j\Big),
  \ea
  \end{subequations}
where
\begin{equation}\label{eq:kHatxyDef}
    \kHat_x \eqdef  \f{\sin(\xi_x/2)}{\dx/2},  \qquad \kHat_y \eqdef  \f{\sin(\xi_y/2)}{\dy/2}, \qquad \xi_x \eqdef \f{\pi k_x}{L_x} \dx , \qquad \xi_y \eqdef \f{\pi k_y}{L_y} \dy .
\end{equation}
  Substituting~\eqref{eq:normalModeWave} into the ME time-stepping schemes for the different orders of accuracy, determined by $p$, leads to a quadratic equation
  for $A$,
  \ba
       A^2 - 2 b_p A + 1 = 0 , \qquad p=1,2,3,\label{eq:stabQuad}
  \ea
  where $b$ depends on the various parameters of the discretization. Stability requires $b_p\in\Real$ and $|b_p|<1$. 
  Note that when $b_p=\pm 1$ there is a double root for $A$ which leads to algebraic growth which we exclude.

  \mni
  For $p=1$, 
  \ba
    b_1 =  1 - 2 \bigl( \lambdaHat_x^2 + \lambdaHat_y^2\bigr)
  \ea
  where
\ba
  \lambdaHat_x \eqdef c\dt  \, \f{\kHat_x}{2}  ,\qquad \lambdaHat_y \eqdef c \dt \, \f{\kHat_y}{2},
\ea 
with 
\ba
  | \lambdaHat_x |  \le \f{c\dt}{\dx} , \qquad |\lambdaHat_y| \le \f{c\dt}{\dy} .
\ea
Note that $b_1<1$ is clearly satisfied, while the condition $b_1>-1$ implies
\ba
   \max_{\{k_x,k_y\}}\,  \bigl(\lambdaHat_x^2 + \lambdaHat_y^2\bigr) < 1  ,
\ea
and this implies the time-step restriction in~\eqref{eq:waveEquationTimeStep}.

\mni
For $p=2$, 
\ba 
  b_2 =   1- 2\Big( \lambdaHat_x^2 + \lambdaHat_y^2 
+\f{\dx^2}{12}\lambdaHat_x^2\,\kHat_x^2 +\f{\dy^2}{12}\lambdaHat_y^2\,\kHat_y^2 \Big) 
+ \f{2}{3} \Big( \lambdaHat_x^2 + \lambdaHat_y^2 \Big)^2   . \label{eq:b2Def}
\ea
From \eqref{eq:stabQuad}, we find $A_{\pm} = b_2 \pm \sqrt{b_2^2 - 1}$ which we express in terms of the four variables \linebreak $(c\dt/\dx,c\dt/\dy,\dx\kHat_x,\dy\kHat_y)$. For each $(c\dt/\dx,c\dt/\dy)$, we define
\begin{equation}\label{eq:AmaxDef}
A_{\text{max}} = \max\left\{\max_{\overset{\scriptstyle-2\leq \dx\kHat_x\leq2}{-2\leq \dy\kHat_y\leq2}} \abs{A_{+}},\max_{\overset{\scriptstyle-2\leq \dx\kHat_x\leq2}{-2\leq \dy\kHat_y\leq2}} \abs{A_{-}}\right\},
\end{equation}
and find the region in the $\left(c\dt/\dx,c\dt/\dy\right)$ plane such that $A_{\text{max}}\leq 1$. 
We run a similar experiment for $p=3$, where $b_3$ takes the form 
\ba
\begin{split}
  b_3 = &  1- 2\left[ \lambdaHat_x^2\left(1 + \f{\dx^2}{12}\kHat_{x}^2 + \f{\dx^4}{90}\kHat_{x}^4\right) + \lambdaHat_y^2 \left(1 + \f{\dy^2}{12}\kHat_{y}^2 + \f{\dy^4}{90}\kHat_{y}^4\right)\right]\\
  & \quad+ \f{2}{3} \left[\lambdaHat_{x}^{4}\left(1+\f{\dx^2}{6}\kHat_{x}^2\right)+\lambdaHat_{y}^{4}\left(1+\f{\dy^2}{6}\kHat_{y}^2\right)+ 2\lambdaHat_{x}^2\lambdaHat_{y}^2\left(1+\f{\dx^2}{12}\kHat_{x}^2\right)\left(1 + \f{\dy^2}{12}\kHat_{y}^2\right)\right]\\
  & \qquad -\f{4}{45}\left(\lambdaHat_{x}^2 + \lambdaHat_{y}^2\right)^3. 
  \end{split}
\ea
Figure \ref{fig:stabRegion} shows that the stability region, $A_{\text{max}}\le1$, for both the fourth-order ($p=2$) and sixth-order ($p=3$) accurate time-stepping schemes.  The stability region for both schemes is found to lie within the unit circle, and thus $\dt$ satisfies the condition in~\eqref{eq:waveEquationTimeStep} when $p=2$ and~$3$.

\noindent{\Huge $\square$}
\end{proof}
{
	\newcommand{\figWidth}{9cm}
	\newcommand{\trimfig}[2]{\trimw{#1}{#2}{0.}{0.}{.0}{.0}}  
	\newcommand{\trimfigb}[2]{\trimw{#1}{#2}{.0}{.05}{.1}{.1}} 
	\begin{figure}[H]
		\begin{center}
			\begin{tikzpicture}[scale=1]
				\useasboundingbox (0,0) rectangle (7,7);  
				\draw (-2,-0.5) node[anchor=south west] {\trimfigb{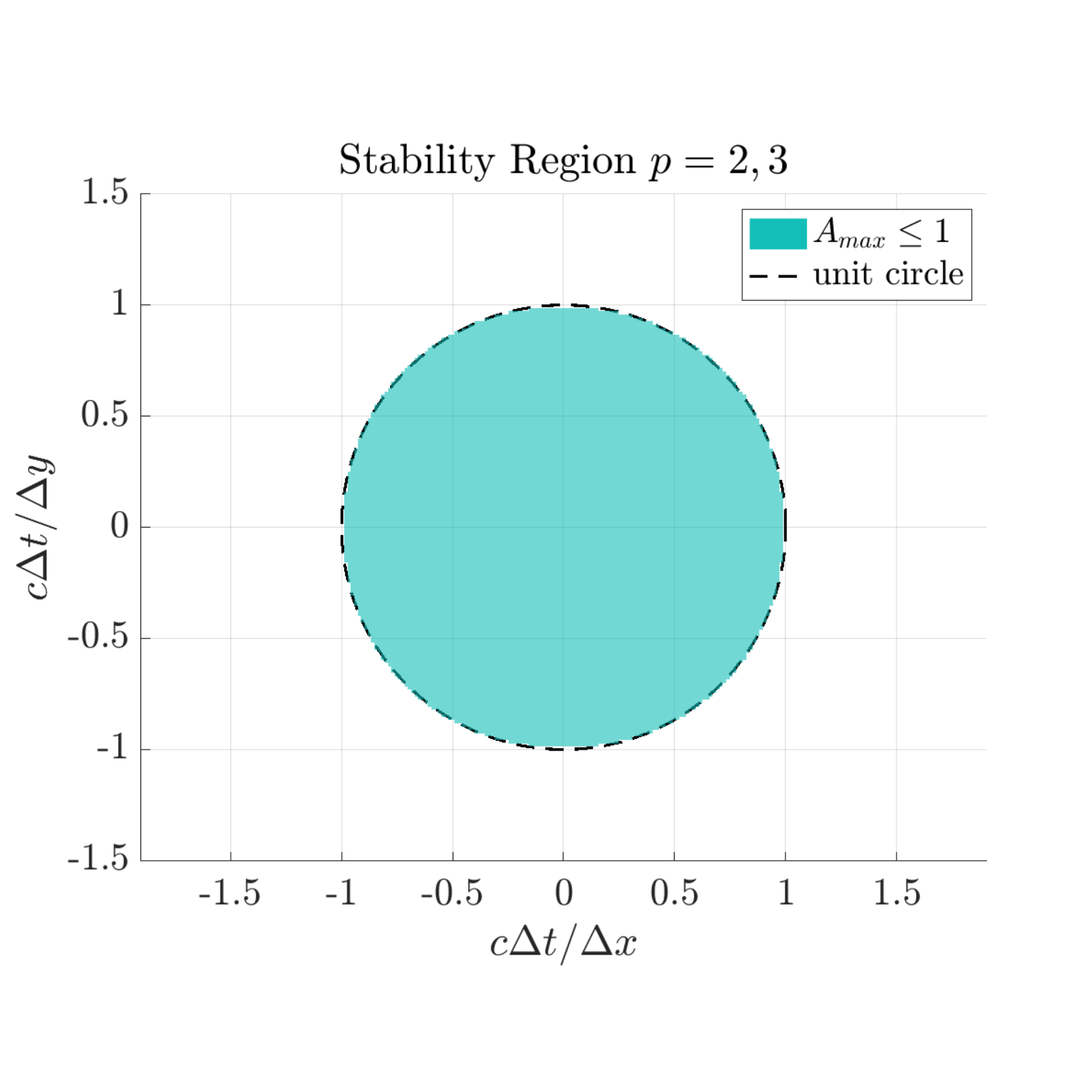}{\figWidth}};
			\end{tikzpicture}
		\end{center}
		\caption{Stability region of the fourth-order and sixth-order accurate ME time-stepping schemes for the wave equation on a Cartesian grid using the LCBC approach.}
		\label{fig:stabRegion}
	\end{figure}
}

In \ref{app:stabilityProof2}, we provide an analytical proof for the stability results observed in Figure~\ref{fig:stabRegion} when $p = 2$.

%% file: tex/results.tex
\section{Numerical results} \label{sec:NumericalResults}

The aim of this section is to demonstrate the accuracy of the LCBC method by considering several examples involving second-order and higher-order accurate discretizations of the model problem in~\eqref{eqn:MainEqn}.  A semi-discrete scheme for this problem is given in~\eqref{eqn:SemiDiscreteEqn}, and our first task is to describe suitable time-stepping schemes for the cases when $q=1$ and~$2$.  We then utilize the method of manufactured solutions to verify the accuracy of the fully discrete schemes with the LCBC method for problems with constant and variable coefficients defined on rectangular and curvilinear domains.  For the latter case, a mapping is used to bring the problem to a unit-square computational domain.  Finally, we consider physically-motivated examples involving heat flow and wave propagation as further tests of the LCBC method.

We note that when solving the wave equation, which has no inherent dissipation, some upwind dissipation is often required to retain a stable time-stepping scheme for problems with variable coefficients or for discretizations on curvilinear or overset grids.  For wave equations in second-order form an effective upwind dissipation is described in~\cite{sosup2012,mxsosup2018}.  For the results presented here, however, no upwind dissipation has been added.

\subsection{Fully discrete schemes}\label{subsec:DiscreteSchemes}

For the case of elliptic problems ($q=0$), the discretization in~\eqref{eqn:SemiDiscreteEqn} becomes
\begin{equation}
	\begin{cases}
		Q_{d,h}U_{\iv} = - f(\xv_{\iv}), & \xv_{\iv}\in\Omega_h, \\
		\Bc_{d,h}U_{\iv} = g(\xv_{\iv}), & \xv_{\iv} \in\p\Omega_h,
	\end{cases}
\end{equation}
where $Q_{d,h}$ is the $d$-th order accurate approximation of the differential operator $Q$ given in~\eqref{eqn:OperatorQ} applied for points $\xv_\iv$ on the grid $\Omega_h$.  The discrete operator $\Bc_{d,h}$ is a $d$-th order accurate approximation of the boundary operator $\Bc$ given in~\eqref{eqn:OperatorB} applied on the boundary $\p\Omega_h$, and the LCBC method is used to obtain ghost values on the extended domain.

For the time-dependent cases ($q=1,2$), we require methods of time-stepping.  For these cases, let $t^n = n\dt$, $n = 0,1, \dots, \Nt$, for a fixed time step $\dt=T/\Nt$, and let $U_{\iv}^{n}\approx u(\xv_{\iv},t^{n})$.  For the parabolic case ($q=1$), we consider an explicit forward-Euler (FE) time-stepping scheme as well as implicit Backward Differentiation Formula (BDF) time-stepping schemes.  For the explicit FE scheme, we choose a stable time step $\dt$ 
given by~\footnote{This condition arises from approximating the region of absolute stability of
the time-stepping scheme by an ellipse with semi-axes of lengths $\alpha_d$ and $\beta_d$.}
\begin{equation}
\dt = \dfrac{{\rm cfl}}{\sqrt{\left(\frac{\alpha_d}{2}\right)^2 + \beta_d^2}},
\end{equation}
where ${\rm cfl}\le1$ is a safety factor, set to ${\rm cfl}=0.9$ for all calculations, and $\alpha_d$ and $\beta_d$ are estimates for the real and imaginary parts of the time-stepping eigenvalue corresponding to $Q_{d,h}$ with frozen coefficient functions.  These values are given by
\begin{equation}
\begin{array}{l}
\displaystyle{
\alpha_d=r_d\left({c_{11,{\rm max}}\over\dx^2}+{c_{22,{\rm max}}\over\dy^2}\right)+s_d{\vert c_{12}\vert_{{\rm max}}\over\dx\dy}+\vert c_0\vert_{{\rm max}}
}\medskip\\
\displaystyle{
\beta_d=q_d\left({\vert c_{1}\vert_{{\rm max}}\over\dx}+{\vert c_{2}\vert_{{\rm max}}\over\dy}\right)
}
\end{array}
\end{equation}
where the maximum values of the coefficient functions (or their absolute values) 
are taken for $\xv\in\Omega$, and the parameter values are taken to be $(r_d,s_d,q_d)=(4,\,2,\,1)$,
$(16/3,\,9/2,\,3/2)$ and $(272/45,\,131/261,\,1199/756)$, for $d=2$, $4$ and $6$, respectively.
Details of the stability analyses leading to these conditions can be found in~\ref{sec:parabolicStepRestrictions}.
For the implicit BDF schemes, we take
\begin{equation}
\dt = \min\{\dx,\dy\},
\end{equation}
and match the temporal accuracy of the BDF scheme with the chosen spatial order of accuracy.
The fully discrete schemes for the hyperbolic case ($q=2$) are based on a modified equation (ME) approach as discussed, for example, in~\cite{max2006b} for the case of Maxwell's equations in second-order form.  
For the model problem here, the baseline second-order accurate ($d=2$) scheme is given by
\begin{equation}\label{eqn:ModifiedEqn2}
	\begin{cases}
		\dpmt U_{\iv}^n  = Q_{2,h}U_{\iv}^n + f(\xv_{\iv},t^n), & \xv_{\iv}\in \Omega_h,\\
		\mathcal{B}_{2,h}U_{\iv}^{n+1} = g(\xv_{\iv},t^{n+1}), & \xv_{\iv}\in\partial\Omega_h,\\
		U_{\iv}^0 = u_0(\xv_{\iv}),& \xv_{\iv}\in \overline{\Omega}_h, \\
		U_{\iv}^1 = u_0(\xv_{\iv}) + \dt u_1(\xv_{\iv}) + \frac{\dt^2}{2}\left[Qu_0(\xv_{\iv})+f(\xv_{\iv},0)\right], & \xv_{\iv}\in \overline{\Omega}_h,
	\end{cases}
\end{equation}
while the fourth and sixth-order accurate schemes ($d=4$ and $6$) are
\begin{equation}\label{eqn:ModifiedEqn4}
	\begin{cases}
		\dpmt U_{\iv}^n  = Q_{4,h}U_{\iv}^n + f(\xv_{\iv},t^n) + \frac{\dt^2}{12}\left[Q^2_{2,h}U_{\iv}^n + \Psi_2 f(\xv_{\iv},t^n)\right], & \xv_{\iv}\in \Omega_h,\\
		\mathcal{B}_{4,h}U_{\iv}^{n+1} = g(\xv_{\iv},t^{n+1}), & \xv_{\iv}\in\partial\Omega_h,\\
		U_{\iv}^0 = u_0(\xv_{\iv}),& \xv_{\iv}\in \overline{\Omega}_h, \\
		U_{\iv}^1 = u_0(\xv_{\iv}) + \dt u_1(\xv_{\iv})  \\[2pt]
\hspace{1.2cm}+ \frac{\dt^2}{2}\left[Qu_0(\xv_{\iv})+f(\xv_{\iv},0)\right] + \frac{\dt^3}{6}\left[Qu_1(\xv_{\iv})+\p_tf(\xv_{\iv},0)\right]  \\[2pt]
		\hspace{1.7cm} +\frac{\dt^4}{24}\left[Q^2u_0(\xv_{\iv}) + \Psi_2 f(\xv_{\iv},0)\right], & \xv_{\iv}\in \overline{\Omega}_h,
	\end{cases}
\end{equation}
and
\begin{equation}\label{eqn:ModifiedEqn6}
	\begin{cases}
		\dpmt U_{\iv}^n  = Q_{6,h}U_{\iv}^n + f(\xv_{\iv},t^n)\\[2pt]
\hspace{2.6cm}+ \frac{\dt^2}{12}\left[{\blue \Qtilde^2_{4,h}}U_{\iv}^n + \Psi_2 f(\xv_{\iv},t^n)\right] + \frac{\dt^4}{360}\left[Q_{2,h}^3U_{\iv}^n + \Psi_3 f(\xv_{\iv},t^n)\right], & \xv_{\iv}\in \Omega_h,\\
		\mathcal{B}_{6,h}U_{\iv}^{n+1} = g(\xv_{\iv},t^{n+1}), & \xv_{\iv}\in\partial\Omega_h,\\
		U_{\iv}^0 = u_0(\xv_{\iv}),& \xv_{\iv}\in \overline{\Omega}_h, \\
		U_{\iv}^1 = u_0(\xv_{\iv}) + \dt u_1(\xv_{\iv}) \\[2pt]
\hspace{1.2cm}+ \frac{\dt^2}{2}\left[Qu_0(\xv_{\iv})+f(\xv_{\iv},0)\right] + \frac{\dt^3}{6}\left[Qu_1(\xv_{\iv})+\p_tf(\xv_{\iv},0)\right]\\[2pt]
\hspace{1.7cm}+\frac{\dt^4}{24}\left[Q^2u_0(\xv_{\iv}) + \Psi_2 f(\xv_{\iv},0)\right] + \frac{\dt^5}{120}\left[Q^2u_1(\xv_{\iv})+\pt\Psi_2 f(\xv_{\iv},0)\right]\\[2pt]
		\hspace{2.2cm} + \frac{\dt^6}{720}\left[Q^3u_0(\xv_{\iv}) + \Psi_3 f(\xv_{\iv},0)\right], & \xv_{\iv}\in \overline{\Omega}_h,
	\end{cases}
\end{equation}
respectively, where ${\blue \Qtilde_{4,h}^2}$ directly approximates the continuous operator $Q^2$ to fourth-order accuracy rather than forming the square of the discrete operator $Q_{4,h}$ which would have a stencil greater than $2p+1$. This approach leads to the compact ME scheme noted previously in Section~\ref{sec:stabilityWaveEquation}.

For all three schemes, an estimate for a stable time step is given by
\begin{equation}
	\dt = \dfrac{{\rm cfl}}{\sqrt{\frac{c_{11,{\rm max}}}{\dx^2} + \frac{c_{22,{\rm max}}}{\dy^2}}},
	   \label{eq:dtWave}
\end{equation}
where ${\rm cfl} < 1$.  We use ${\rm cfl} = 0.9$ for all numerical examples.   Details of the stability analysis leading to~\eqref{eq:dtWave} can 
be found in~\ref{sec:hyperbolicStepRestrictions}.

\subsection{Convergence tests using the method of manufactured solutions}\label{subsec:MSmethod}

We first verify the accuracy of the fully discrete schemes with the LCBC method for the problem in~\eqref{eqn:MainEqn} using manufactured solutions.  In this approach, we set the forcing functions of the PDE and boundary conditions, and the initial conditions, to be
\begin{equation}
	\begin{cases}
		f(\xv,t) = \Lc_q u_e(\xv,t) - Qu_e(\xv,t), & \xv \in \Omega, \quad t\in(0,T], \quad q = 0,1,2,\\
		g(\xv,t) = \Bc u_e(\xv,t), & \xv \in\p\Omega, \quad t\in[0,T], \\
		u_{\alpha - 1}(\xv) = \pt^{\alpha - 1}u_e(\xv,0), & \xv\in \overline{\Omega}, \quad \alpha = 1, \dots, q, \quad q  = 1,2,
	\end{cases}
\end{equation}
so that $u=u_e(\xv,t)$ is an exact solution of the problem in~\eqref{eqn:MainEqn}.  For all tests, we select an exact solution of the form 
\begin{equation}\label{eqn:ManufacturedSoln}
	u_e(\xv,t) = \cos(k_x x)\cos(k_y y)\phi(t), \qquad \xv\in \overline{\Omega}, \quad t\in[0,T],
\end{equation}
where $(k_x,k_y)=(2\pi,\sqrt2\pi)$ are wave numbers, and $\phi(t)$ is taken to be
\begin{equation}
\phi(t) = \begin{cases}
1, & \hbox{if $q = 0$}, \\
t + 1, & \hbox{if $q = 1$ and explicit FE time-stepping}, \\
\cos(\pi t), & \hbox{otherwise}.
\end{cases}
\end{equation}
When $q\neq 0$, we choose the time-stepping schemes and $\phi(t)$ such that the temporal accuracy is at least equal to the spatial accuracy of the discretized problem. For example, the FE scheme, used for parabolic problems, is first-order accurate in time.  For this case, we choose $\phi(t)$ to be a linear polynomial since the FE time-stepping scheme is exact for polynomials of degree one or less and thus the convergence rates of the fully discrete scheme reflect the spatial accuracy. On the other hand, the implicit BDF methods for problems with $q=1$ and the explicit schemes for problems with $q = 2$ have temporal orders of accuracy balanced with the spatial accuracy.  For such high-order accurate schemes, $\phi(t)$ may be any smooth function.  Note that we use exact values taken from the manufactured solution to start the multi-step BDF schemes.

We now consider a sequence of three tests of the fully discrete schemes.

\subsubsection{Test 1. Square domain and constant-coefficient functions}\label{NumRes:example1}

The first test considers the set-up discussed previously in Section~\ref{sec:wave2d}.  The domain for the problem is the unit square, $\Omega=(0,1)^2$, with $Q=\Delta$, the usual Laplacian operator.  The previous discussion focused on the wave equation, $q=2$, but here we also consider the cases when $q=0$ and~$1$.  The boundary conditions on the left ($x=0$) and right ($x=1$) are taken to be Dirichlet conditions, while Neumann conditions are assumed on the bottom ($y=0$) and top ($y=1$), as given in~\eqref{eq:wave2d}.  Numerical solutions are computed on uniform grids with $\dx=\dy=h=h_0/2^{j}$, $h_0=1/10$, $j=0,1,2,3$, and the maximum error between the numerical solution and the exact (manufactured) solution is determined for each value of $h$, and at a time $T=1$ for the time-dependent cases, $q=1,2$.  

\input tex/MSexample1

Figure~\ref{fig:MSexample1} shows maximum errors versus grid spacing for the elliptic case ($q=0$), the parabolic case ($q=1$) for both FE and BDF time-stepping, and the hyperbolic case ($q=2$).  Calculations are performed using schemes with order of accuracy $d=2$, $4$ and~$6$.  In each case, the computed errors can be compared with reference lines indicating the design order of accuracy.  We observe that the computed errors show the expected convergence rates.

\subsubsection{Test 2. Square domain and variable-coefficient functions}\label{NumRes:example2}

The second test also uses a unit-square domain and boundary conditions as in the previous test, but now the coefficient functions in the operator $Q$ defined in~\eqref{eqn:OperatorQ} are taken to be
\begin{alignat*}{5}
	&c_{11}(\xv) = \dfrac{e^{x}}{2 + x + y},&\qquad& c_{12}(\xv) = xye^{-x-y-1},& \qquad& c_{22}(\xv) = \dfrac{e^{x+y}}{2 + y}, \\
	&c_1(\xv) =1+xy, &\qquad& c_2(\xv) = 1 + \frac{x}{1+y},&\qquad & c_0(\xv) = 2+x-y. 
\end{alignat*}
Note that for this choice, we have
\begin{equation*}
\max_{\xv\in\overline{\Omega}}\dfrac{\abs{c_{12}(\xv)}}{\sqrt{c_{11}(\xv)c_{22}(\xv)}} \approx 0.3772,
\end{equation*} 
so that the cross term is significant, but small enough so that the corner conditions are sufficiently well conditioned.  Numerical solutions are computed on uniform grids with $\dx=\dy=h=h_0/2^{j}$, $h_0=1/10$, $j=0,1,2,3$, as before, and maximum errors are determined at $T=1$ for the time-dependent cases.

\input tex/MSexample2

Figure~\ref{fig:MSexample2} shows maximum errors versus grid spacing for the elliptic case ($q=0$), the parabolic case ($q=1$) for both FE and BDF time-stepping, and the hyperbolic case ($q=2$).  Comparing the computed errors with the reference lines, we again observe the expected rate of convergence of the fully discrete schemes for $d=2$, $4$ and~$6$.

\subsubsection{Test 3. Curvilinear domain and constant-coefficient functions}\label{NumRes:example3}

As a final test using the method of manufactured solutions, we consider a problem defined on an annular domain, $\Omega^{P}=\{\xv=(\rho\cos\theta,\rho\sin\theta)\,\vert\,\rho_{1}<\rho<\rho_{2},\,0<\theta<\theta_2\}$, where $(\rho,\theta)$ are polar coordinates.  The parameters $(\rho_1,\rho_2)$ define the inner and outer radii of the problem domain and $\theta_2$ defines its angular extent.  In the annular domain, the problem is taken to be
\begin{equation}\label{eq:MSexampleP}
	\begin{cases}
		\Lc_qu = \Delta u + f(\xv,t),  & \xv\in \Omega^{P}, \quad q = 0,1,2, \\[5pt]
		u(\xv ,t) = g_{\ell}(\xv,t),  			& \xv\in\p\Omega^{P}_{\ell}, \\ 
		u(\xv ,t) = g_{b}(\xv,t), 				& \xv\in\p\Omega^{P}_{b}, \\
		\p_n u(\xv,t) = g_{r}(\xv,t), 		&   \xv\in\p\Omega^{P}_{r}, \\
		\p_n u(\xv,t) = g_{\tau}(\xv,t), 	&   \xv\in\p\Omega^{P}_{\tau}, \\[5pt]
		\pt^{\alpha-1}u(\xv,0) = u_{\alpha-1}(\xv), & \xv \in\overline{\Omega}^{P}, \quad \alpha = 1,\dots, q, \quad q = 1,2,
	\end{cases}
\end{equation}
where $\p\Omega_{{\rm s}}^{P}$ for ${\rm s} = \ell,r,b,\tau$ correspond to the boundaries at $\rho= \rho_1$, $\rho = \rho_2$, $\theta = 0$ and $\theta = \theta_2$ respectively.  The problem is transformed from the ``physical'' domain $\xv\in\Omega^{P}$ to the computational domain $\rv=(r,s)\in\Omega=(0,1)^2$ using the mapping
\begin{equation} \label{eq:mapping}
\xv = \Gv(\rv)= \begin{bmatrix}
\bigl(\rho_1+r(\rho_2-\rho_1)\bigr)\cos(s\theta_2)\\[5pt] \bigl(\rho_1+r(\rho_2-\rho_1)\bigr)\sin(s\theta_2)
\end{bmatrix}, \qquad \rv = (r,s)\in \Omega.
\end{equation}
In terms of the mapped coordinates, the spatial operator of the PDE becomes
\begin{equation}
Qu=\p_r^2u+{4\over\pi^2(1+r)^2}\p_s^2u+{1\over(1+r)}\p_ru,
\end{equation}
for the choice $\rho_1=1$, $\rho_2=2$ and $\theta_2=\pi/2$.  The mapped problem is discretized in the computational domain using $\dr = \dr_0/2^j$, $\dr_0=1/10$, $j=0,1,2,3$, $\ds = \dr/2$, and $h = \min(\dr,\ds)=\ds$.\linebreak  Note that for each grid resolution~$j$ the corresponding grid spacings in physical space are approximately equal.

\input tex/MSexample3

Figure~\ref{fig:MSexample3} shows maximum errors at $T=1$ versus grid spacing.  The calculations are performed using various schemes for $q=0,1,2$ with order of accuracy $d=2$, $4$ and~$6$.  In each case, we observe that the computed errors agree with the expected convergence rates.

\subsection{Scattering of a plane wave from a cylinder}

Moving on from tests involving the method of manufactured solutions, we now consider three physically-motivated problems.  The first of these problems involves the scattering of an incident plane wave from a cylinder of radius equal to one.  The incident wave, traveling from left to right, is given by
\begin{equation}
u_{\text{inc}}(\xv,t) =\cos[k(x - ct)],
\end{equation}
where $k$ is the wave number and $c$ is the wave velocity.  An exact solution of the wave equation 
for the scattered field $u(\xv,t)$ from a hard cylinder with $u=-u_{\text{inc}}(\xv,t)$ for $\vert\xv\vert=1$ is given by
\begin{equation}\label{eqn:ScatteringExact}
	u_e(\xv,t) = -\Re \left[e^{-ikct}\left(\dfrac{J_0(k)}{H_{0}^{(1)}(k)} H_{0}^{(1)}(k\rho) + 2\sum_{n = 1}^{\infty}i^{n}\dfrac{J_n(k)}{H_n^{(1)}(k)}H_{n}^{(1)}(k\rho)\cos(n\theta)\right)\right], \quad \vert\xv\vert>1,
\end{equation}
where $J_n$ and $H_n^{(1)}$ are Bessel and Hankel functions of the first kind, respectively, and $(\rho,\theta)$ are polar coordinates (see~\cite{martin_2006}).  Our aim is to compute the scattered field for this problem numerically by solving a corresponding half-plane problem with a symmetry condition applied along the mid-line, $\theta=0$ and $\pi$ for $\rho>1$.  For this half-plane problem, we consider a finite domain given by
\begin{equation*}
    \Omega^{P}=\{\xv=(\rho\cos\theta,\rho\sin\theta)\,\vert\,1<\rho<\rho_{\infty},\,0<\theta<\pi\},
\end{equation*}
and let $u(\xv,t)$ satisfy
\begin{equation}\label{eqn:ScatteringMain}
	\begin{cases}
		\pt^2u = c^2\Delta u, & \xv\in \Omega^P,\quad t\in(0,T],\\[5pt]
		u(\xv,t) = - u_{\text{inc}}(\xv,t), & \rho = 1,\quad \theta \in[0,\pi], \\
		u(\xv,t) = u_{e}(\xv,t), & \rho = \rho_\infty,\quad \theta \in[0,\pi], \\
		\p_n u(\xv,t) = 0, & \theta = 0, \; \theta = \pi,\quad \rho>1, \\[5pt]
		u(\xv,0) = u_e(\xv,0), & \xv\in\overline{\Omega}^P,\\
		\pt u(\xv,0) = \pt u_e(\xv,0), & \xv \in \overline{\Omega}^P.
	\end{cases}
\end{equation}
The exact solution of the problem in~\eqref{eqn:ScatteringMain} is given by~\eqref{eqn:ScatteringExact} for any choice of the outer radius $\rho_\infty>1$ of the domain~$\Omega^P$.

We solve \eqref{eqn:ScatteringMain} numerically to a final time $T = 1$ for the case $k=30$, $c=1$ and $\rho_\infty=2$.  This is done by first mapping the problem to a unit-square computational domain using~\eqref{eq:mapping} with $\rho_1=1$, $\rho_2=\rho_\infty$ and $\theta_2=\pi$.  On the computational domain $\rv=(r,s)\in\Omega$, we use the explicit time-stepping schemes in~\eqref{eqn:ModifiedEqn2}, \eqref{eqn:ModifiedEqn4} and~\eqref{eqn:ModifiedEqn6} corresponding to orders of accuracy equal to $d=2$, $4$ and $6$, respectively, in both time and space.  The LCBC method is used to handle the boundary and corner conditions for each scheme.  Numerical solutions are computed using grids with $\dr=\dr_0/2^j$, $\dr_0=1/20$, $j=0,1,2,3$, and with $\ds=\dr/5$ so that the grid spacings in physical space are approximately equal.  Figure~\ref{fig:example4} illustrates the results of the calculations.  Maximum errors are computed for each time-stepping scheme and for each grid resolution, and the results are shown in the upper-left plot in the figure (using $h=\min\{\dr,\ds\}=\ds$).  Here, we observe that the errors verify the expected convergence rates of the three schemes.  The coarsest grid used for the calculations is shown in the lower-left plot for reference.  The plots in the right column of the figure show the scattered field (top), the error in the scattered field (middle), and the total field (bottom), all at $T=1$ computed using the sixth-order ($d=6$) accurate scheme on the finest grid.  We note that the error in the computed scattered field is smooth on the grid as expected based on the convergence behavior of the solutions.

\input tex/Example4.tex

\subsection{Heat flow in a wavy channel}\label{NumRes:HeatFlowWavyChannel}

For the next case, we consider a heat flow problem in channel domain denoted by $\xv\in\Omega^P$.  Assuming the temperature $u(\xv,t)$ is specified on the boundary of the domain, the problem is given by
\begin{equation}\label{eqn:2Dheat}
\begin{cases}
u_{t} = \Dc \Delta u - \vv\cdot\grad{u} + \gamma u, & \xv\in\Omega^P, \quad t\in(0,T], \\
u(\xv,t) = g(\xv,t), & \xv\in\p\Omega^P, \\
u(\xv,0) = u_{0}(\xv), & \xv\in{\bar\Omega}^P,
\end{cases}
\end{equation}
where $\Dc$ is a diffusivity, $\vv$ is a convection velocity and $\gamma$ is a reaction rate, all taken to be constants, and $u_{0}(\xv)$ is the initial temperature.  An exact solution for this problem can be constructed by first considering the free-space solution of the PDE given by
\begin{equation}\label{eqn:2DheatGsol}
	u_e(\xv,t) = \dfrac{e^{\gamma t}}{4\pi\Dc t}\int\!\!\!\int_{\mathbb{R}^2} u_{0}(\xv^\prime)\exp\left(\frac{-\vert\xv-\xv^\prime-\vv t\vert^2}{4\Dc t}\right)\,d\xv^\prime.
\end{equation}
Assuming a Gaussian initial condition of the form
\begin{equation}\label{eqn:2DheatGIC}
	u_0(\xv) = \exp(-\sigma\vert\xv\vert^2), \qquad\sigma=\hbox{constant}>0,
\end{equation}
the solution in~\eqref{eqn:2DheatGsol} reduces to
\begin{equation}\label{eqn:2DheatSol}
u_e(\xv,t) = \dfrac{e^{\gamma t}}{1+ 4\sigma\Dc t}\exp\left(\frac{-\sigma\vert\xv-\vv t\vert^2}{1 + 4\sigma\Dc t}\right).
\end{equation}
We now set the Dirichlet boundary forcing $g(\xv,t)$ in~\eqref{eqn:2Dheat} to equal the free-space solution so that~\eqref{eqn:2DheatSol} becomes the exact solution of the heat flow problem for any channel domain $\Omega^P$.

We let $\xv\in\Omega^P$ be a wavy channel domain given by
\begin{equation}\label{eqn:flagregion}
\left.
\begin{array}{l}
x = G_{1}(r,s) = (x_2 - x_1)r + x_1 \smallskip\\
y = G_{2}(r,s) = (y_2 - y_1)s + y_1 + h(r)
\end{array}
\right\}\qquad (r,s)\in\Omega=[0,1]^2,
\end{equation}
where $(x_1,x_2)=(y_1,y_2)=(-1,1)$.  The shape of the bottom and top walls of the wavy channel is specified by
\begin{equation*}
	h(r) = A_c\sin(2\pi r), 
\end{equation*}
where $A_c$ is an amplitude.  The heat flow problem in~\eqref{eqn:2Dheat} is mapped to the unit-square $\Omega$, and numerical solutions are computed using the implicit BDF schemes with $d=2,4,6$ for the parameter values $\Dc = 0.2$, $\vv=(0.5,0.3)$, $\gamma=1$, $\sigma=6$ and $A_c=0.1$.  Note that the mapping in~\eqref{eqn:flagregion} is not orthogonal unlike the annular mappings considered in previous problems.

\input tex/Example5.tex

Figure~\ref{fig:example5} shows the results of the calculations at the final time $T=0.5$.  The upper-left plot shows the maximum error in the solutions with $d=2,4,6$ for grids with $h=\dr=\ds=h_0/2^j$, $h_0=1/40$, $j=0,1,2,3$.  These errors verify the expected convergence rates.  The lower-left plot illustrates the grid used for the calculations, here for $h=1/40$, and the plots in the right column show the sixth-order accurate solution and its error at $T=0.5$, both for the finest grid with $h=1/320$.  In particular, we observe that the error is smooth throughout the domain, including the boundaries and corners.

\newcommand{\pow}{z}
\subsection{Pulse propagation in curved channel}

As a final test, we consider a wave propagation problem for $u(\xv,t)$ in curved channel with closed ends.  The domain of the channel $\Omega^{P}$ is taken to be the quarter annulus described previously by the mapping in~\eqref{eq:mapping} with $\rho_1=1$, $\rho_2=2$ and $\theta_2=\pi/2$.  We assume that $u=\pt u=0$ initially, and that a pulse is generated in the channel by setting $u(\xv,t)=g_{{\rm p}}(\xv,t)$ on the vertical boundary $x=0$, $y\in[1,2]$ denoted by $\xv\in\p\Omega^P_{{\rm pulse}}$.  The form of this boundary forcing is taken to be
\begin{equation}
g_{{\rm p}}(\xv,t)=e^{-\beta(t - t_0)^2}\sin(\omega t)\left[\frac{1}{2}\cos\left(2\pi\left(y-\frac{3}{2}\right)\right)-\frac{1}{2}\right]^{\pow+1},
\end{equation}
where $\beta$, $t_0$ and $\omega$ are constants, and $\pow=d/2$ for a calculation using a scheme with order of accuracy equal to~$d$.  Assuming zero Neumann conditions on the remaining boundaries, the problem to solve is
\begin{equation}\label{eqn:AnnularPulseIBVP}
	\begin{cases}
		\pt^2u = c^2\Delta u, & \xv \in \Omega^{P}, \quad t\in(0,T],\\
		u(\xv,t) = g_{{\rm p}}(\xv,t), & \xv \in \p\Omega^P_{{\rm pulse}},\\
		\p_nu(\xv,t) = 0, & \xv \in \p\Omega^{P}\backslash\p\Omega^P_{{\rm pulse}},\\[2pt]
		u(\xv,0) = \pt u(\xv,0) = 0, & \xv \in \overline{\Omega}^{P}.
	\end{cases}
\end{equation}

\input tex/Example6.tex

The pulse propagation problem in~\eqref{eqn:AnnularPulseIBVP} is solved numerically for $c=1$, $\beta = 50$, $t_0=1$ and $\omega = 6\pi$.  This is done by mapping the problem to the unit-square $\Omega$, and then solutions are computed using the second, fourth and sixth-order accurate time-stepping schemes given by~\eqref{eqn:ModifiedEqn2}, \eqref{eqn:ModifiedEqn4} and~\eqref{eqn:ModifiedEqn6}, respectively.  The grid spacings in the unit-square coordinates $(r,s)$ are taken to be $\dr = \dr_0/2^j$, $\dr_0=1/40$, $j=0,1,2,3$, $\ds = \dr/2$, and $h = \min(\dr,\ds)= \ds$.  Figure~\ref{fig:AnnularPulseEvolution} shows the evolution of the pulse for the times $t=1,2,\ldots,6$ as determined by the sixth-order accurate scheme on the finest grid.  Here we observe the forward propagation of the pulse at an early time ($t=1$) and the subsequent interaction with the curved walls of the channel at later times ($t=2,3$).  By $t=4$ the pulse has reached the lower flat wall of the channel, and then reflects backwards into the channel at the final times ($t=5,6$).

\input tex/Example6b.tex

Figure~\ref{fig:AnnularPulse} shows the behavior of the computed error in the solutions at $t=3$ and~$6$.  Since an exact solution for the problem is not available, 
we compute the error approximately using a Richardson extrapolation approach~\cite{mog2006,BanksAslamRider2008}.
 Let $U_\iv^{(j)}$ denote the solution computed on a grid with mesh spacing $h_j$ at a time $t^{n}$.  The error is computed at a fixed time $T=t^{n}$, and so the dependence on $t^n$ is suppressed for notational convenience.  Assuming convergent schemes, we have
\[
U_\iv^{(j)}\approx u_e\bigl(\xv_\iv^{(j)}\bigr)+C_\iv^{(j)}h_j^{\sigma},
\]
where $u_e(\xv)$ is the exact solution (at $t^n=T$), $\sigma$ is the rate of convergence (for $h\rightarrow0$), and $C_\iv^{(j)}$ is a grid function (whose values are independent of $h_j$ and $\sigma$).  Subtracting solutions at grid resolutions~$j$ and~$j+1$, assuming $h_{j+1}=\half h_j$, gives
\[
U_\iv^{(j)}-R_jU_\iv^{(j+1)}\approx C_\iv^{(j)}(1-2^{-\sigma})h_j^{\sigma},
\]
where $R_j$ is an operator that restricts a grid function at resolution~$j+1$ to one at~$j$.  Taking the maximum norm of both sides gives
\[
\Delta U_j\eqdef\Vert U_\iv^{(j)}-R_jU_\iv^{(j+1)}\Vert_\infty\approx C(1-2^{-\sigma})h_j^{\sigma},
\]
where $C$ is a positive constant.  A similar expression involving solutions at grid resolutions~$j+1$ and~$j+2$ can be used to give
\begin{equation}
\sigma\approx{\log_2\left({\Delta U_j\over \Delta U_{j+1}}\right)},\qquad C\approx{\Delta U_j\over(1-2^{-\sigma})h_j^{\sigma}}.
\label{eq:estimates}
\end{equation}
The maximum error in the numerical solution at grid resolution $j$ is then estimated as
\begin{equation}
\Vert U_\iv^{(j)} - u_e\bigl(\xv_\iv^{(j)}\bigr)\Vert_\infty\approx Ch_j^{\sigma},
\label{eq:errorestimates}
\end{equation}
where the convergence rate $\sigma$ and constant $C$ are given in~\eqref{eq:estimates}.  The convergence plots at $t^n=3$ and $t^n=6$ in the left column of Figure~\ref{fig:AnnularPulse} are given by~\eqref{eq:errorestimates} using values for $\sigma$ and $C$ obtained from solutions at grid resolutions $j=1,2,3$.  The two plots show convergence rates that agree with the expected order of accuracy of the explicit time-stepping schemes for $d=2$, $4$ and~$6$.  The plots in the right column show representative behaviors of the computed error.

%% file: tex/MSexample1.tex

{
	\newcommand{\figWidth}{5.6cm}
	\newcommand{\trimfig}[2]{\trimw{#1}{#2}{0.}{0.}{.0}{.0}}  
	\newcommand{\trimfigb}[2]{\trimwb{#1}{#2}{.2}{0.1}{.2}{.15}} 
	\begin{figure}
		\begin{center}
			\begin{tikzpicture}[scale=1]
				\useasboundingbox (0,0) rectangle (11,9.5);  
				\draw ( -0.8,4.40) node[anchor=south west] {\trimfig{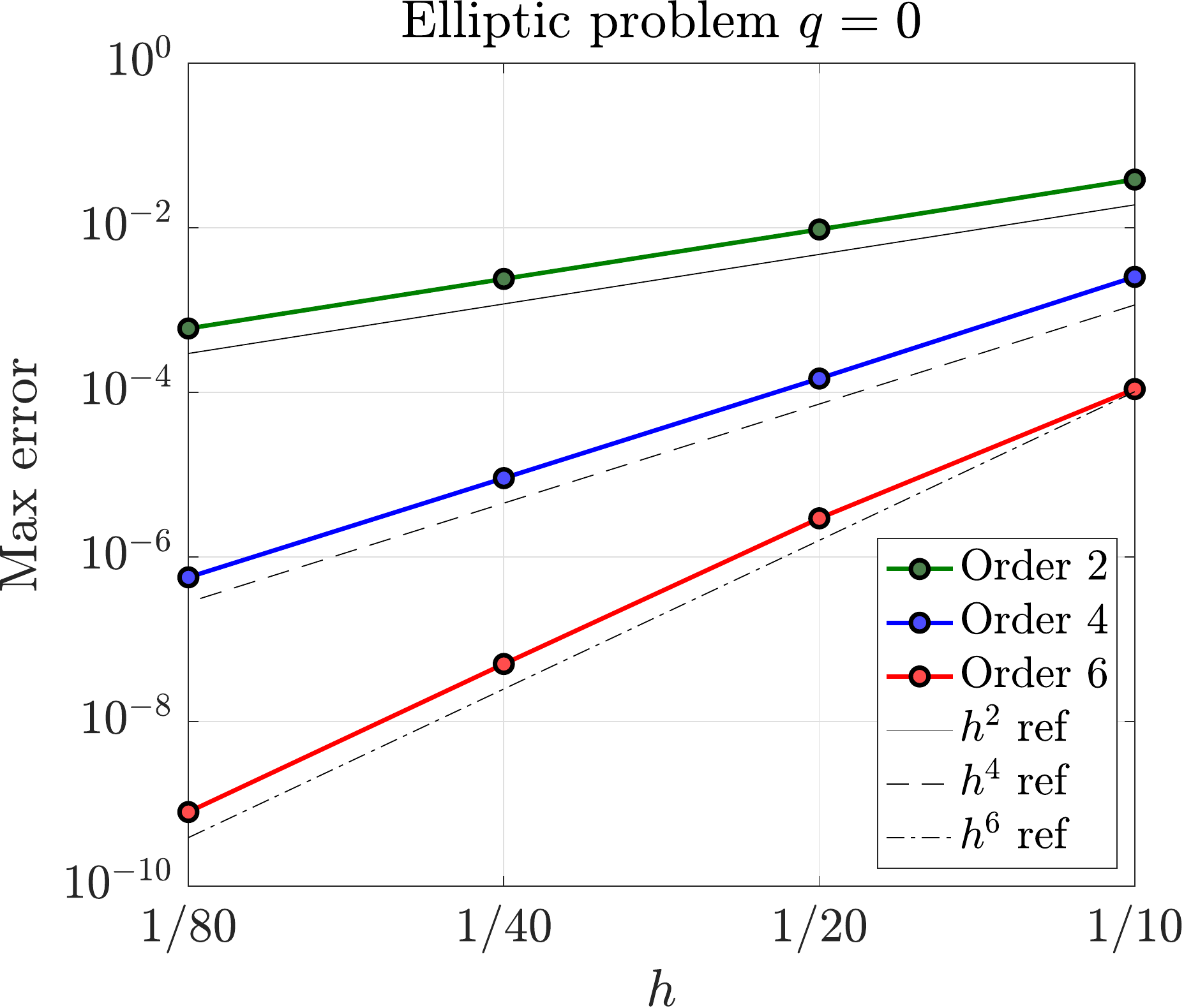}{\figWidth}};
				\draw ( 4.8,4.40) node[anchor=south west] {\trimfig{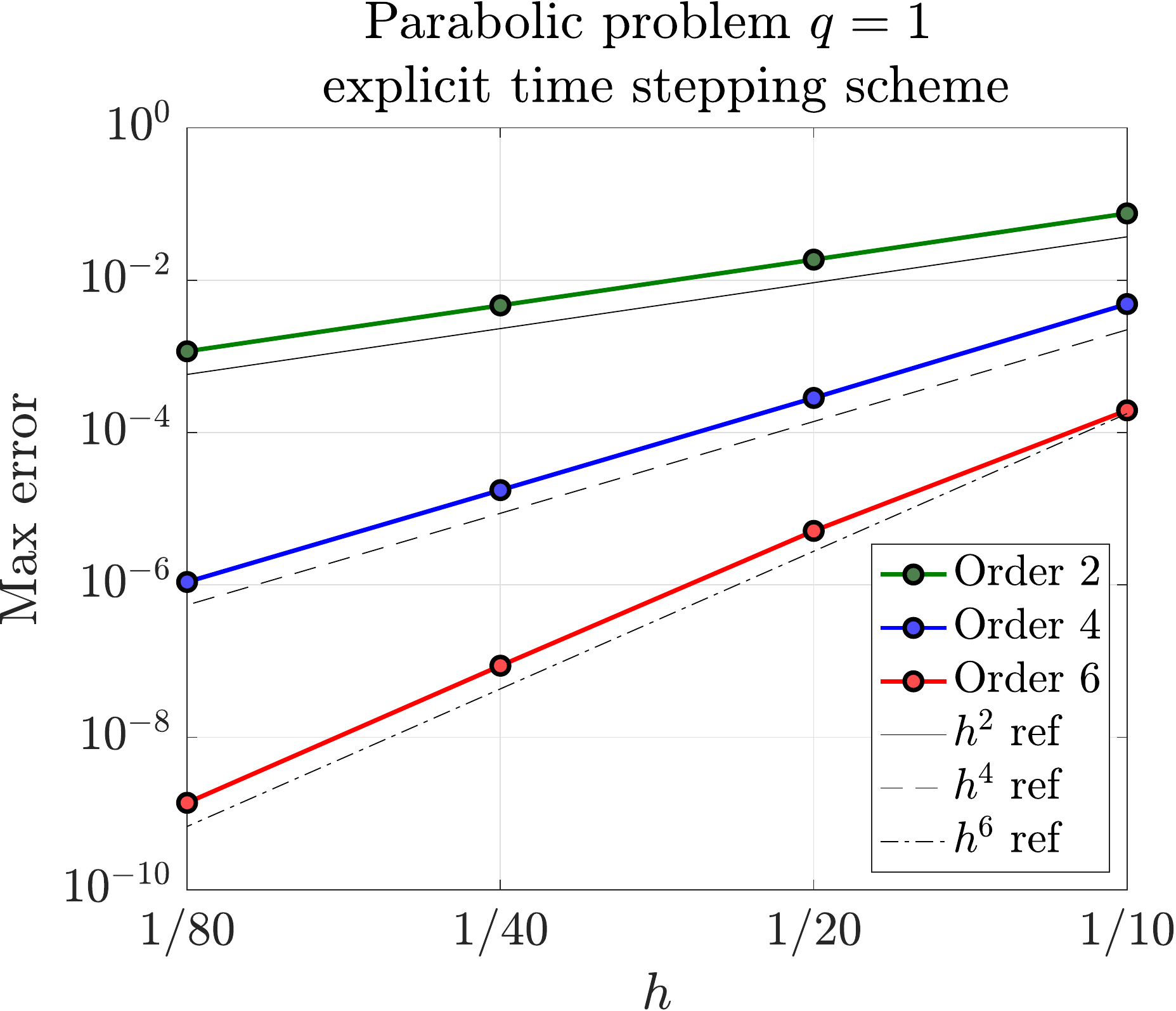}{\figWidth}};
				\draw ( -0.8,-0.65) node[anchor=south west] {\trimfig{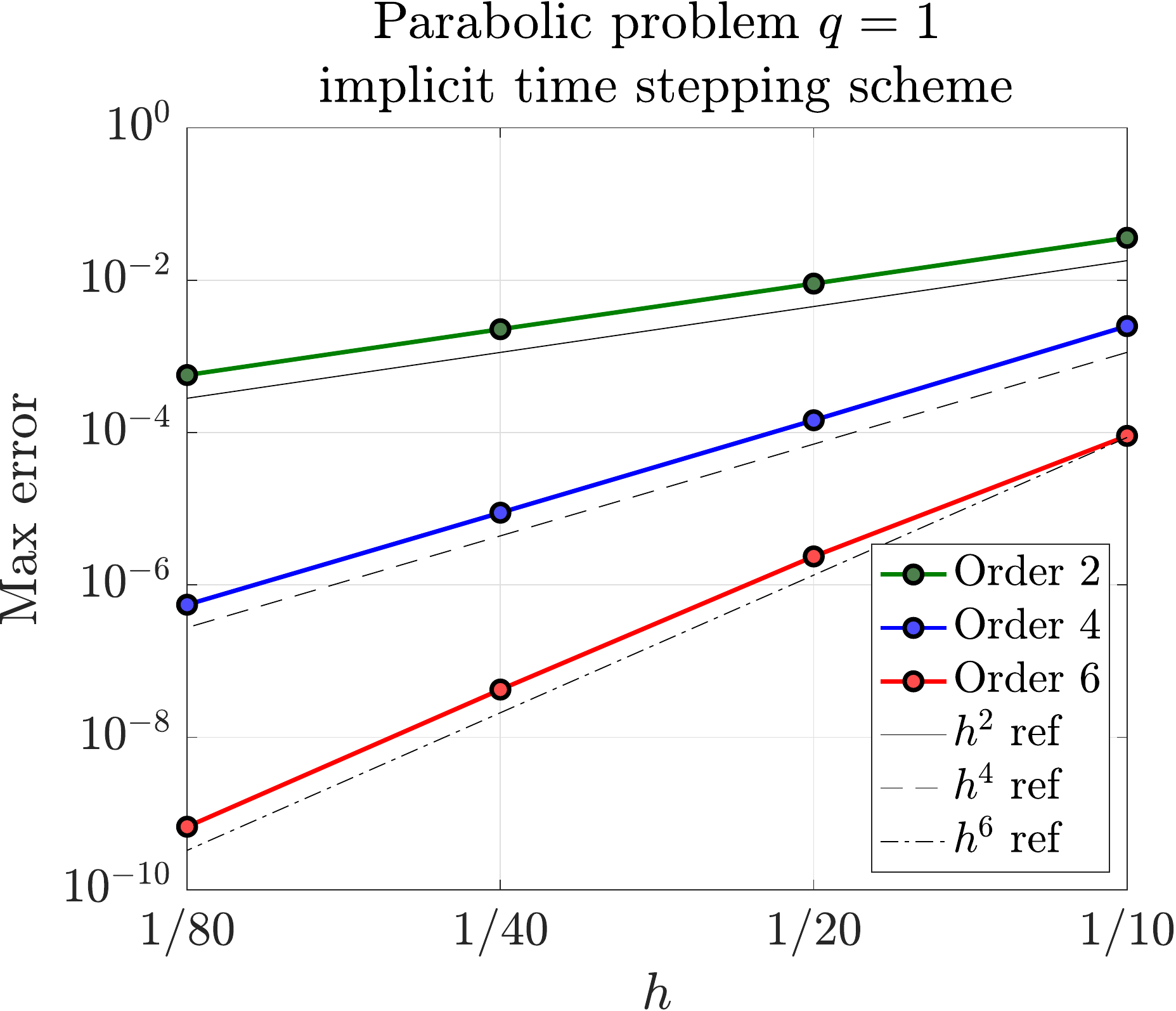}{\figWidth}};
				\draw ( 4.8,-0.65) node[anchor=south west] {\trimfig{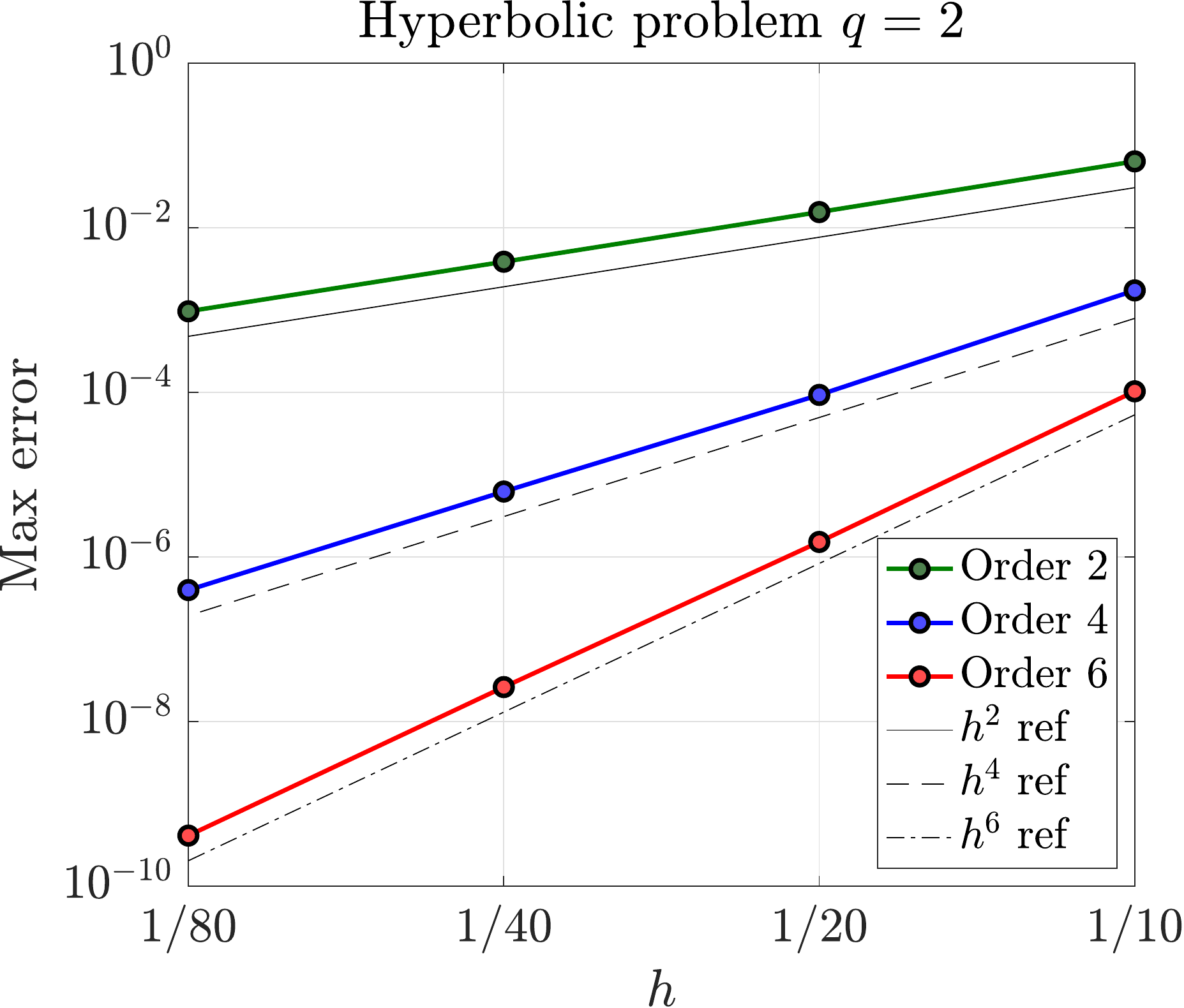}{\figWidth}};
			\end{tikzpicture}
		\end{center}
		\caption{Maximum errors versus grid spacing $h$ for Test 1 using manufactured solutions.  Fully discrete schemes for $q=0$ (upper left), $q=1$ and FE time-stepping (upper right), $q=1$ and BDF time-stepping (lower left), and $q=2$ (lower right).}
		\label{fig:MSexample1}
	\end{figure}
}

%% file: tex/MSexample2.tex

{
	\newcommand{\figWidth}{5.6cm}
	\newcommand{\trimfig}[2]{\trimw{#1}{#2}{0.}{0.}{.0}{.0}}  
	\newcommand{\trimfigb}[2]{\trimwb{#1}{#2}{.2}{0.1}{.2}{.15}} 
	\begin{figure}
		\begin{center}
			\begin{tikzpicture}[scale=1]
				\useasboundingbox (0,0) rectangle (11,9.5);  
				\draw ( -0.8,4.40) node[anchor=south west] {\trimfig{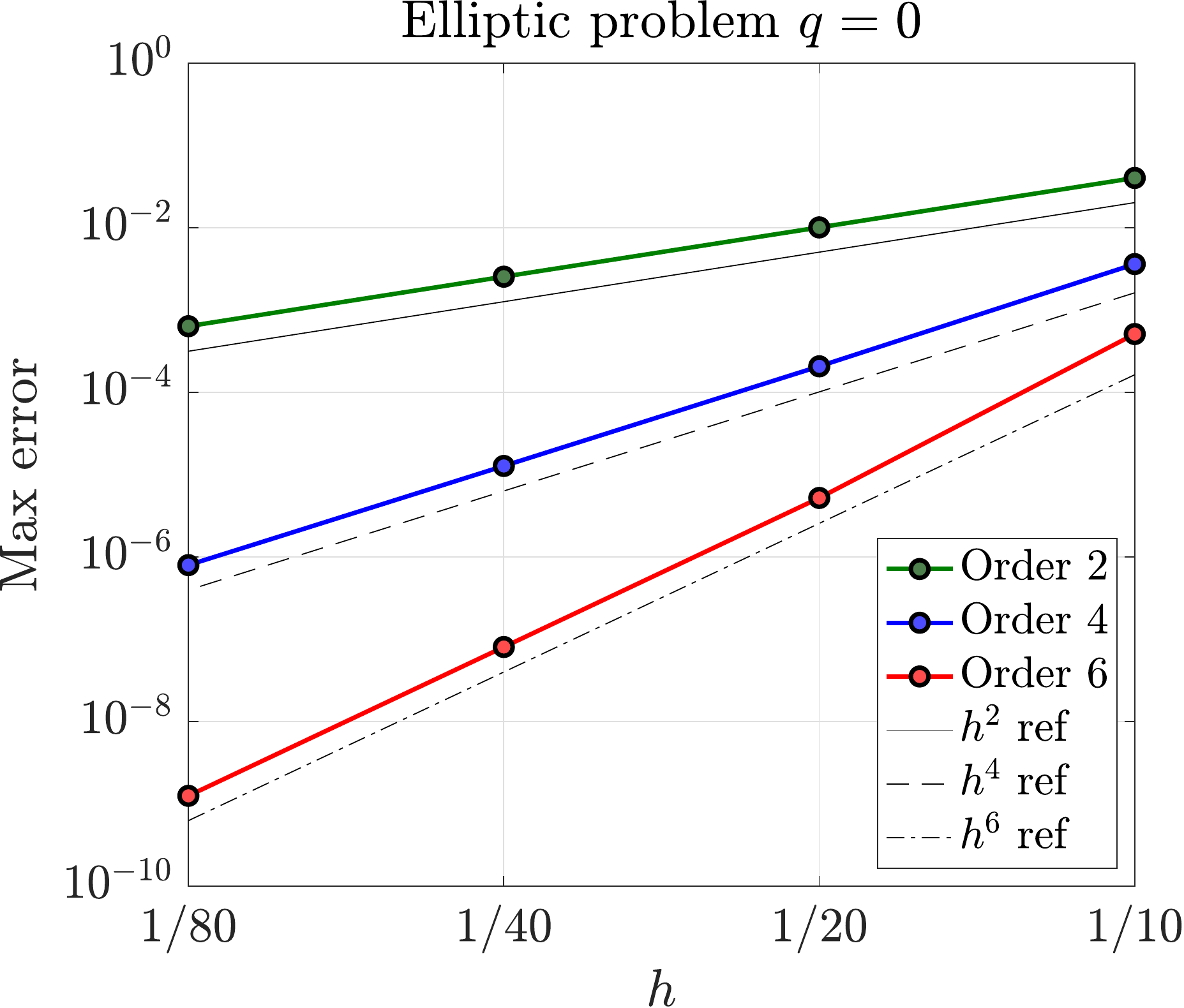}{\figWidth}};
				\draw ( 4.8,4.40) node[anchor=south west] {\trimfig{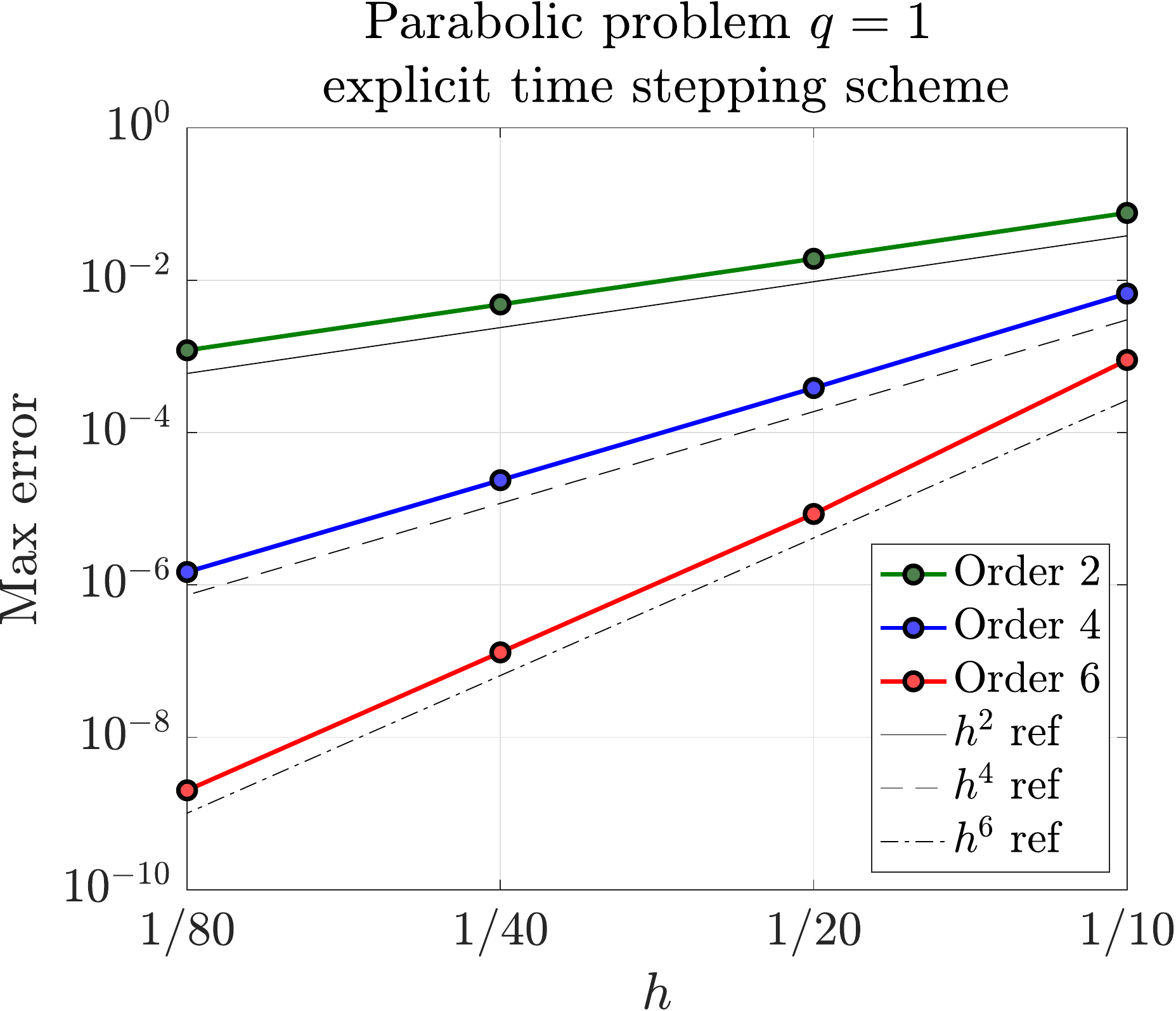}{\figWidth}};
				\draw ( -0.8,-0.65) node[anchor=south west] {\trimfig{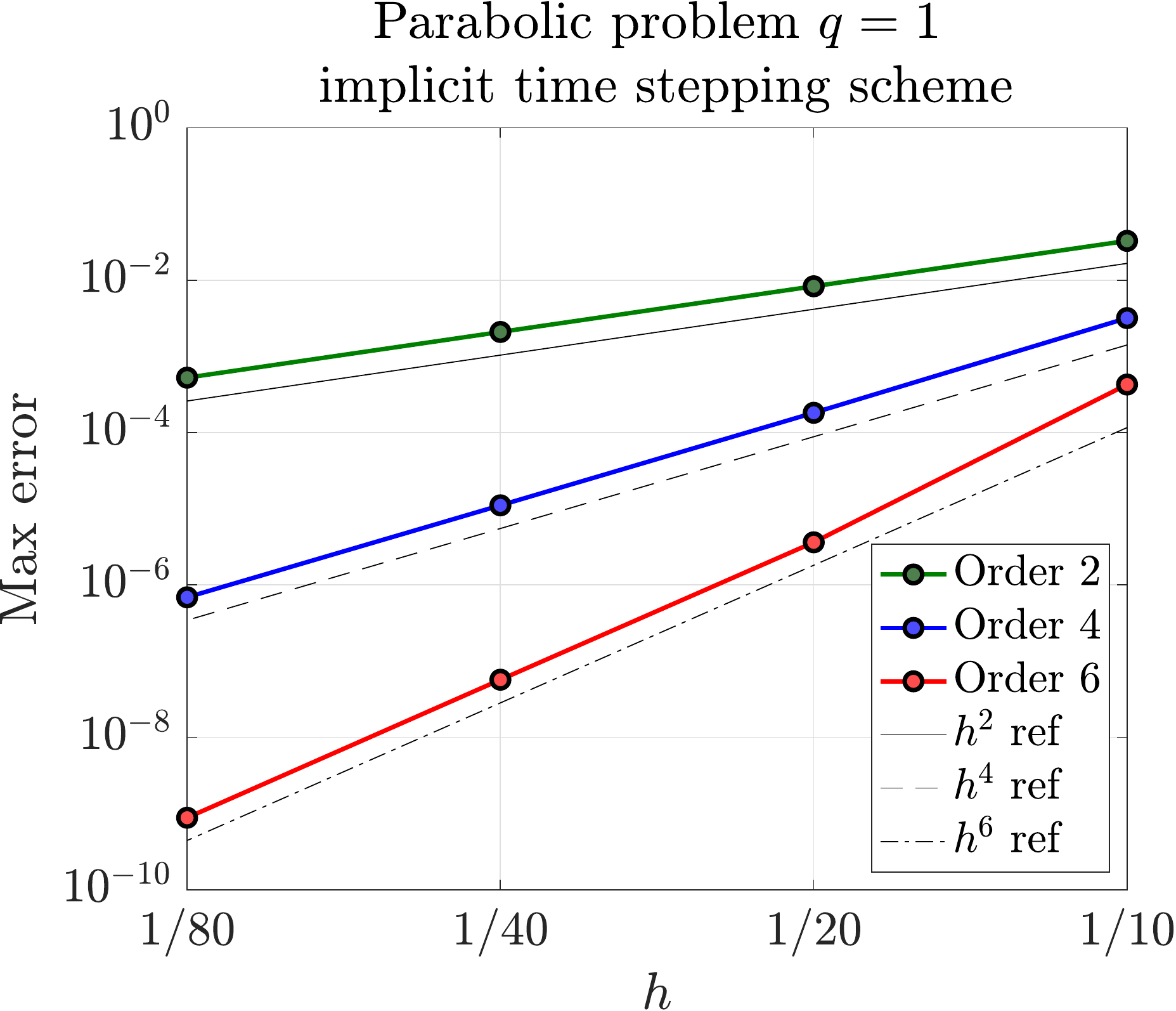}{\figWidth}};
				\draw ( 4.8,-0.65) node[anchor=south west] {\trimfig{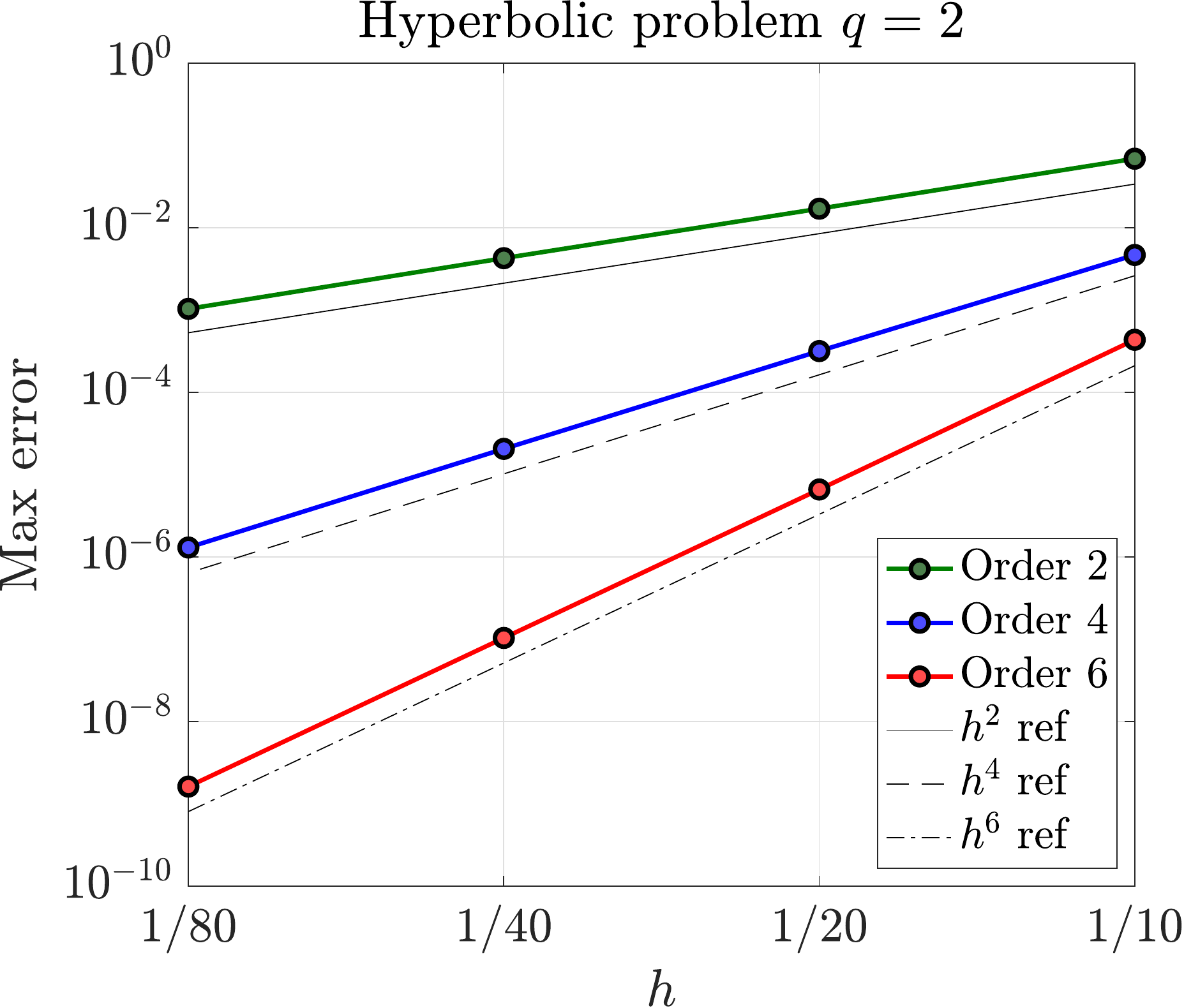}{\figWidth}};
			\end{tikzpicture}
		\end{center}
		\caption{Maximum errors versus grid spacing $h$ for Test 2 using manufactured solutions.  Fully discrete schemes for $q=0$ (upper left), $q=1$ and FE time-stepping (upper right), $q=1$ and BDF time-stepping (lower left), and $q=2$ (lower right).}
		\label{fig:MSexample2}
	\end{figure}
}

%% file: tex/MSexample3.tex

{
	\newcommand{\figWidth}{5.6cm}
	\newcommand{\trimfig}[2]{\trimw{#1}{#2}{0.}{0.}{.0}{.0}}  
	\newcommand{\trimfigb}[2]{\trimwb{#1}{#2}{.2}{0.1}{.2}{.15}} 
	\begin{figure}
		\begin{center}
			\begin{tikzpicture}[scale=1]
				\useasboundingbox (0,0) rectangle (11,9.5);  
				\draw ( -0.8,4.40) node[anchor=south west] {\trimfig{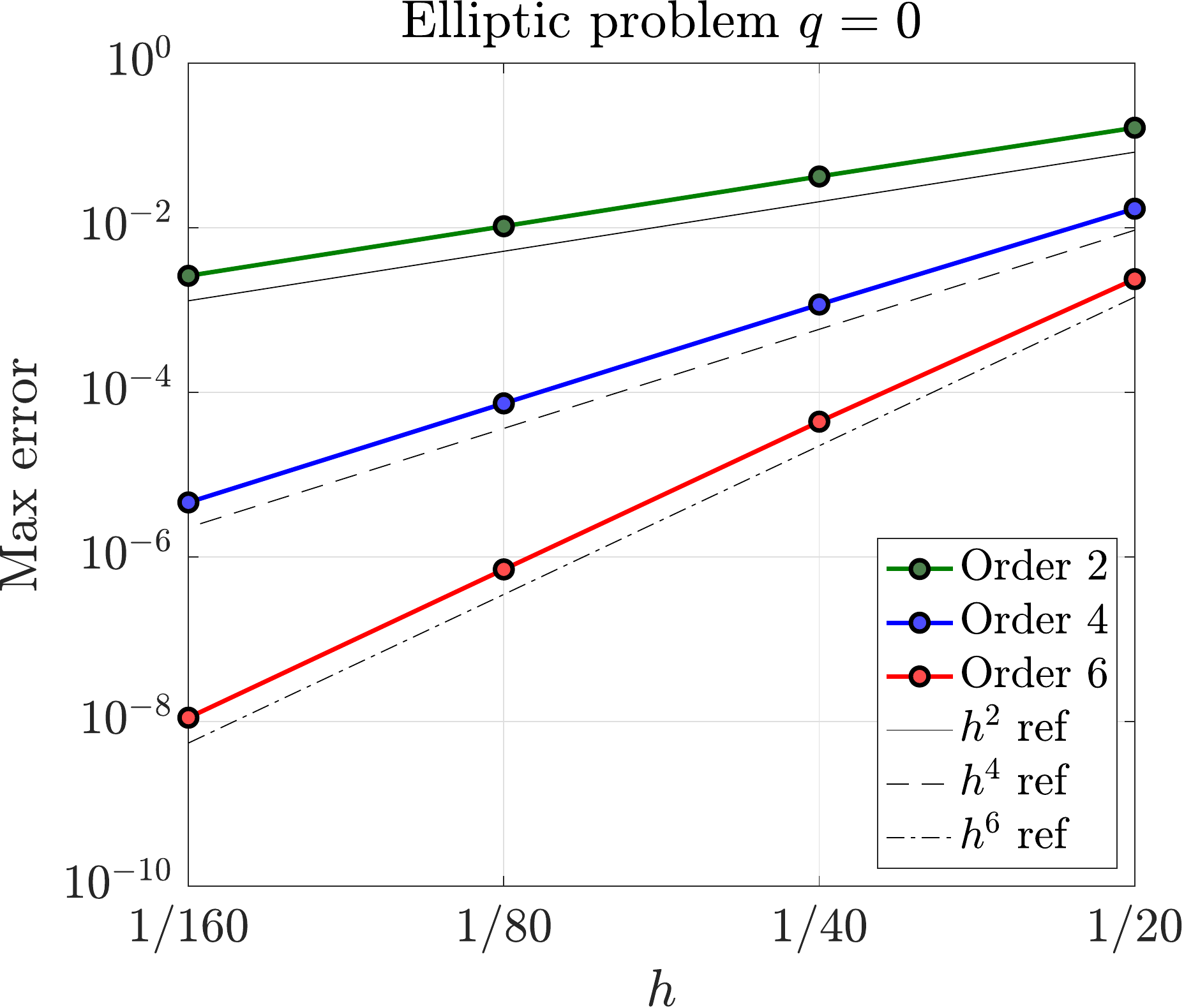}{\figWidth}};
				\draw ( 4.8,4.40) node[anchor=south west] {\trimfig{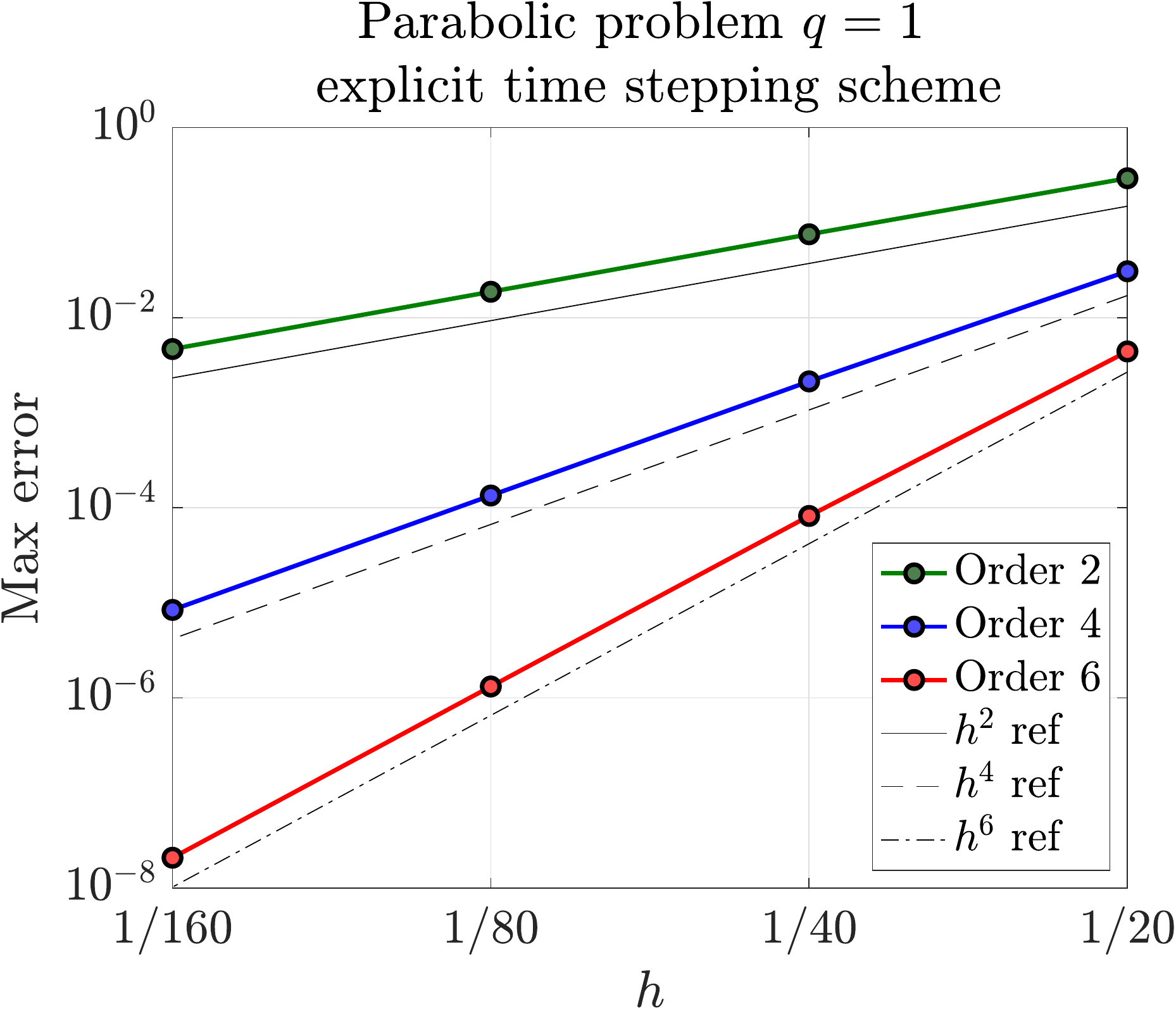}{\figWidth}};
				\draw ( -0.8,-0.65) node[anchor=south west] {\trimfig{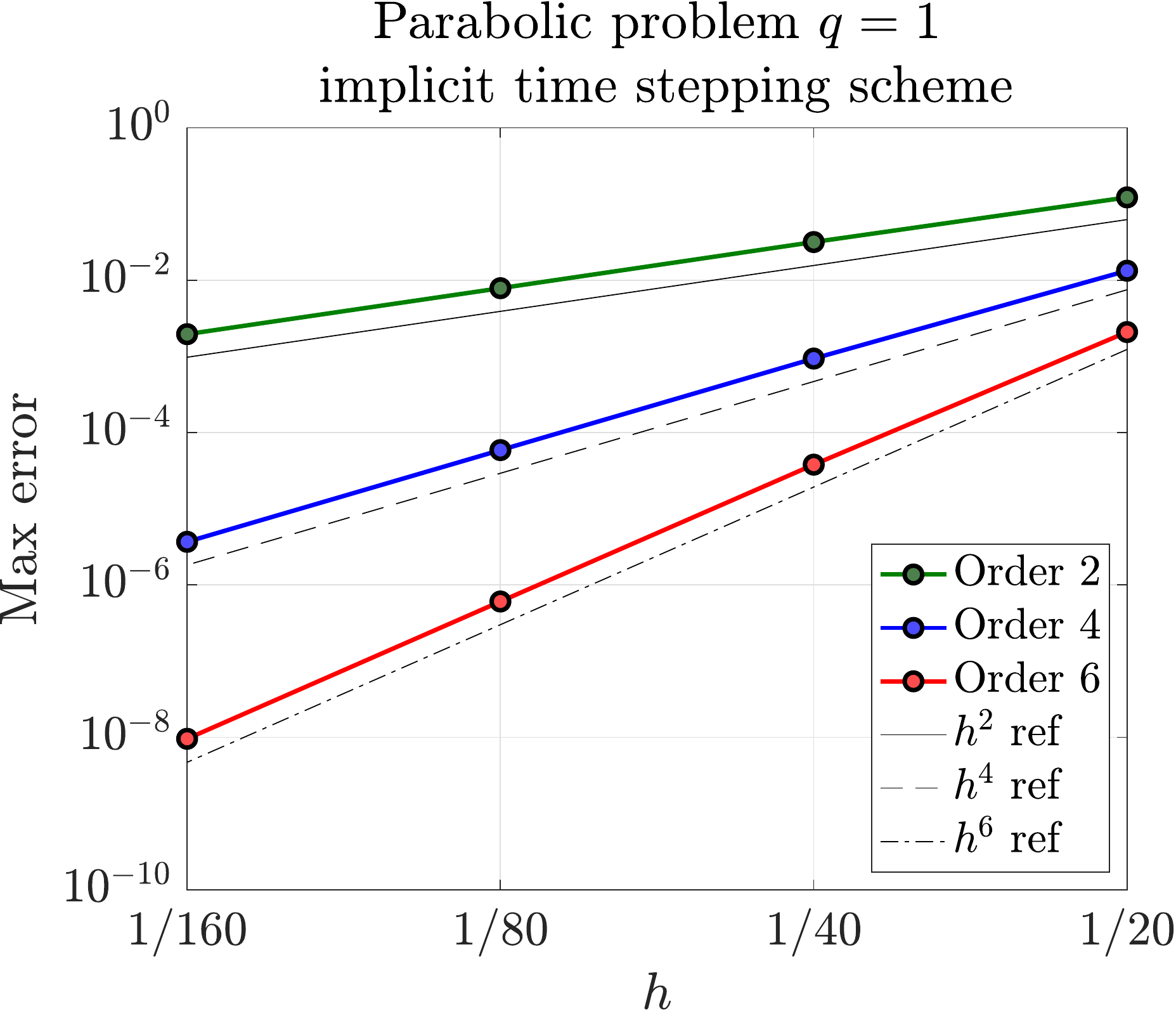}{\figWidth}};
				\draw ( 4.8,-0.65) node[anchor=south west] {\trimfig{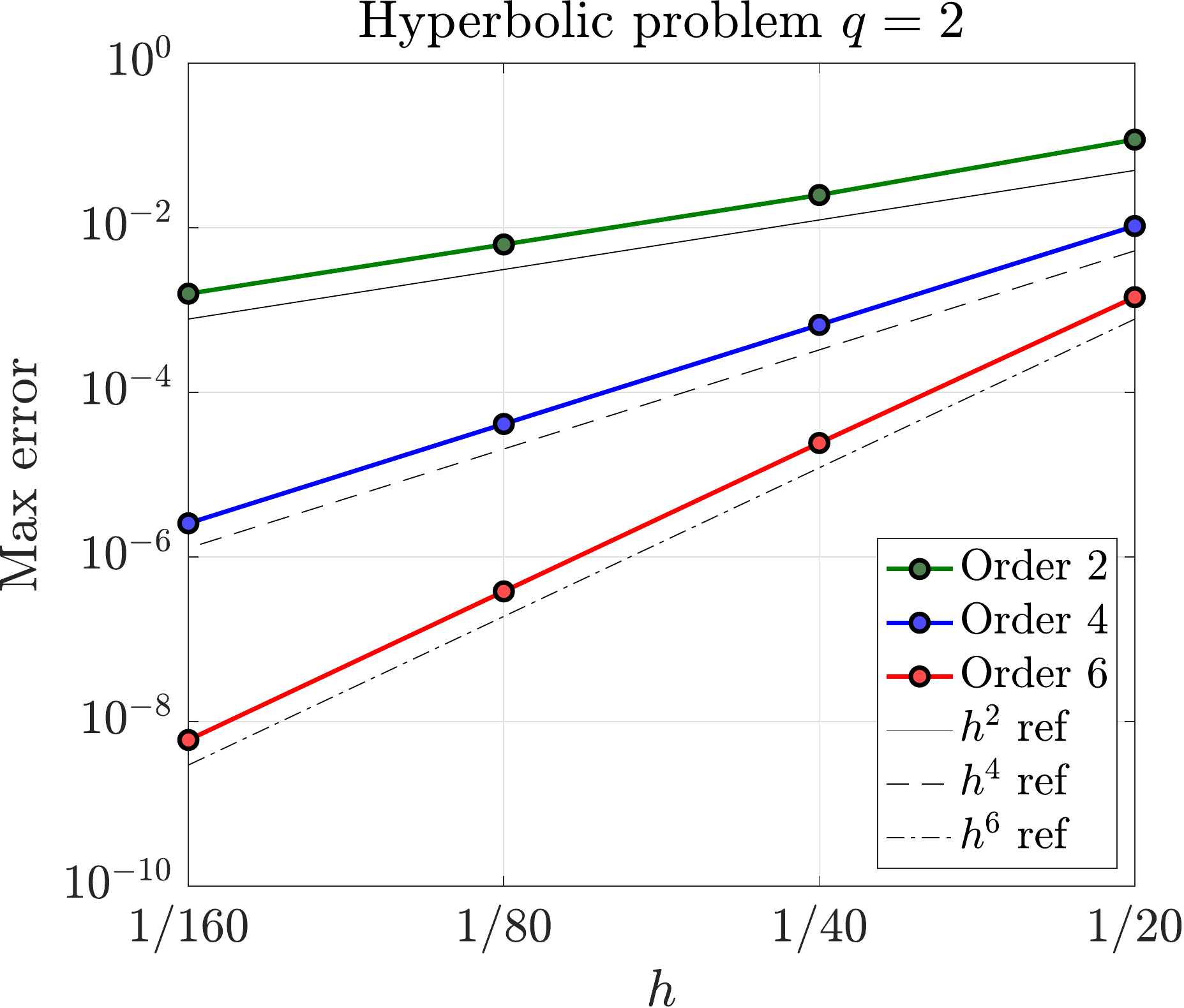}{\figWidth}};
			\end{tikzpicture}
		\end{center}
		\caption{Maximum errors versus an approximate grid spacing $h$ for Test 3 using manufactured solutions.  Fully discrete schemes for $q=0$ (upper left), $q=1$ and FE time-stepping (upper right), $q=1$ and BDF time-stepping (lower left), and $q=2$ (lower right).}
		\label{fig:MSexample3}
	\end{figure}
}

%% file: tex/Example4.tex


{
	\newcommand{\figWidth}{5.6cm}
	\newcommand{\trimfig}[2]{\trimw{#1}{#2}{0.}{0.}{.0}{.0}}  
	\newcommand{\trimfigb}[2]{\trimw{#1}{#2}{0}{0}{.2}{.2}} 
	\begin{figure}
		\begin{center}
			\begin{tikzpicture}[scale=1]
				\useasboundingbox (0,0) rectangle (11,11);  
				\draw ( -1.3,4.5) node[anchor=south west] {\trimfig{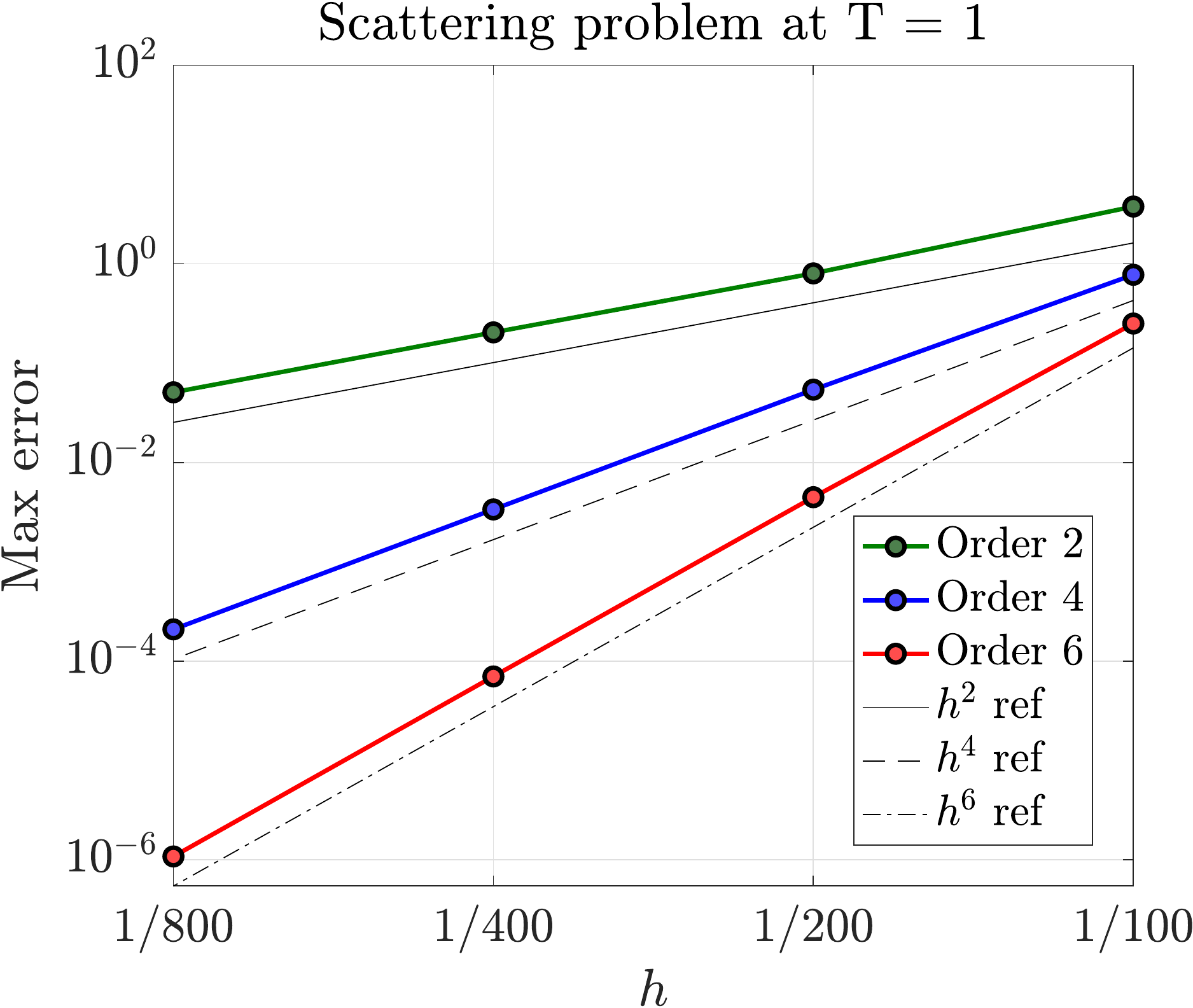}{\figWidth}};
				\draw ( -1.3,.5) node[anchor=south west] {\trimfig{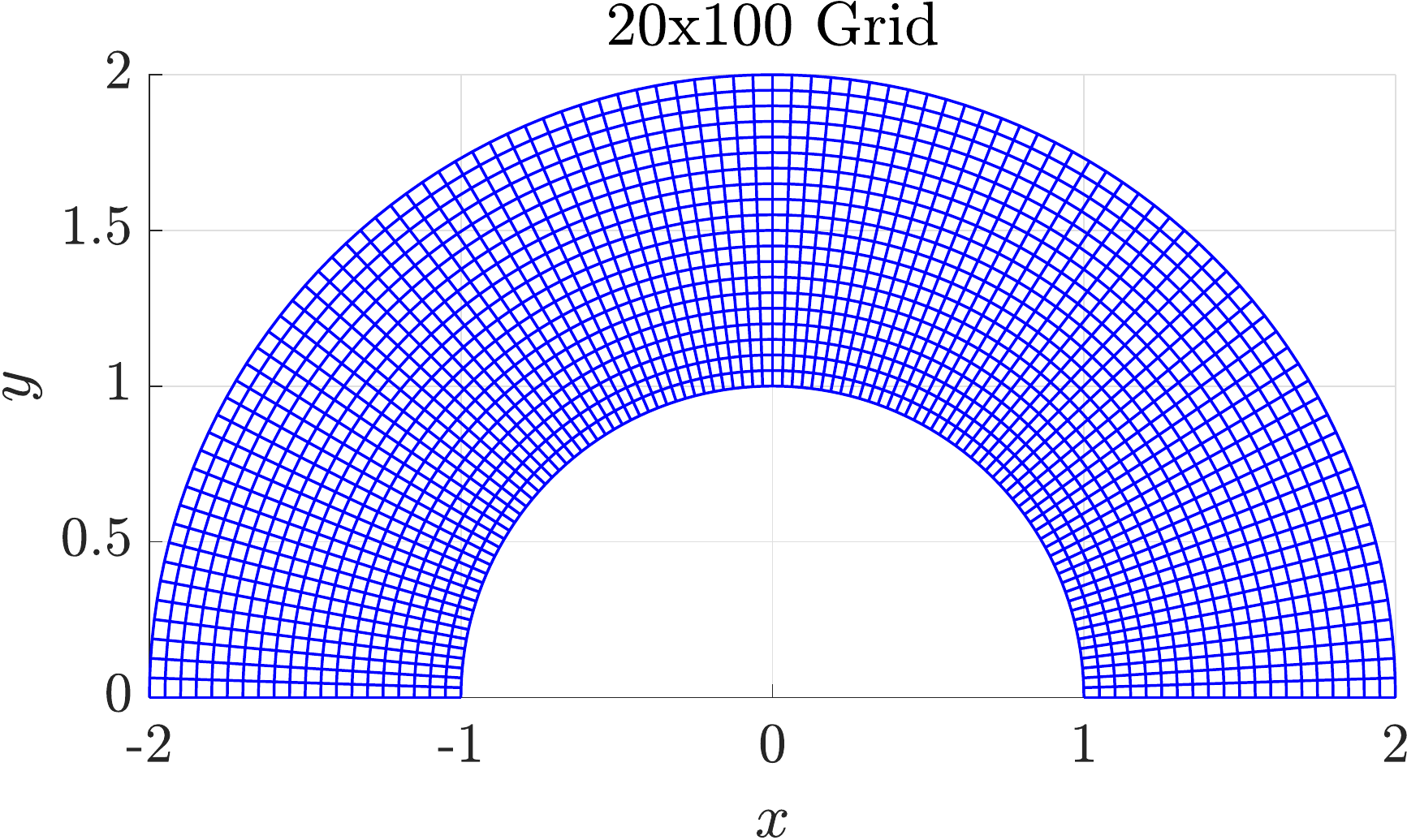}{5.7cm}};
				\draw ( 4.6,6.8) node[anchor=south west] {\trimfigb{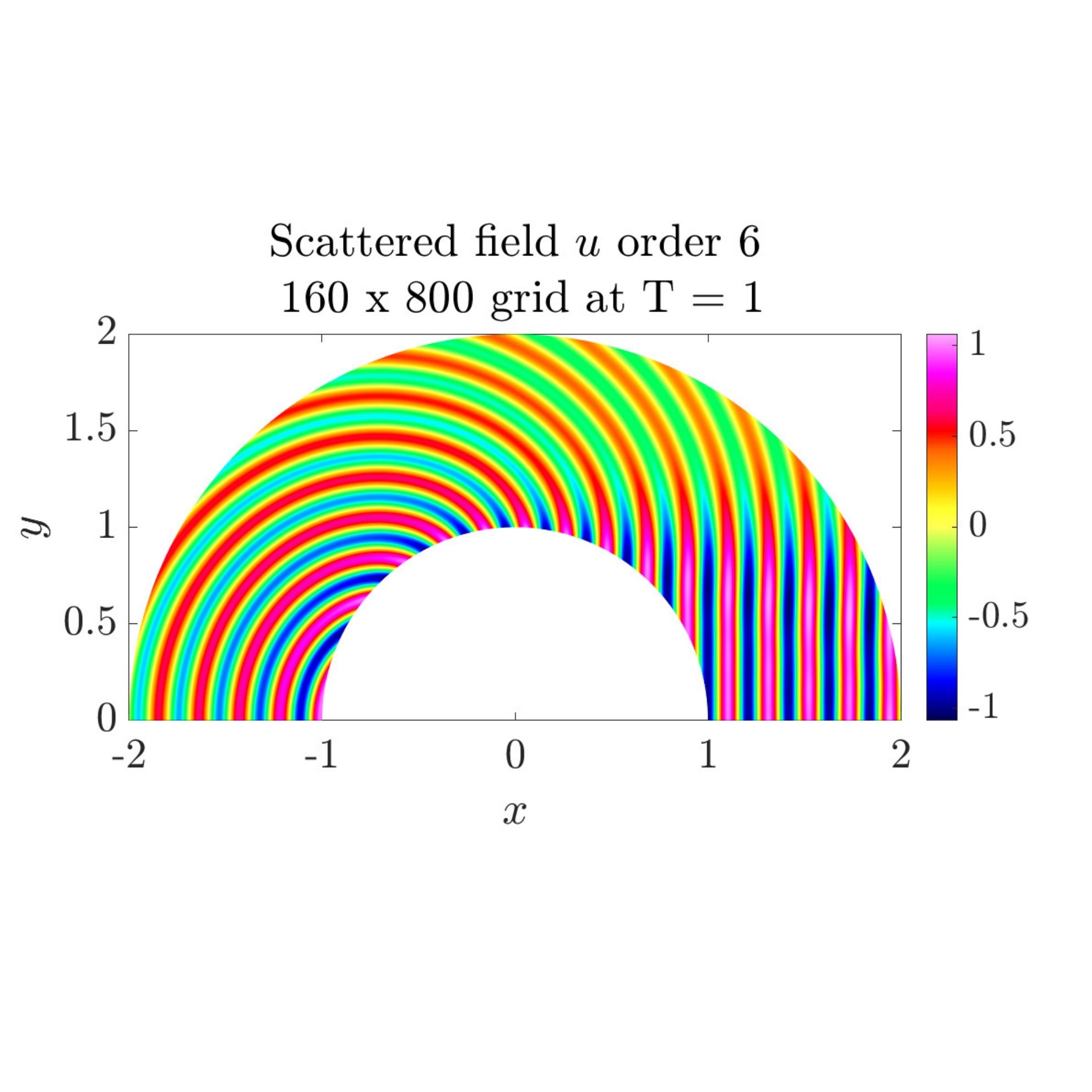}{6.5cm}};
				\draw ( 4.6,3) node[anchor=south west] {\trimfigb{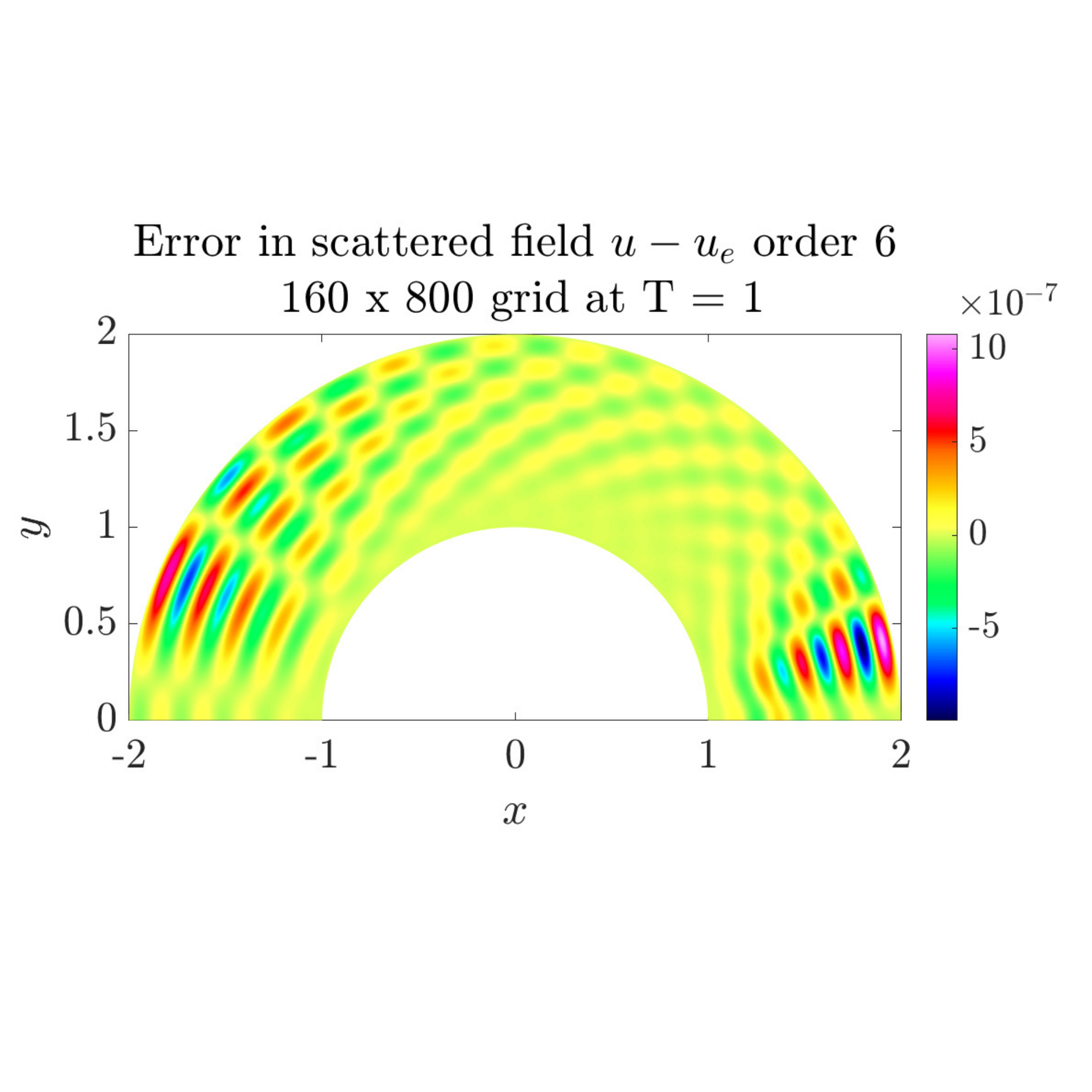}{6.5cm}};
				\draw ( 4.6,-.7) node[anchor=south west] {\trimfigb{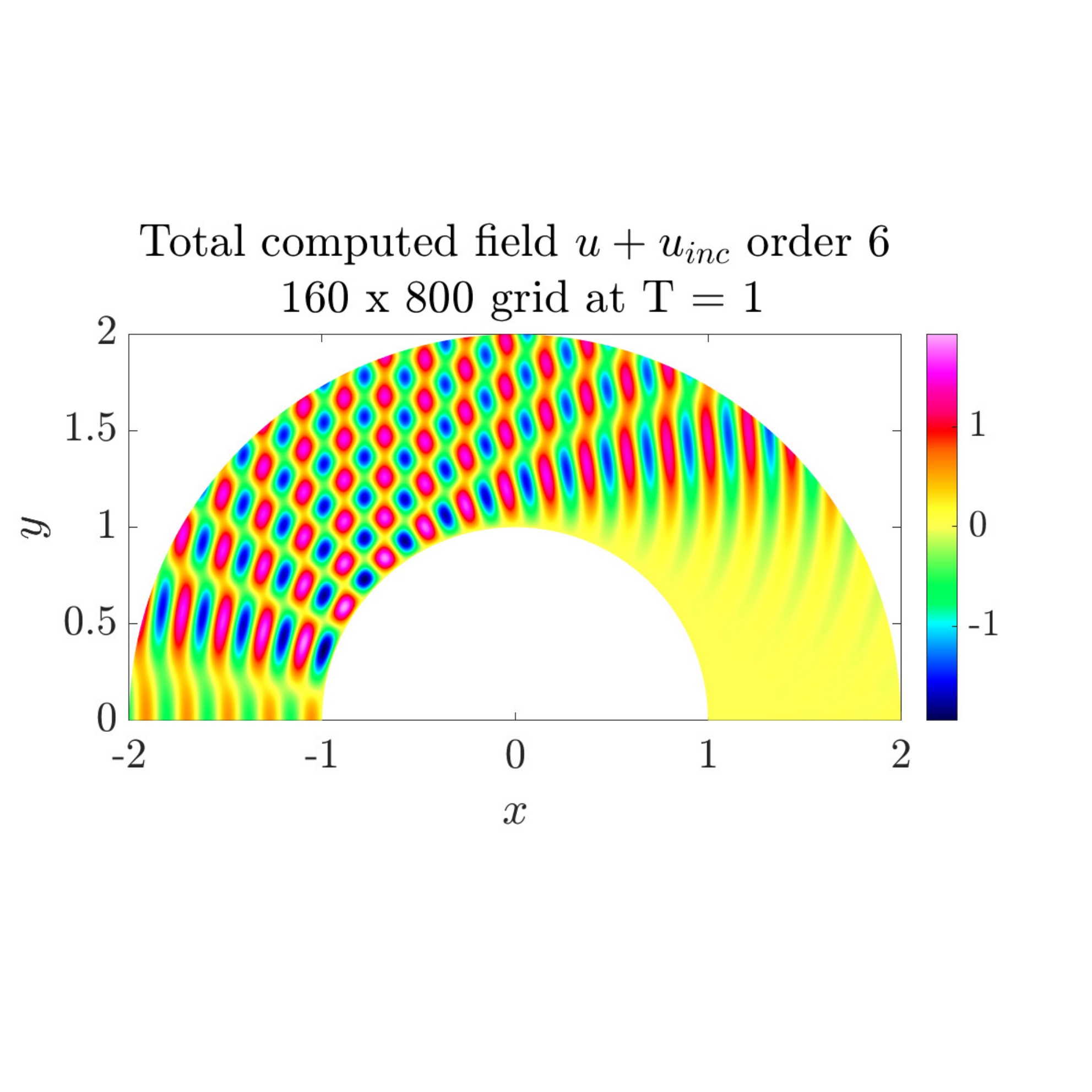}{6.3cm}};
			\end{tikzpicture}
		\end{center}
		\caption{Plane-wave scattering from a cylinder.  Maximum errors at $T=1$ for solutions computed using the explicit time-stepping schemes with $d=2$, $4$ and $6$ (upper left) and the coarsest grid for $h=1/100$ (lower left).  Right column shows the scattered field (top), error in the scattered field (middle) and the total field (bottom) at $T=1$ computed using the sixth-order accurate scheme on the finest grid.}
		\label{fig:example4}
	\end{figure}
}

%% file: tex/Example5.tex

{
	\newcommand{\figWidth}{6.0cm}
	\newcommand{\trimfig}[2]{\trimw{#1}{#2}{0.}{0.}{.0}{.0}}  
	\newcommand{\trimfigb}[2]{\trimwb{#1}{#2}{.2}{0.1}{.2}{.15}} 
	\begin{figure}[H]
		\begin{center}
			\begin{tikzpicture}[scale=1]
				\useasboundingbox (0,0) rectangle (10,11);  
				\draw ( -1.5,5.4) node[anchor=south west] {\trimfig{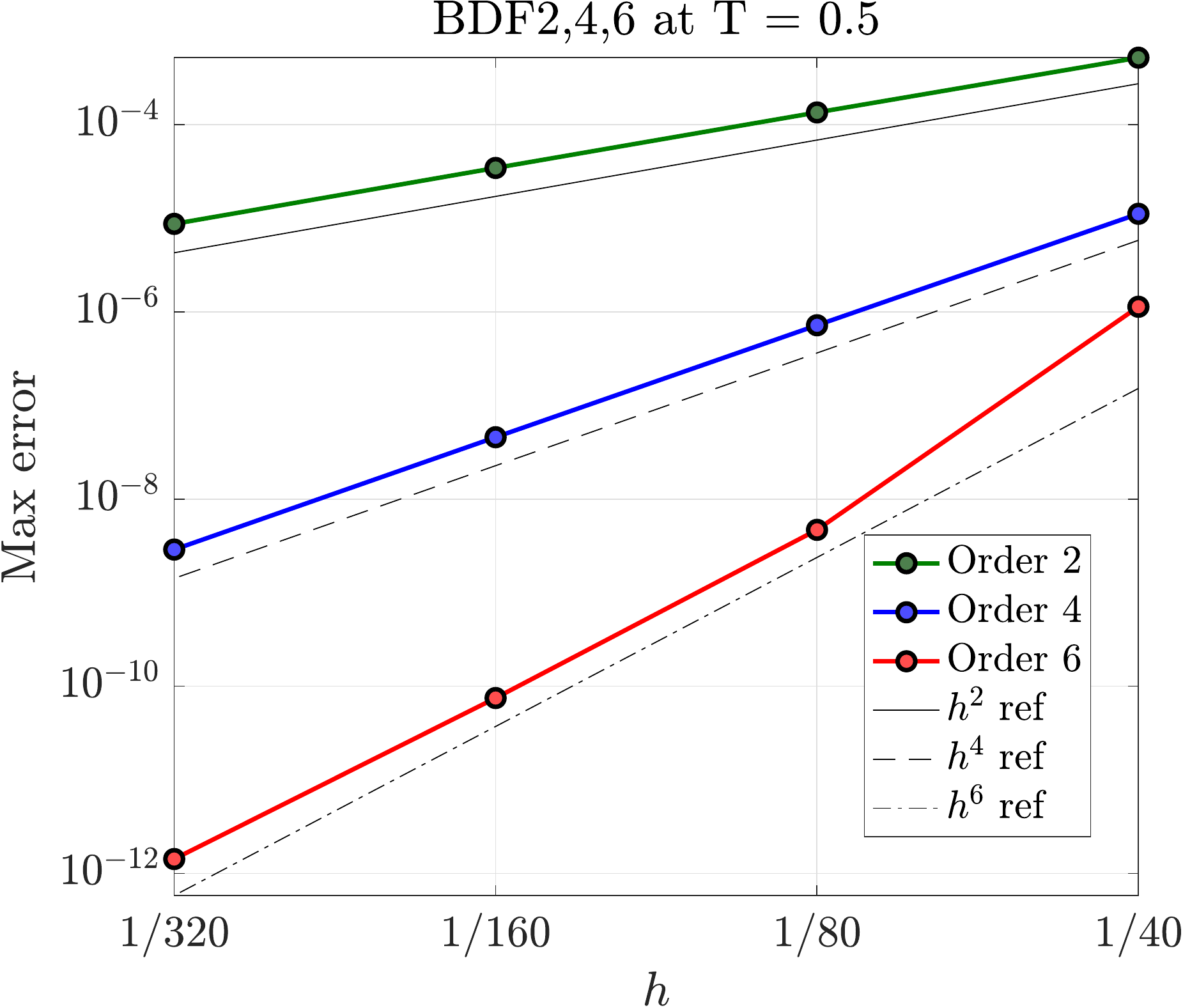}{\figWidth}};
				\draw ( 4.7,5) node[anchor=south west] {\trimfig{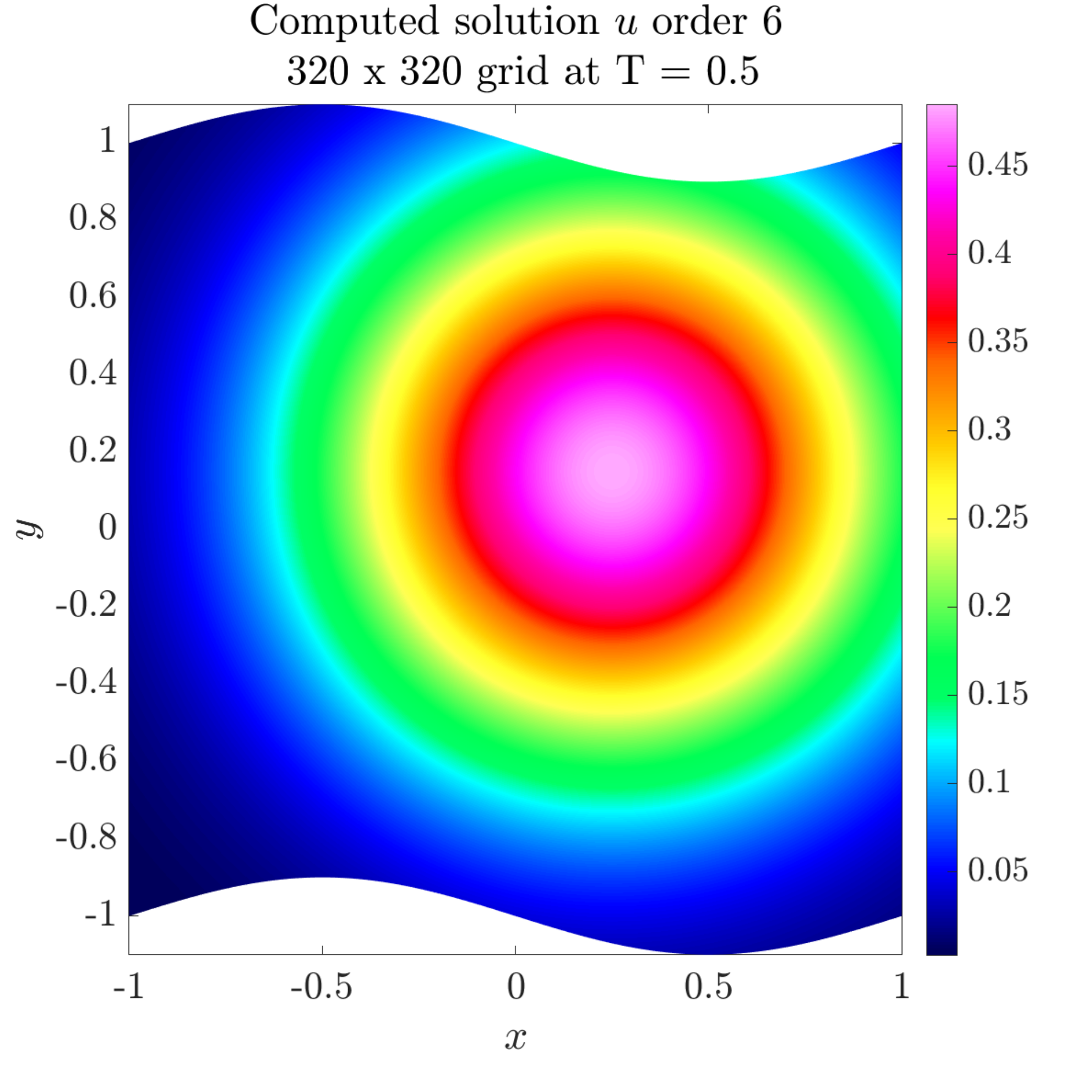}{\figWidth}};
				\draw ( 4.7,-.8) node[anchor=south west] {\trimfig{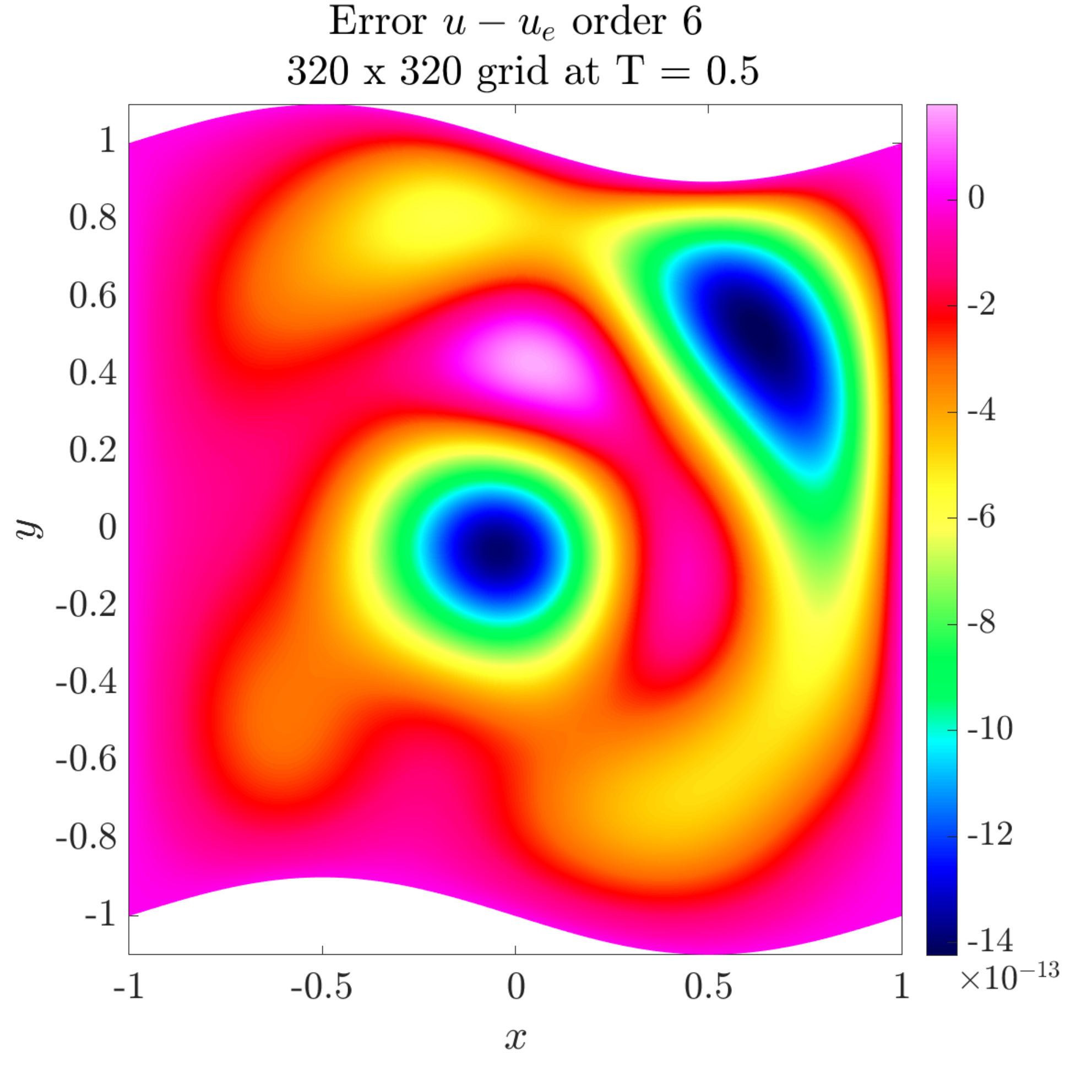}{\figWidth}};
				\draw ( -.7,-.5) node[anchor=south west] {\trimfig{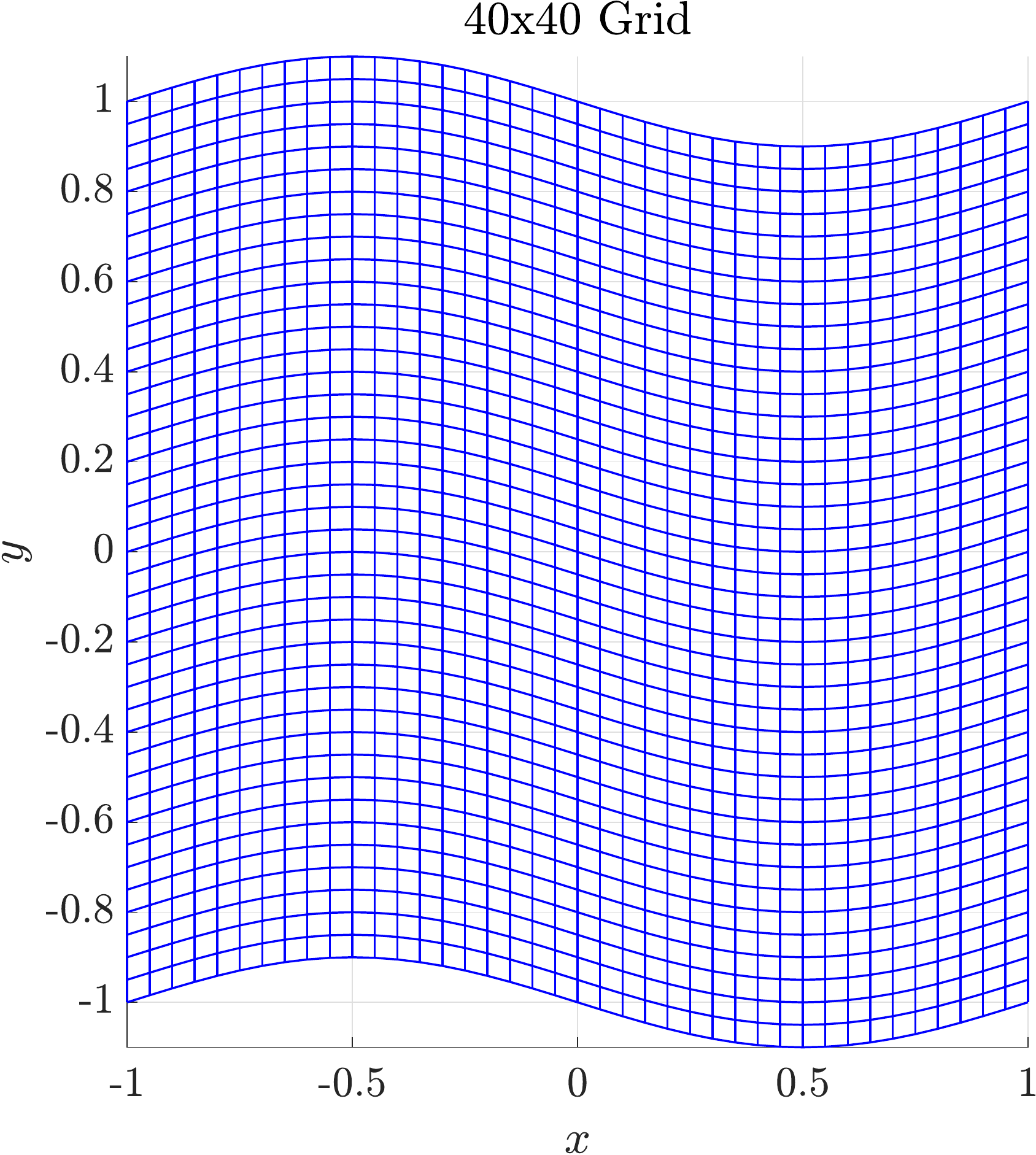}{4.8cm}};
			\end{tikzpicture}
		\end{center}
		\caption{Heat flow in a wavy channel.  Maximum errors at $T=0.5$ for solutions computed using the BDF time-stepping schemes with $d=2$, $4$ and $6$ (upper left) and the coarsest grid for $h=1/40$ (lower left).  Right column shows the temperature (top) and its error (bottom) at $T=0.5$ computed using the sixth-order accurate scheme on the finest grid.}
		\label{fig:example5}
	\end{figure}
}

%% file: tex/Example6.tex

{
	\newcommand{\figWidth}{5.6cm}
	\newcommand{\trimfig}[2]{\trimw{#1}{#2}{0.}{0.}{.0}{.0}}  
	\newcommand{\trimfigb}[2]{\trimwb{#1}{#2}{.2}{0.1}{.2}{.15}} 
	\begin{figure}[H]
		\begin{center}
			\begin{tikzpicture}[scale=1]
				\useasboundingbox (0,0) rectangle (18,10.5);  
				\draw ( -0.5,4.9) node[anchor=south west] {\trimfig{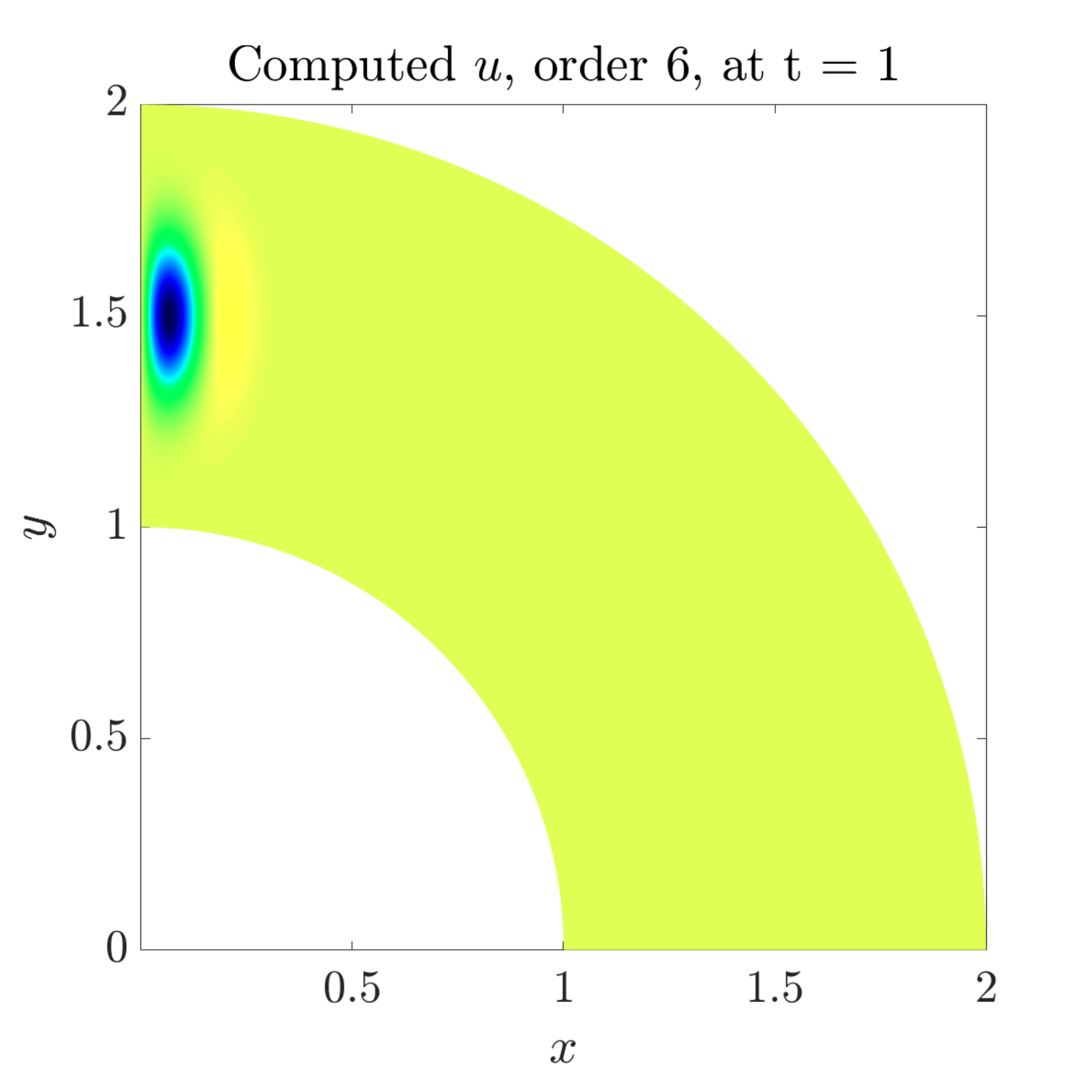}{\figWidth}};
				\draw ( 4.7,4.9) node[anchor=south west] {\trimfig{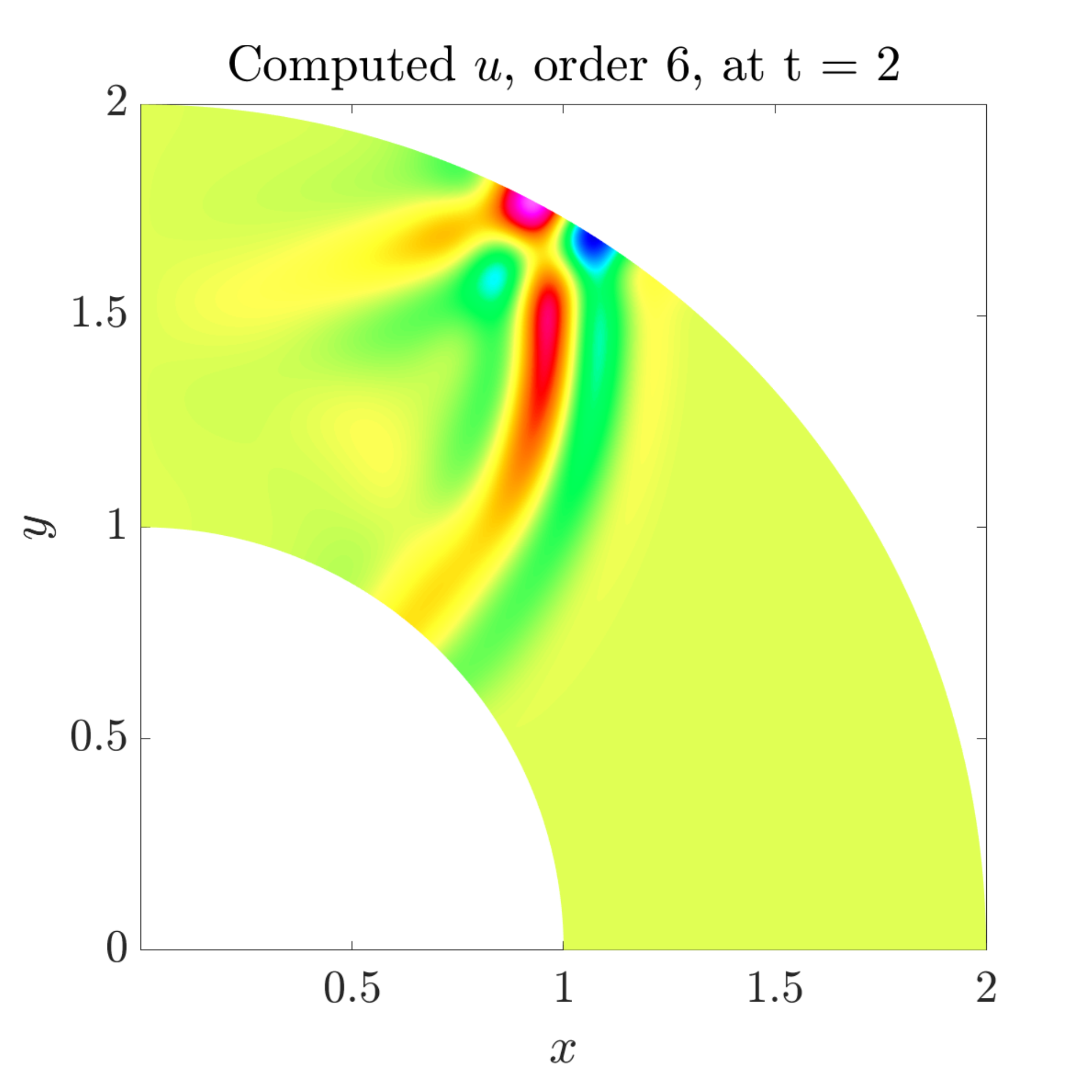}{\figWidth}};
				\draw ( 9.9,4.9) node[anchor=south west] {\trimfig{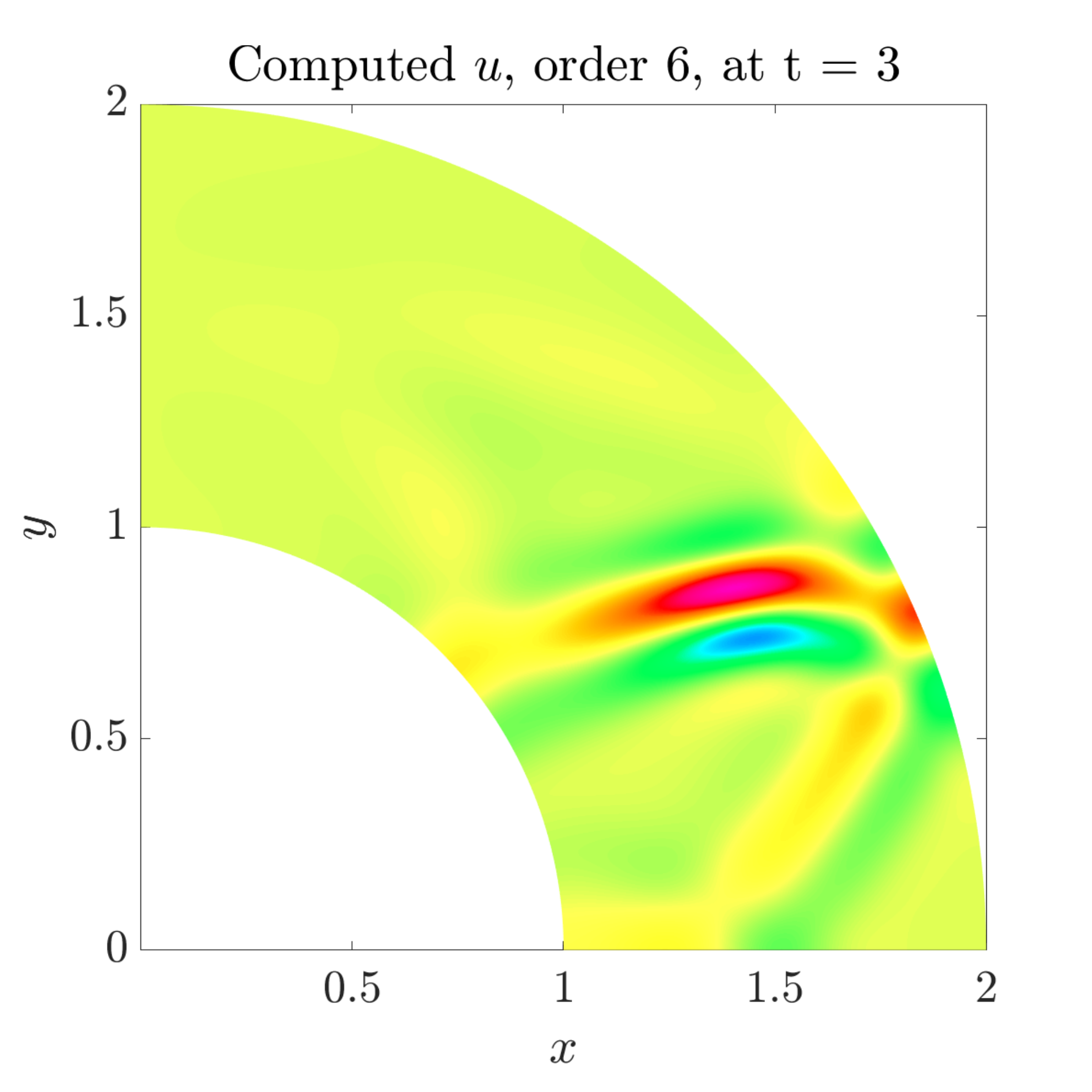}{\figWidth}};
				\draw ( 15,5.61) node[anchor=south west] {\trimfig{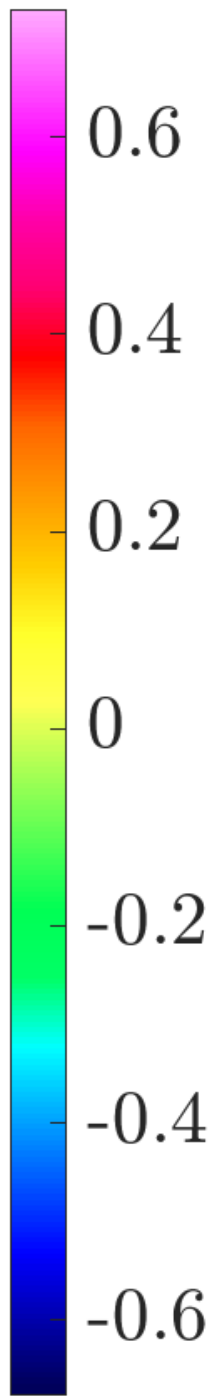}{.665cm}};
				\draw ( -.5,-.7) node[anchor=south west] {\trimfig{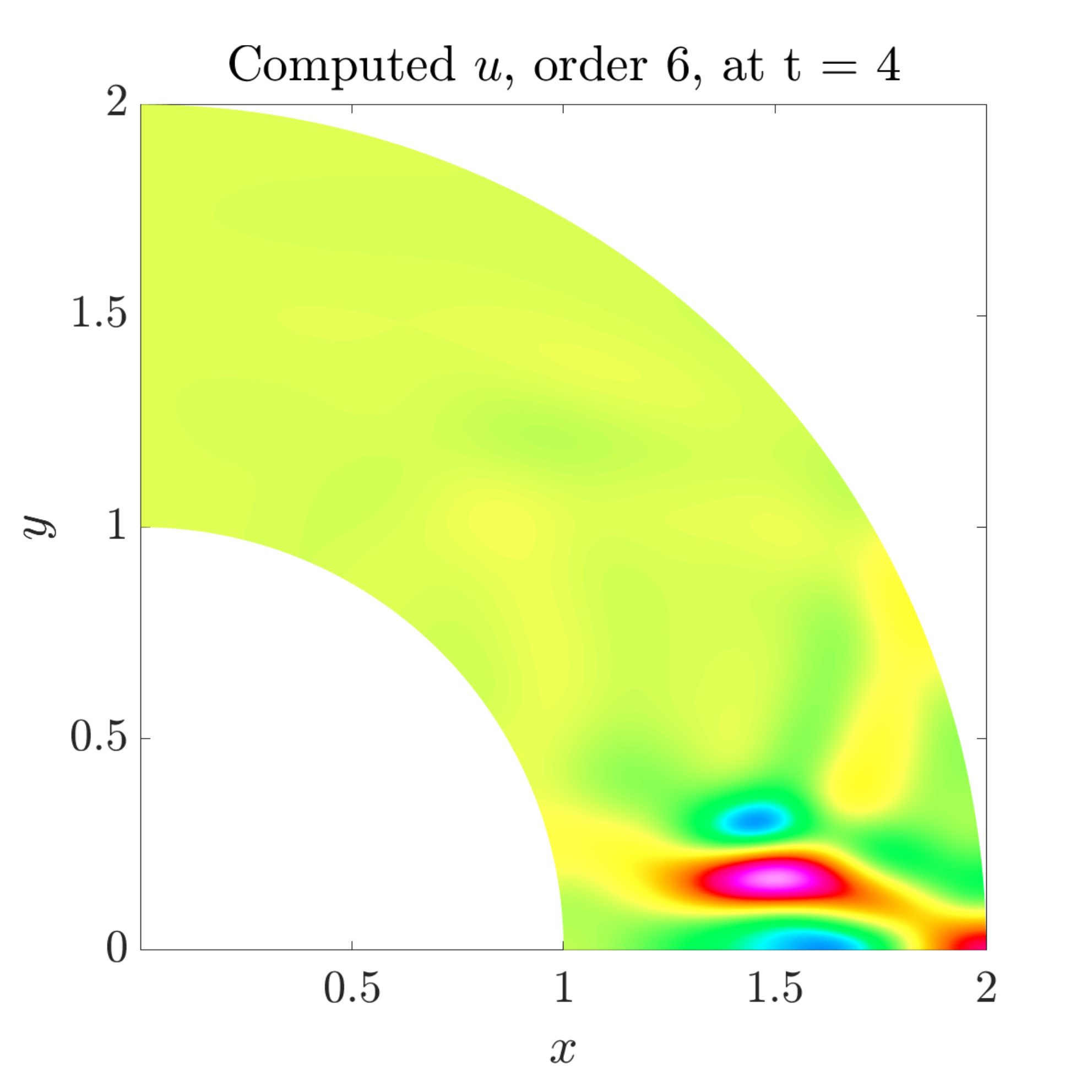}{\figWidth}};
				\draw ( 4.7,-.7) node[anchor=south west] {\trimfig{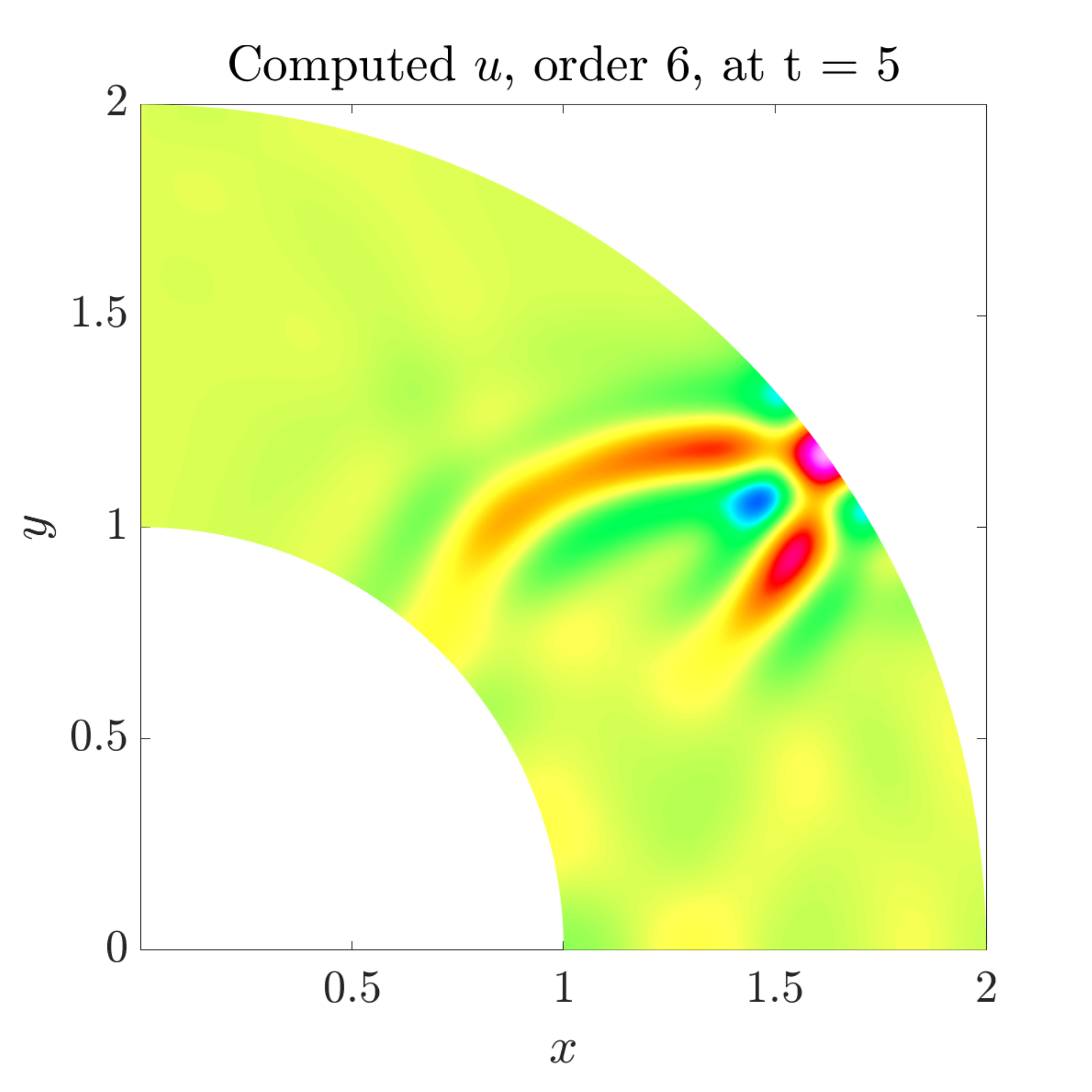}{\figWidth}};
				\draw ( 9.9,-.7) node[anchor=south west] {\trimfig{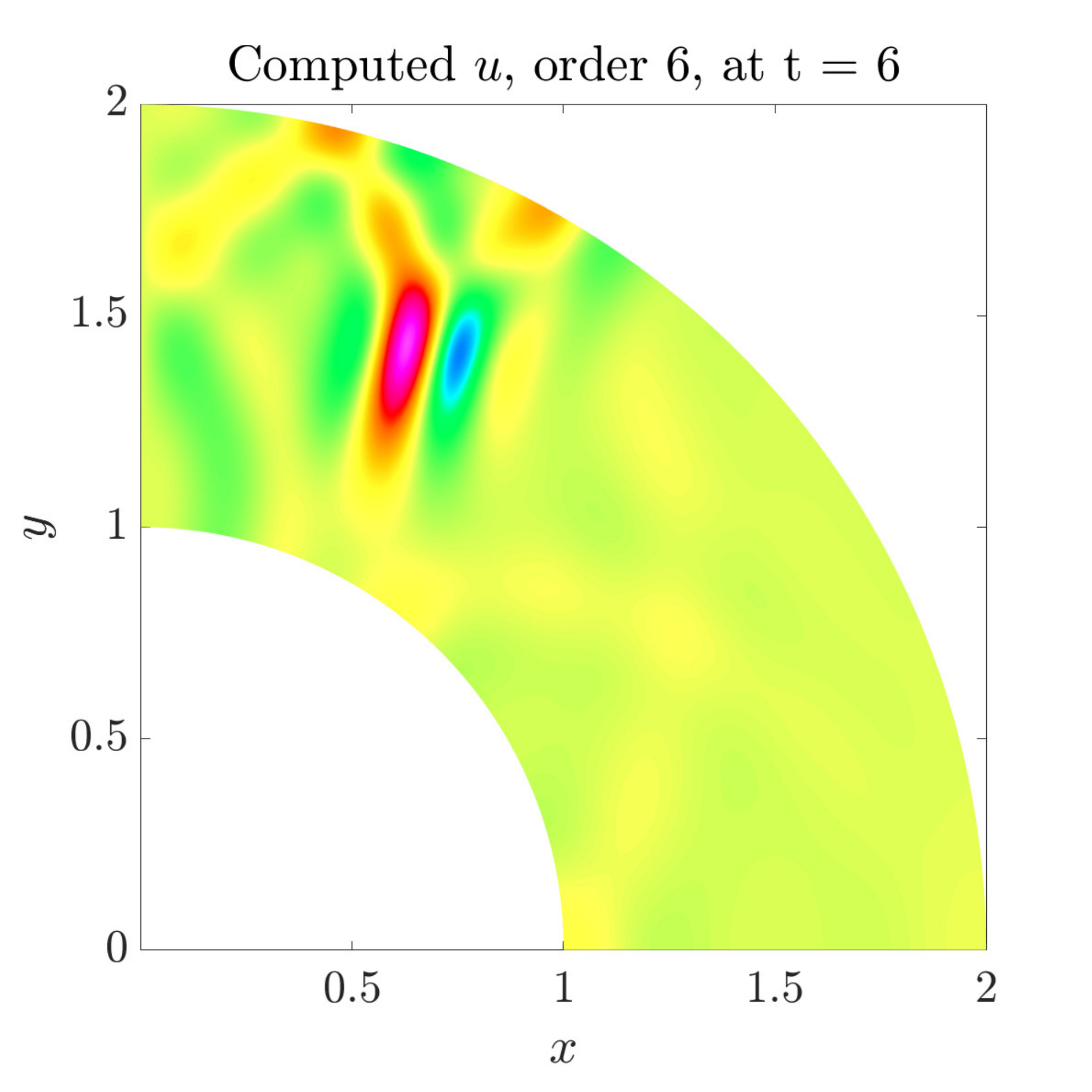}{\figWidth}};
				\draw ( 15,0.01) node[anchor=south west] {\trimfig{Resultsfig/pulseColorbar}{.665cm}};
				
			\end{tikzpicture}
		\end{center}
		\caption{Evolution of a pulse in a curved channel with closed ends.  Solutions at times $t=1,2,\ldots,6$ are computed using the sixth-order accurate scheme in~\eqref{eqn:ModifiedEqn6} for a grid with $h=1/640$, grid resolution $j=3$.}
		\label{fig:AnnularPulseEvolution}
	\end{figure}
}

%% file: tex/Example6b.tex

{
	\newcommand{\figWidth}{6.0cm}
	\newcommand{\trimfig}[2]{\trimw{#1}{#2}{0.}{0.}{.0}{.0}}  
	\newcommand{\trimfigb}[2]{\trimwb{#1}{#2}{.2}{0.1}{.2}{.15}} 
	\begin{figure}[H]
		\begin{center}
			\begin{tikzpicture}[scale=1]
				\useasboundingbox (-2,0) rectangle (10,10.5);  
				\draw ( -2.5,5.1) node[anchor=south west] {\trimfig{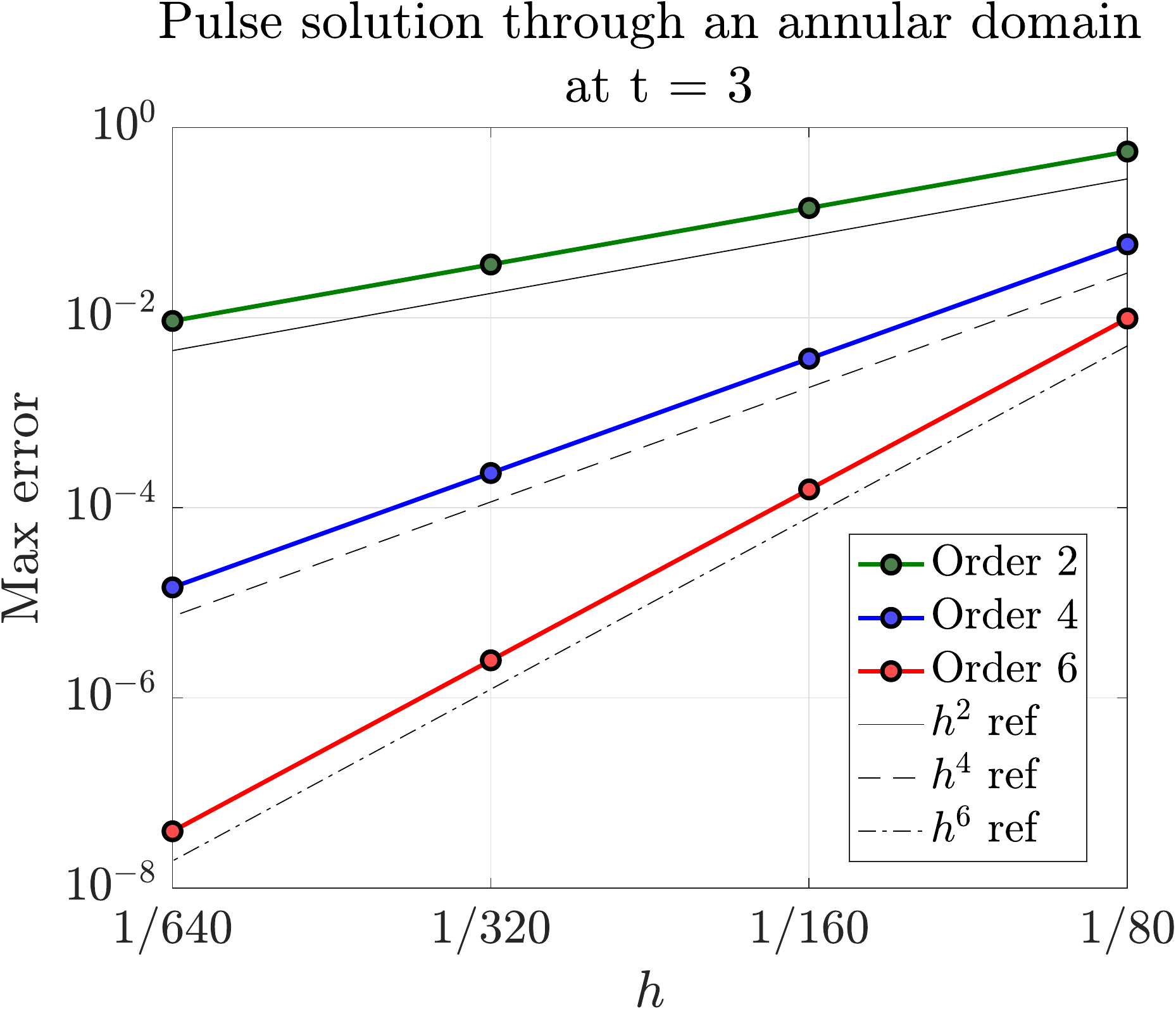}{\figWidth}};
				\draw ( 3.5,4.6) node[anchor=south west] {\trimfig{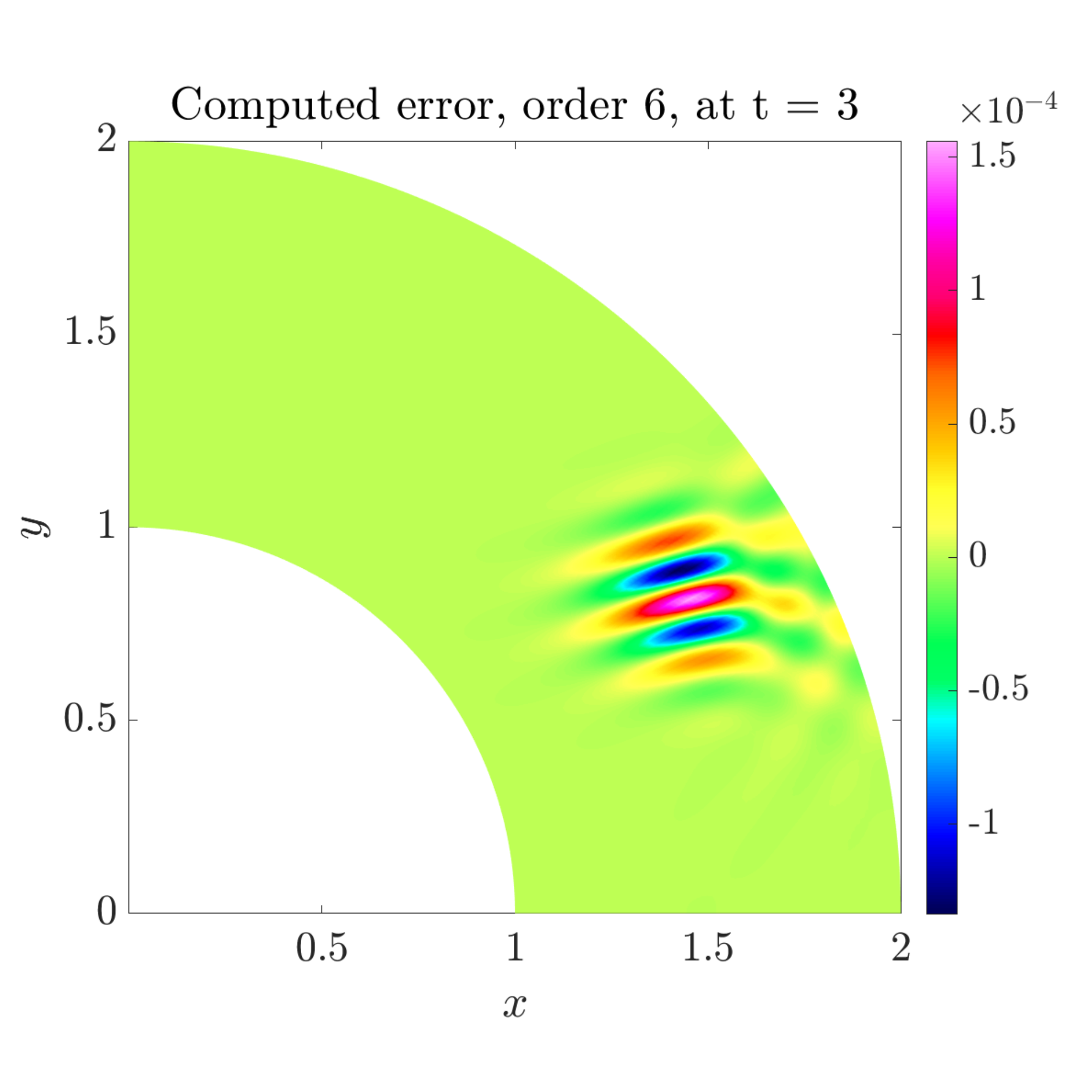}{\figWidth}};
				\draw ( -2.5,-0.6) node[anchor=south west] {\trimfig{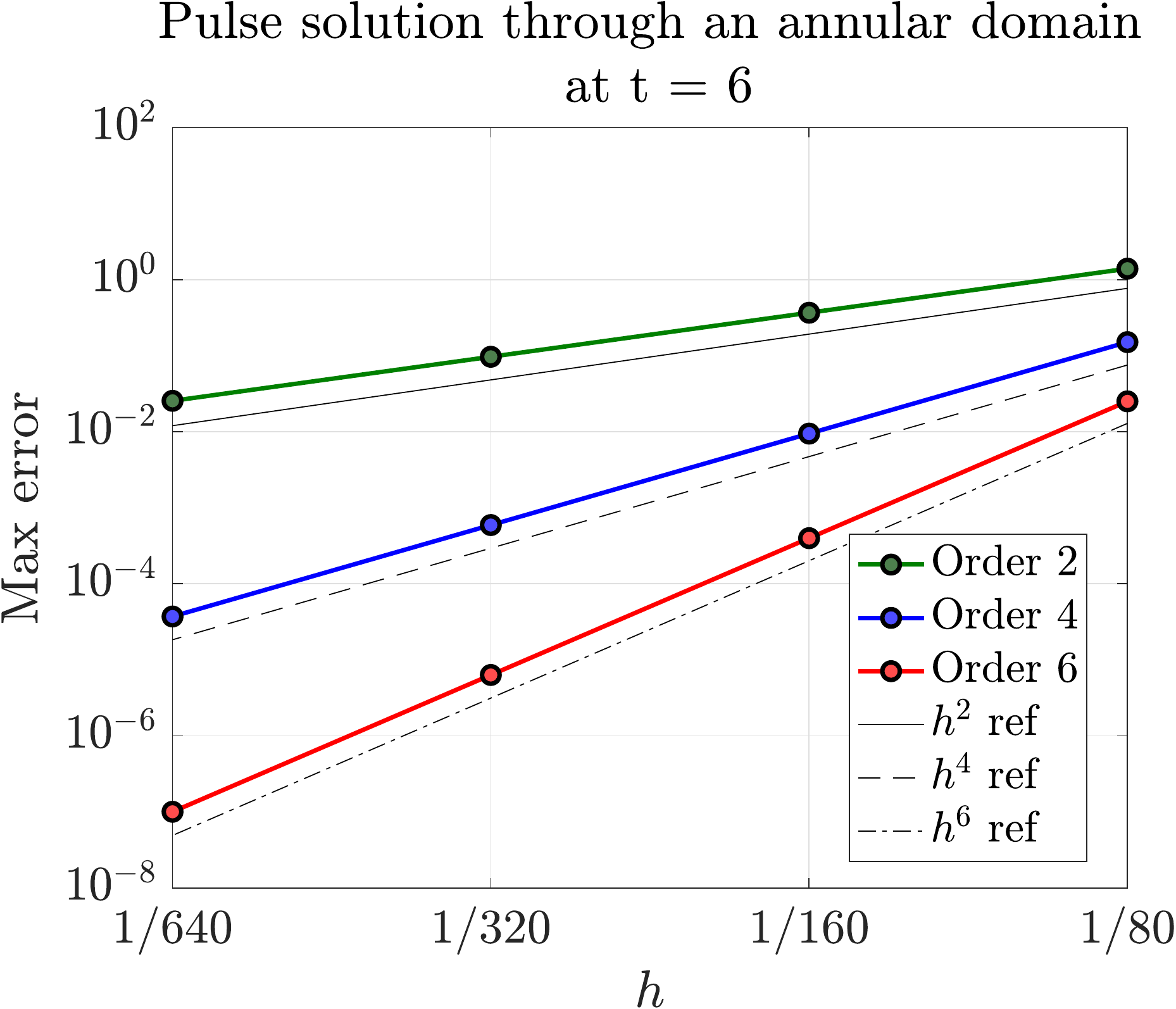}{\figWidth}};
				\draw ( 3.5,-1) node[anchor=south west] {\trimfig{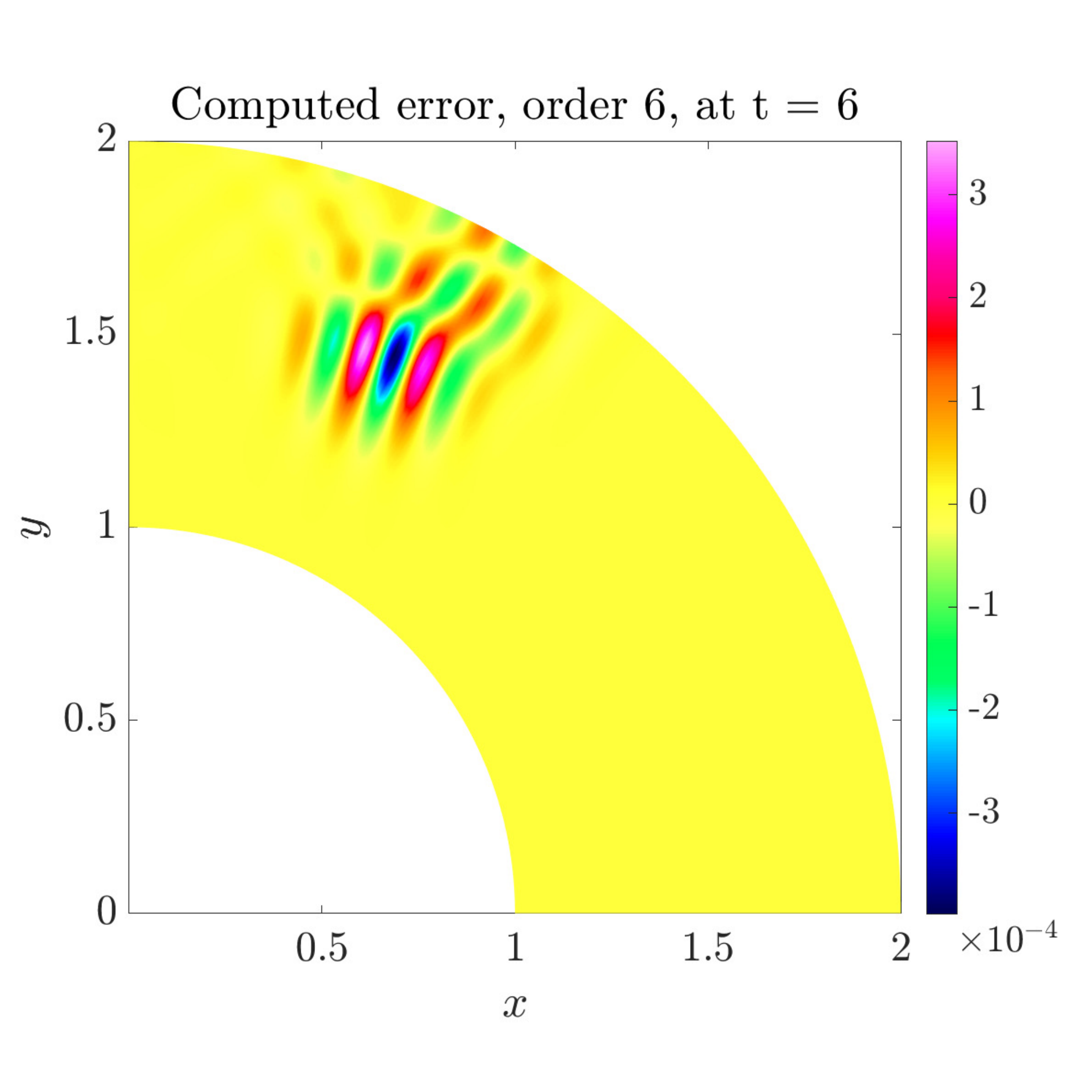}{\figWidth}};

			\end{tikzpicture}
		\end{center}
		\caption{Evolution of a pulse in a curved channel with closed ends.  Behavior of the computed maximum-norm errors at times $t=3$ (top) and $t=6$ (bottom).  Estimated errors versus $h$ (left) and errors in the sixth-order accurate solution for the grid with $h=1/160$, grid resolution $j=1$.}
		\label{fig:AnnularPulse}
	\end{figure}
}

%% file: tex/conclusion.tex
\section{Conclusion}\label{sec:Conclusion}

We have described a new approach for developing numerical approximations to boundary conditions for high-order accurate finite difference approximations.  In contrast to traditional one-sided approximations, this new local compatibility boundary condition (LCBC) approach results in centered approximations that are generally more accurate and stable than one-sided approximations.  The LCBC approach uses a local polynomial representation of the solution whose coefficients are defined using the given boundary conditions and additional compatibility boundary conditions derived from the governing equations, as well as using (known) solution data at grid points in the interior and on the boundary.  The LCBC approach uses the boundary conditions and CBCs defined at continuous level, which enables automatic construction of the local polynomial, an important consideration when developing arbitrarily high-order accurate schemes.  Values of the discrete solution at ghost points are then determined by evaluation of the local polynomial.  Algorithms have been given for computing the local polynomial as well as for forming the discrete stencil approximations for the ghost point values which can be used to efficiently assign the ghost point values. The LCBC approach at corners has also been described. In this case the boundary conditions and compatibility conditions from the two sides adjacent to the corner are used to define the local polynomial representation. 

The LCBC approach was developed for general boundary-value problems (BVPs) and initial-boundary-value problems (IBVPs) for second-order scalar PDEs, including elliptic, parabolic and hyperbolic equations.  The approach was formulated for two-dimensional domains that are described by a smooth mapping from the reference parameter space unit square.  Properties of the LCBC approximations were analyzed.  Conditions for the solvability of the LCBC equations were derived.  It was also shown that the LCBC approach results in even and odd symmetry conditions in some important special cases.  The stability of high-order approximations to the wave equation were also analyzed.

Numerical results were presented that demonstrate the accuracy and stability of the approach.  These numerical tests included manufactured solutions, problems with known solutions, and problems where the errors were estimated using a grid refinement self-convergence approach.  To test the LCBC conditions, high-order accurate schemes were presented for discretizing the BVPs and IBVPs. For time-dependent problems, formulae for choosing the maximal stable time-step were given.  Results were presented for orders of accuracy equal to two, four and six. In all cases the numerical solutions converged at close to the expected rates.

In future work we will consider extensions of the LCBC approach to BVPs and IBVPs in three dimensions, problems with interfaces, problems involving vector PDEs such as those that appear in electromagnetics or elasticity, and nonlinear problems.

%% file: tex/CornerAppendix.tex
\section{LCBC corner conditions} \label{sec:LCBCcornerGeneral}

\subsection{Neumann-Neumann Corner}
Consider the bottom-left corner, $\tilde x=(0,0)$, where two Neumann boundaries meet. The physical (primary) boundary are taken to be
\bse
\begin{align}
	\px u(\xv,t)&=g_\ell(y,t), \qquad\xv\in\partial\Omega_\ell, \\
	\py u(\xv,t)&=g_b(x,t), \qquad \xv\in\partial\Omega_b,
\end{align}
\ese
for some fixed time $t$.  We start by specifying the interpolating polynomial $\tilde u(\xv)$ at known interior data given by
\bse \label{eq:NNcornerconditions}
\begin{equation}
	\tilde{u}\bigl(\ih\dx,\jh\dy\bigr) = U_{\ih,\jh}(t),\qquad \ih = 1,\dots,p, \quad \jh = 1,\dots, p.
	\label{eq:interiorNNC}
\end{equation}
Next, we apply tangential derivatives of the primary boundary conditions and compatibility conditions given by
\begin{equation}
	\left.
	\begin{array}{l}
		\py^{\mu}\px Q^{\nu}\tilde{u}(0,0) = \py^{\mu}S_{\ell,\nu}(0,t) \smallskip\\
		\px^{\mu}\py Q^{\nu}\tilde{u}(0,0) = \px^{\mu}S_{b,\nu}(0,t)
	\end{array}\;\right\}
	\qquad \nu=0,\ldots,p-1, \quad \mu\in{\cal M}_\nu,
	\label{eq:cbcNNC}
\end{equation}
\ese
respectively, where $S_{\ell,\nu}(y,t)$ is defined in~\eqref{eq:Sdef} and
\begin{equation} \label{eq:Sbdef}
	S_{b,\nu}(x,t)\eqdef
	\begin{cases}
		-\px Q^{\nu - 1}f(x,0), & q = 0,  \\[5pt]
		\pt^{q\nu}g_b(x,t) - \px\Psi_{\nu}f(x,0,t), & q=1,2.
	\end{cases}
\end{equation}
The sets ${\cal M}_\nu$, $\nu=0,\ldots,p-1$, chosen to eliminate redundant constraints, are given by
\begin{equation}
	{\cal M}_\nu=\left\{
	\begin{array}{ll}
		0,1,2,3,\ldots,2p-1,2p, &\hbox{if $\nu=0$, with an average for $\mu=1$,}\smallskip\\
		0,\quad2,3,4,\ldots,2p-1,2p, &\hbox{if $\nu=1$, with an average for $\mu=3$,}\smallskip\\
		0,\quad2,\quad 4,5,\ldots,2p-1,2p, &\hbox{if $\nu=2$, with an average for $\mu=5$,}\smallskip\\
		\quad\quad\quad\vdots       &\quad\quad\vdots \smallskip\\
		0,\quad 2,\quad 4,\quad 6,\ldots,2p-1,2p, &\hbox{if $\nu=p-1$, with an average for $\mu=2p-1$.}
	\end{array}
	\right.
\end{equation}

\subsection{Dirichlet-Neumann Corner}
Consider the bottom-left corner, $\tilde x=(0,0)$, where a Dirichlet boundary on the left meets a Neumann boundary at the bottom. The physical (primary) boundary are taken to be
\bse
\begin{align}
	 u(\xv,t)&=g_\ell(y,t), \qquad\xv\in\partial\Omega_\ell, \\
	\py u(\xv,t)&=g_b(x,t), \qquad \xv\in\partial\Omega_b,
\end{align}
\ese
for some fixed time $t$.  We specify the interpolating polynomial $\tilde u(\xv)$ at known interior data:
\bse \label{eq:DNcornerconditions}
\begin{equation}
	\tilde{u}\bigl(\ih\dx,\jh\dy\bigr) = U_{\ih,\jh}(t),\qquad \ih = 1,\dots,p, \quad \jh = 1,\dots, p.
	\label{eq:interiorDNC}
\end{equation}
Next, we apply tangential derivatives of the primary boundary conditions and compatibility conditions given by
\begin{equation}
	\left.
	\begin{array}{l}
		\py^{\mu_1} Q^{\nu}\tilde{u}(0,0) = \py^{\mu_1}R_{\ell,\nu}(0,t) \smallskip\\
		\px^{\mu_2}\py Q^{\nu}\tilde{u}(0,0) = \px^{\mu_2}S_{b,\nu}(0,t)
	\end{array}\;\right\}
	\qquad \nu=0,\ldots,p-1, \quad (\mu_1,\mu_2)\in{\cal M}_\nu,
	\label{eq:cbcDNC1}
\end{equation}
and 
\begin{equation}
	\left.
	\begin{array}{l}
		\py^{\mu_1} Q^{\nu}\tilde{u}(0,0) = \py^{\mu_1}R_{\ell,\nu}(0,t) \smallskip\\
	\end{array}\;\right\}
	\qquad \nu=p, \quad \mu_1=0, 2, 4, \dots, 2(p-1), 2p,
	\label{eq:cbcDNC2}
\end{equation}
\ese
respectively, where $R_{\ell,\nu}(y,t)$ is defined in~\eqref{eq:Rdef} and $S_{b,\nu}(x,t)$ is defined in~\eqref{eq:Sbdef}. 

The sets ${\cal M}_\nu$, $\nu=0,\ldots,p-1$, chosen to eliminate redundant constraints, are given by
\begin{equation}
	{\cal M}_\nu=\left\{
	\begin{array}{ll}
		(\mu_1,\mu_2) \in {\cal M}_0^{(1)}\times{\cal M}_0^{(2)}, &\hbox{if $\nu=0$, with an average for $(\mu_1,\mu_2) = (1,0)$,}\smallskip\\
		(\mu_1,\mu_2) \in {\cal M}_1^{(1)}\times{\cal M}_1^{(2)}, &\hbox{if $\nu=1$, with an average for $(\mu_1,\mu_2) = (3,2)$,}\smallskip\\
		(\mu_1,\mu_2) \in {\cal M}_2^{(1)}\times{\cal M}_2^{(2)}, &\hbox{if $\nu=2$, with an average for $(\mu_1,\mu_2) = (5,4)$,}\smallskip\\
		\quad\quad\quad\vdots       &\quad\quad\vdots \smallskip\\
		(\mu_1,\mu_2) \in {\cal M}_{p-1}^{(1)}\times{\cal M}_{p-1}^{(2)}, &\hbox{if $\nu=p-1$, with an average for $(\mu_1,\mu_2) = (2p-1,2p-2)$,}
	\end{array}
	\right.
\end{equation}
where ${\cal M}_\nu^{(1)}$, ${\cal M}_\nu^{(2)}$, $\nu = 0,\dots, p-1$ are defined as 
\begin{equation}
	{\cal M}_\nu^{(1)}=\left\{
\begin{array}{ll}
	0,1,2,3,\ldots,2p-1,2p, &\hbox{if $\nu=0$,}\smallskip\\
	0,\quad2,3,4,\ldots,2p-1,2p, &\hbox{if $\nu=1$,}\smallskip\\
	0,\quad2,\quad 4,5,\ldots,2p-1,2p, &\hbox{if $\nu=2$, }\smallskip\\
	\quad\quad\quad\vdots       &\quad\quad\vdots \smallskip\\
	0,\quad 2,\quad 4,\quad 6,\ldots,2p-1,2p, &\hbox{if $\nu=p-1$,}
\end{array}
\right.
\end{equation}
and
\begin{equation}
	{\cal M}_\nu^{(2)}=\left\{
	\begin{array}{ll}
		0,1,2,3,\ldots,2p-2,2p-1,2p, &\hbox{if $\nu=0$,}\smallskip\\
		\quad 1,2,3,4,\ldots,2p-2,2p-1,2p, &\hbox{if $\nu=1$,}\smallskip\\
		\quad 1,\quad 3,4,5,\ldots,2p-2,2p-1,2p, &\hbox{if $\nu=2$,}\smallskip\\
		\quad\quad\quad\vdots       &\quad\quad\vdots \smallskip\\
		\quad 1,\quad 3,\quad 5,\ldots,2p-2,2p-1,2p, &\hbox{if $\nu=p-1$.}
	\end{array}
	\right.
\end{equation}

%% file: tex/fourthOrderWaveStability.tex
\section{Analytical proof of the stability of fourth-order LCBC approximations of the wave equation} \label{app:stabilityProof2}
In this appendix, we provide an analytical proof of the stability of the fourth-order accurate LCBC approximations of the wave equation. We have already demonstrated the stability of such approximations in section~\ref{sec:stabilityWaveEquation} via a numerical approach. In the proof of theorem~\ref{th:stabilityWaveEquation}, we arrived at the quadratic equation \eqref{eq:stabQuad} for the amplification factor $A$ which involved a parameter $b_{2}$ defined in \eqref{eq:b2Def}. We have also defined the parameters $\kHat_x$, $\kHat_y$ in \eqref{eq:kHatxyDef} and required that $\abs{b_{2}}<1$ for all values of $\kHat_x$, $\kHat_y$ for the approximations to be stable.\\

We now set $\Lambda_x = c\dt/\dx$, $\Lambda_y = c\dt/\dy$, $X = (\dx\kHat_{x})^2$ and $Y = (\dy\kHat_{y})^2$,  $X,Y\in[0,4]$. We write $b_2$ as 
\begin{equation*}
	b_2 =  1 - \frac{1}{2}\left[\Lambda_x^2X\left(1 + \frac{X}{12}\right)+\Lambda_y^2Y\left(1 + \frac{Y}{12}\right) \right]+ \frac{1}{24}\left(\Lambda_x^2X + \Lambda_y^2Y \right)^2, 
\end{equation*}
and define the variable $z = \Lambda_x^2 + \Lambda_y^2$. \\

From the numerical results in figure~\ref{fig:stabRegion}, we observe that 
\begin{center}
\textit{Given $\dt$ such that $z<1$, $\abs{b_2}<1$ for all $X,Y\in[0,4]$.}
\end{center}
\begin{proof}
Assuming $z<1$, we aim to show $\abs{b_2}<1$ or $b_2^2<1$. For this, we write $b_2^2 -1 = F(F + 2)$ where
\begin{align}
	F(X,Y) =& -\frac{1}{2}\left[\Lambda_x^2X\left(1+\f{X}{12}\right) + \Lambda_y^2Y\left(1+\f{Y}{12}\right)\right] + \f{1}{24}\left(\Lambda_x^2X + \Lambda_y^2Y\right)^2,\nonumber \\
	=&\underbrace{\f{1}{24}\Lambda_x^2(\Lambda_x^2-1)}_{\Ac}X^2 + \underbrace{\f{1}{12}\Lambda_x^2\Lambda_{y}^2}_{\Bc}XY + \underbrace{\f{1}{24}\Lambda_y^2(\Lambda_y^2-1)}_{\Cc}Y^2 -\f{1}{2}\Lambda_{x}^2X - \f{1}{2}\Lambda_y^2Y. \label{eq:F}
\end{align}
For stability we require $-2<F(X,Y)<0$ for all values of $X$, $Y$. The function $F(X,Y)$ in \eqref{eq:F} is a bi-variate quadratic function in $X$ and $Y$. We find the discriminant of $F$ to be
\begin{align*}
	\Bc^2 - 4\Ac\Cc & = \f{1}{144}\Lambda_x^4\Lambda_y^4- \f{1}{144}\Lambda_x^2\Lambda_y^2(\Lambda_x^2 - 1)(\Lambda_y^2-1), \\
	&=  \f{1}{144}\Lambda_x^2\Lambda_y^2\left[\Lambda_x^2\Lambda_y^2 - \left(\Lambda_x^2\Lambda_y^2 -\Lambda_x^2 - \Lambda_y^2 + 1\right)\right], \\
	& = \f{1}{144}\Lambda_x^2\Lambda_y^2(z - 1)<0,
\end{align*}
which shows that $F$ is an elliptic paraboloid. Also, with $\Lambda_x^2,\Lambda_{y}^2<z<1$, we find that $\Ac,\Bc<1$ and $F$ is concave downwards. The maximum of $F$ occurs at $X = Y = 0$ where $F(0,0) = 0$. The minimum of $F$ occurs at $X = Y = 4$ where $$f(z)\defeq F(4,4) = -\f{8}{3}z + \f{2}{3}z^2.$$
The function $f(z)$ is a monotonically decreasing function with increasing $z\in(0,1)$.  We find $\lim_{z\rightarrow 1}f(z) = -2$ which is within the bound needed for stability. Hence, when $z<1$, $b_2^2<1$ and the fourth-order accurate LCBC approximation of the wave equation is stable. \\
\noindent{\Huge $\square$}
\end{proof}

%% file: tex/timeStepRestrictions.tex
\section{Time-step restrictions} \label{sec:timeStepRestrictions}

In this section we provide details of the determination of the time-step for the schemes
used in Section~\ref{sec:NumericalResults}.
For variable coefficient problems we freeze coefficients and consider the constant coefficient operator
\ba 
   Q  = c_{11} \p_x^2  + 2c_{12} \p_x\p_y + c_{22} \p_y^2 + c_1 \p_x + c_2 \p_y + c_0, \label{eq:Qcc}
\ea 
where $c_{11}, c_{12}, c_{22}, c_1, c_2, c_0$ are constants.
We then perform a von-Neumann
stability analysis to arrive at a time-step restriction for the frozen coefficients.
The final time-step is chosen to enforce stability over the range of possible coefficients
in the problem at hand.

\subsection{Parabolic time-step restrictions for method of lines schemes} \label{sec:parabolicStepRestrictions}

We consider the parabolic problem
\begin{equation}
  u_{t} = Qu, \qquad \xv \in \Omega = [0,2\pi]^2, \label{eq:parabolicCC}
\end{equation}
with $Q$ given by~\eqref{eq:Qcc}
together with periodic boundary conditions in both $x,y$ directions and an initial condition. 
We consider a method of lines approximation.
Let
\ba 
  \xv_\jv=(j_1\dx,j_2\dy),  \qquad j_1=0,1,\ldots,N_x,  \quad j_2=0,1,\ldots,N_y,
  \label{eq:TSgrid}
\ea 
denote the grid points for $\Omega$ with $N_x+1$ and $N_y+1$ points in the $x$ and $y$ directions, 
where $\dx=2\pi/N_x$ and $\dy=2\pi/N_y$.
Let $U_{\jv}(t)\approx u(\xv_{\jv}, t)$ denote the semi-discrete grid function. 
 
\paragraph{Second-order accurate discretization}
We discretize~\eqref{eq:parabolicCC} to second-order accuracy in space
\begin{equation}\label{eqn:DiscreteParabolic2}
\frac{d}{dt}U_{\jv}(t) = c_{11}\dpmx U_{\jv} + 2c_{12}\dzx\dzy U_{\jv} + c_{22}\dpmy U_{\jv} + c_{1}\dzx U_{\jv} + c_{2}\dzy U_{\jv} + c_0U_{\jv} .
\end{equation}
We perform a von-Neumann analysis to find the time-step restriction.
The semi-discrete solution is expanded in a discrete Fourier series in space. We thus make the ansatz
$$
   U_{\jv}(t) = \Vhat_{\kv}(t)e^{i\kv\cdot\xv_{\jv}},  \qquad k_x = - \Nx/2, \dots, (\Nx/2) -1, \ k_y = - \Ny/2, \dots, (\Ny/2)-1, 
$$
where $\kv = (k_x,k_y)$ and where we have assumed $N_x$ and $N_y$ are even for convenience.
Substitute this ansatz into~\eqref{eqn:DiscreteParabolic2} and simplify to arrive at the standard test equation 
\ba
  & \frac{d\Vhat_{\kv}}{dt} = \lambda_{\kv}^{(2)} \Vhat_{\kv}, 
\ea
where the time-stepping eigenvalue is 
\bas 
  & \lambda_{\kv}^{(2)} 
   = \alpha + i \beta , \\
  & \alpha \eqdef -c_{11}\, \kHat_x^2 - 2 c_{12} \,\kHatzx \, \kHatzy  - c_{22} \, \kHat_y^2 +c_0 , \\
  & \beta \eqdef c_1 \, \kHatzx + c_2 \,\kHatzy  ,
\eas
and where
\bat
  & \kHat_x \eqdef  \f{\sin(\xi_x/2)}{\dx/2},  \qquad&& \kHat_y \eqdef  \f{\sin(\xi_y/2)}{\dy/2}, \label{eq:kHatDef} \\
  & \kHatzx \eqdef  \f{\sin(\xi_x)}{\dx},      \qquad&& \kHatzy \eqdef  \f{\sin(\xi_y)}{\dy}, \label{eq:kHatzDef}\\
  & \xi_x \eqdef k_x\dx,                       \qquad&& \xi_y \eqdef k_y\dy.
\eat
Note that we have the bounds,
\bat
   &  -\pi \le \xi_x \le \pi,   \qquad&& -\pi \le \xi_y \le \pi, \\
   &  |\kHat_x| \le \f{2}{\dx}, \qquad&& |\kHat_y| \le \f{2}{\dy}, \\
   &  |\kHatzx| \le \f{1}{\dx}, \qquad&& |\kHatzy| \le \f{1}{\dy}.
\eat
We now un-freeze the coefficients to obtain an estimate of the worst case time-stepping eigenvalue 
$\lambdaMax = \alphaMin + i \betaMax$ as 
\begin{align*}
  \alphaMin & = \min_{\xv\in\Omega^P} \left[ -\dfrac{4\abs{c_{11}(\xv)}}{\dx^2} - \dfrac{2\abs{c_{12}(\xv)}}{\dx\dy} - \dfrac{4\abs{c_{22}(\xv)}}{\dy^2} + c_0(\xv) \right],\\
  \betaMax & = \max_{\xv\in\Omega^P} \left[ \frac{\abs{c_1(\xv)}}{\dx} + \frac{\abs{c_2(\xv)}}{\dy} \right],
\end{align*}
where $\Omega^P$ is the physical domain under consideration.
We approximate the region of absolute stability of 
the time-stepping scheme as an ellipse in the complex plane for $z=x+iy$, 
(for schemes not stable on the imaginary axis such as forward-Euler 
this assumes we are not too close to the imaginary axis)
\ba
    \left( \f{x}{\alpha_0}\right)^2 +  \left( \f{y}{\beta_0}\right)^2 \le 1 ,
\ea
and ensure that $z=\lambdaMax\dt$ lies within the approximate stability region.
This leads to the condition,
\ba
    \left( \f{\alphaMin\dt}{\alpha_0}\right)^2 +  \left( \f{\betaMax\dt}{\beta_0}\right)^2 \le 1 ,
\ea
which can be re-arranged to give an inequality for $\dt$.

For forward-Euler we approximate the region of absolute stability 
using ellipse parameters $\alpha_0=2$ and $\beta_0=1$.

\paragraph{Fourth-order accurate discretization}
We now use fourth-order accurate centered differences 
in~\eqref{eqn:DiscreteParabolic2}. The time stepping eigenvalue is 
\begin{align*}
\lambda_{\kv}^{(4)} =& 
  - c_{11} \, \kHat_x^2\left[1 + \frac{\dx^2}{12}\kHat_x^2 \right]
  -  2c_{12}\, \kHatzx\left[1 + \frac{\dx^2}{6}\kHat_x^2\right]
            \, \kHatzy\left[1 + \frac{\dy^2}{6}\kHat_y^2\right]\\
  & -c_{22} \, \kHat_y^2\left[1 +\frac{\dy^2}{12}\kHat_y^2\right]
  + c_1 \, i\kHatzx\left[1 + \frac{\dx^2}{6}\kHat_x^2\right]
  + c_2 \, i\kHatzy\left[1 + \frac{\dy^2}{6}\kHat_y^2\right] + c_0.
\end{align*}
We shall use the bounds
\ba  
   &  \kHat_x^2\left[1 + \frac{\dx^2}{12}\kHat_x^2 \right] \le \frac{16}{3\dx^2}, \\
   &  \abs{\kHatzx} \, \left[1 + \frac{\dx^2}{6}\kHat_x^2\right] \le \frac{3}{2\dx},
\ea 
and similar ones for the $y$ variable to obtain an estimate for 
the worst case time-stepping eigenvalue with real and imaginary parts given by 
\begin{align*}
    \alphaMin & = \min_{\xv\in\Omega^P} \left[
       -\frac{16\abs{c_{11}(\xv)}}{3\dx^2} - \frac{9\abs{c_{12}(\xv)}}{2\dx\dy} - \frac{16\abs{c_{22}(\xv)}}{3\dy^2} + c_{0}(\xv) \right],\\
  \betaMax & = \max_{\xv\in\Omega^P} \left[ \frac{3\abs{c_1(\xv)}}{2\dx} + \frac{3\abs{c_2(\xv)}}{2\dy} \right].
\end{align*}
\paragraph{Sixth-order accurate discretization}
For sixth-order accurate centered differences 
in~\eqref{eqn:DiscreteParabolic2}, the time stepping eigenvalue is 
\begin{align*}
	\lambda_{\kv}^{(6)} =
	& - c_{11} \, \kHat_x^2\left[1 + \frac{\dx^2}{12}\kHat_x^2 + \frac{\dx^4}{90}\kHat_{x}^4\right]
   -c_{22} \, \kHat_y^2\left[1 +\frac{\dy^2}{12}\kHat_y^2+\frac{\dy^4}{90}\kHat_y^4\right]\\
	&-  2c_{12}\, \kHatzx\left[1 + \frac{\dx^2}{6}\kHat_x^2 + \frac{\dx^4}{30}\kHat_x^4\right]
	\, \kHatzy\left[1 + \frac{\dy^2}{6}\kHat_y^2+ \frac{\dy^4}{30}\kHat_y^4\right]\\
&+ c_1 \, i\kHatzx\left[1 + \frac{\dx^2}{6}\kHat_x^2 + \frac{\dx^4}{30}\kHat_x^4\right]
	+ c_2 \, i\kHatzy\left[1 + \frac{\dy^2}{6}\kHat_y^2+ \frac{\dy^4}{30}\kHat_y^4\right] + c_0.
\end{align*}
We use the bounds
\ba  
&  \kHat_x^2\left[1 + \frac{\dx^2}{12}\kHat_x^2 + \frac{\dx^4}{90}\kHat_x^4 \right] \le \frac{272}{45\dx^2}, \\
&  \abs{\kHatzx} \, \left[1 + \frac{\dx^2}{6}\kHat_x^2+ \frac{\dx^4}{30}\kHat_x^4\right] \le \frac{1199}{756\dx},
\ea 
and similar ones for the $y$ variable to estimate
the worst case time-stepping eigenvalue with real and imaginary parts as follows
\begin{align*}
	\alphaMin & = \min_{\xv\in\Omega^P} \left[
	-\frac{272\abs{c_{11}(\xv)}}{45\dx^2} - \frac{1313\abs{c_{12}(\xv)}}{261\dx\dy} - \frac{272\abs{c_{22}(\xv)}}{45\dy^2} + c_{0}(\xv) \right],\\
	\betaMax & = \max_{\xv\in\Omega^P} \left[ \frac{1199\abs{c_1(\xv)}}{756\dx} + \frac{1199\abs{c_2(\xv)}}{756\dy} \right].
\end{align*}
\subsection{Time-step restriction for hyperbolic problems} \label{sec:hyperbolicStepRestrictions}

We consider the hyperbolic problem 
\begin{equation}\label{eq:hyperbolicCC}
  u_{tt} = Qu, \qquad \xv \in [0,2\pi]^2,
\end{equation}
with $Q$ given in~\eqref{eq:Qcc} together with periodic boundary conditions in both $x,y$ directions and initial conditions. We set $U_{\jv}^n\approx u(\xv_{\jv},t^n)$ where $\xv_{\jv}$ is defined in~\eqref{eq:TSgrid}.

\paragraph{Second-order accurate discretization}
We descretize \eqref{eq:hyperbolicCC} to second-order accuracy in space and time
\begin{equation}\label{eq:hyperbolicDS}
  \dpmt U_{\bm{j}}^{n} = c_{11}\dpmx U_{\bm{j}}^{n} + 2c_{12}\dzx\dzy U_{\bm{j}}^{n} + c_{22}\dpmy U_{\bm{j}}^{n} + c_{1}\dzx U_{\bm{j}}^{n} + c_{2}\dzy U_{\bm{j}}^{n} + c_0U_{\bm{j}}^n, \quad \xv_{\bm{j}}\in\Omega_h.
\end{equation}

We perform a Von-Neumann analysis to find the time-step restriction where we make the ansatz
$$
U_{\jv}^n = A^ne^{i\kv\cdot\xv_{\jv}},  \qquad k_x = - \Nx/2, \dots, (\Nx/2) -1, \ k_y = - \Ny/2, \dots, (\Ny/2)-1, 
$$
with $\kv = (k_x,k_y)$ and $\Nx$,$\Ny$ assumed even. We plug the ansatz in \eqref{eq:hyperbolicDS} and simplify to arrive at a quadratic equation for $A$
\begin{equation}\label{eq:Aquad}
  A^2 - 2bA + 1 = 0,
\end{equation}
where
\begin{equation}
b= 1 + \dfrac{\dt^2}{2}\hat{Q}_h,
\end{equation}
and
\begin{equation}
  \hat{Q}_{h} = -c_{11}\kHat_x^2-2c_{12}\kHatzx\kHatzy- c_{22}\kHat_y^2+ic_{1}\kHatzx + ic_2\kHatzy + c_0.
\end{equation}
The parameters $\kHat_{x},\kHat_y,\kHatzx,\kHatzy$ are defined in~\eqref{eq:kHatDef} and~\eqref{eq:kHatzDef}.\\

For stability, we wish to find a bound on $\dt$ such that the roots of \eqref{eq:Aquad}, $A_{\pm}$,  satisfy 
$\abs{A_{+}},\abs{A_{-}} \le 1$, and that there be no double roots with $|A|=1$, to avoid
algebraic growth in time. From equation \eqref{eq:Aquad}, we see that $A_{+}A_{-} = 1$. Therefore, $A_+ = e^{i\theta}$ and $A_- = e^{-i\theta}$, for $\theta\in\mathbb{R}$. Furthermore, $b = (A_+ + A_-)/2 = \cos\theta$. 
Since there are double roots when $b=\pm1$, we require $\dt$ such that $b\in\mathbb{R}$ and $\abs{b}<1$.  \\

For this problem, $b$ is complex if $c_1 \ne 0$ or $c_2 \ne 0$ and then we cannot satisfy $|A|\le 1$. 
In this case some dissipation should therefore be added to the scheme so that there is no growth in time.
We note that we could generalize our definition of stability to allow bounded growth in time
in which case the scheme could formally be stable using this extended definition of stability. 
This is because as we refine the mesh, and assuming $\dt \le ({\rm const.}) \min \{ \dx, \dy\}$ then 
\begin{equation}
b\sim \underbrace{1 - \f{\dt^2}{2}\left(c_{11}\kHat_x^2+2c_{12}\kHatzx\kHatzy+ c_{22}\kHat_y^2\right)}_{\tilde{b}} 
    + \Oc(\dt) , \qquad \text{as } \dt\rightarrow 0,  \label{eq:waveGenStabilityB}
\end{equation} 
since $\dt^2 |\kHatzx| \le ({\rm const.}) \dt$  and $\dt^2 |\kHatzy| \le ({\rm const.}) \dt$.
The contribution of the complex part of the symbol to $b$ in~\eqref{eq:waveGenStabilityB},
 which comes from the lower order terms, goes to zero like $\dt$ as $\dt\rightarrow 0$.
 Therefore, provided $|\tilde{b}|<1$, we have  $|A_\pm| = 1 + \Oc(\dt)$
and the scheme will thus be stable, with bounded growth.
Even in this case it is still recommended that dissipation be added to the scheme.

Consider then
\begin{align*}
\tilde{b}^2 =& 1 - \dt^2\left(c_{11}\kHat_{x}^2 + 2c_{12}\kHatzx\kHatzy + c_{22}\kHat_{y}^2 \right)\\
&+\frac{\dt^4}{4}\left(c_{11}\kHat_{x}^2 + 2c_{12}\kHatzx\kHatzy + c_{22}\kHat_{y}^2\right)^2 <1.
\end{align*} 
We look for a time-step $\dt$ such that
\begin{equation}
\frac{\dt^2}{4}\left(c_{11}\kHat_{x}^2 + 2c_{12}\kHatzx\kHatzy + c_{22}\kHat_{y}^2\right) <1.
\end{equation}

With the necessary condition that $c_{11}c_{22}>c_{12}^2$, we find that
\begin{equation}\label{eq:hyperbolicBound}
	\max_{-\pi\leq\xi_{x},\xi_y<\pi}\left[c_{11}\kHat_{x}^2 +2c_{12}\kHatzx\kHatzy + c_{22}\kHat_{y}^2\right] = 4\left(\frac{c_{11}}{\dx^2} + \frac{c_{22}}{\dy^2}\right),
\end{equation}
occurs when $(\xi_{x},\xi_{y}) = (-\pi,-\pi)$. 
Hence, 
\begin{equation}
	\dt < \dfrac{1}{\sqrt{\frac{c_{11}}{\dx^2} + \frac{c_{22}}{\dy^2}}}.
\end{equation}
We obtain the bound in \eqref{eq:hyperbolicBound} by considering the maximum of the function 
\begin{equation}
\chi(\xi_{x},\xi_{y}) = c_{11}\kHat_{x}^2(\xi_x) + 2c_{12}\kHatzx(\xi_x)\kHatzy(\xi_y) + c_{22}\kHat_y^2(\xi_y) . 
\end{equation}
We find the critical points of $\chi$ using the equations
\begin{align}
  &\chi_{\xi_x} = \frac{c_{11}}{\dx^2}\sin\xi_x + \frac{c_{12}}{\dx\dy}\cos\xi_x\sin\xi_y = 0,\label{eqn:Heq1} \\
  &\chi_{\xi_y} = \frac{c_{22}}{\dy^2}\sin\xi_y + \frac{c_{12}}{\dx\dy}\sin\xi_x\cos\xi_y = 0.\label{eqn:Heq2}
\end{align}
Seeing that $\cos(\xi_y)\neq0$ (otherwise we arrive at a contradiction), equation \eqref{eqn:Heq2} gives
\begin{equation}\label{eqn:sinxiExp}
  \sin\xi_x = - \frac{\dx}{\dy}\frac{c_{22}}{c_{12}}\tan\xi_y.
\end{equation}
We use~\eqref{eqn:sinxiExp} to rewrite \eqref{eqn:Heq1} in terms of $\xi_y$,
\begin{equation}
  \sin\xi_y\left(-\dfrac{c_{11}c_{22}}{c_{12}}\frac{1}{\cos\xi_y}\pm c_{12}\sqrt{1 - \dfrac{\dx^2}{\dy^2}\frac{c_{22}^2}{c_{12}^2}\tan^2\xi_y}\right)  = 0. 
\end{equation}

If we suppose that 
\begin{equation*}
  \pm c_{12}\sqrt{1 - \dfrac{\dx^2}{\dy^2}\frac{c_{22}^2}{c_{12}^2}\tan^2\xi_y} = \frac{c_{11}c_{22}}{c_{12}}\frac{1}{\cos\xi_y},
\end{equation*}
square both sides and simplify to obtain
\begin{align*}
  \cos^2\xi_y - \dfrac{\dx^2}{\dy^2}\dfrac{c_{22}^2}{c_{12}^2}\sin^2\xi_y &= \dfrac{c_{11}^2c_{22}^2}{c_{12}^4}, \\
  1 - \sin^2\xi_y - \dfrac{\dx^2}{\dy^2}\dfrac{c_{22}^2}{c_{12}^2}\sin^2\xi_y &= \dfrac{c_{11}^2c_{22}^2}{c_{12}^4}, \\
  \sin^2\xi_y\left(-c_{12}^4 - \frac{\dx^2}{\dy^2}c_{22}^2c_{12}^2\right) &= c_{11}^2c_{22}^2 - c_{12}^4,
\end{align*}
we find that $c_{11}^2c_{22}^2 - c_{12}^4<0 \implies c_{11}c_{22}<c_{12}^2$ which violates the necessary condition $c_{11}c_{22}>c_{12}^2$. 

Therefore, the function  $\chi(\xi_x,\xi_y)$ has critical points when $\sin(\xi_y) = \sin(\xi_x) = 0$, or at $(-\pi,-\pi)$,$(-\pi,0)$, $(0,\pi)$, $(0,0)$ and is maximum at $(-\pi,-\pi)$.